\documentclass[reqno,11pt]{amsart}

\usepackage{amsmath}
\usepackage{amsthm}
\usepackage{amssymb}
\usepackage{amscd}
\usepackage{xypic}
\usepackage{mathabx}
\usepackage{graphicx}
\usepackage{palatino}
\usepackage{bbm}
\usepackage{extarrows}

\usepackage{pdflscape}
\usepackage{stmaryrd}

\usepackage{multirow}

\usepackage[plainpages=false,pagebackref,colorlinks,hyperindex,pdfpagemode=None,bookmarksopen,linkcolor=red,citecolor=blue,urlcolor=blue]{hyperref}

\theoremstyle{plain}
\newtheorem{theorem}[subsubsection]{Theorem}

\newtheorem*{theorem*}{Theorem}
\newtheorem{proposition}[subsubsection]{Proposition}
\newtheorem*{proposition*}{Proposition}
\newtheorem{lemma}[subsubsection]{Lemma}
\newtheorem*{lemma*}{Lemma}

\newtheorem{corollary}[subsubsection]{Corollary}
\newtheorem*{corollary*}{Corollary}

\theoremstyle{definition}
\newtheorem{definition}[subsubsection]{Definition}

\theoremstyle{remark}
\newtheorem{remark}[subsubsection]{Remark}
\newtheorem{remarks}[subsubsection]{Remarks}

\newtheorem{example}[subsubsection]{Example}

\newcommand{\comment}[1] {  }

\DeclareFontFamily{OT1}{rsfs}{}
\DeclareFontShape{OT1}{rsfs}{n}{it}{<-> rsfs10}{}
\DeclareMathAlphabet{\mathscr}{OT1}{rsfs}{n}{it}

\newcommand{\Ad}{\mathrm{Ad}}
\newcommand{\Res}{\mathrm{Res}}
\newcommand{\CC}{\mathbb{C}}

\newcommand{\ZZ}{\mathbb{Z}}
\newcommand{\RR}{\mathbb{R}}

\newcommand{\C}{\mathfrak{C}}
\newcommand{\B}{\mathcal{B}}
\newcommand{\R}{\mathcal{R}}

\newcommand{\ad}{{\operatorname{ad}}}
\renewcommand{\P}{\mathbbm{P}}

\newcommand{\Ind}{\operatorname{Ind}}
\newcommand{\Hom}{\operatorname{Hom}}

\newcommand{\Aut}{{\operatorname{Aut}}}
\newcommand{\Planch}{{\operatorname{Planch}}}
\newcommand{\Haar}{{\operatorname{Haar}}}

\newcommand{\Gm}{\mathbbm{G}_m}
\newcommand{\Ga}{\mathbbm{G}_a}

\newcommand{\GL}{\operatorname{GL}}
\newcommand{\Spin}{\operatorname{Spin}}

\newcommand{\Sym}{\operatorname{Sym}}
\newcommand{\PGL}{\operatorname{PGL}}
\newcommand{\SL}{\operatorname{SL}}
\newcommand{\Sp}{\operatorname{Sp}}

\newcommand{\SO}{{\operatorname{SO}}}

\newcommand{\Spec}{\operatorname{Spec}}

\newcommand{\diag}{{\operatorname{diag}}}

\newcommand{\pr}{{\operatorname{pr}}}

\newcommand{\reg}{{\operatorname{reg}}}
\newcommand{\gr}{{\operatorname{gr}}}

\newcommand{\rk}{{\operatorname{rk}}}

\newcommand{\Std}{{\operatorname{Std}}}

\newcommand{\gl}{{\operatorname{gl}}}
\newcommand{\codim}{{\operatorname{codim}}}

\newcommand{\inv}{{\operatorname{inv}}}

\newcommand{\temp}{{\operatorname{temp}}}

\newcommand{\RTF}{{\operatorname{RTF}}}

\renewcommand{\c}{\mathfrak{c}}
\renewcommand{\a}{\mathfrak{a}}
\renewcommand{\u}{\mathfrak{u}}
\newcommand{\g}{\mathfrak{g}}
\newcommand{\h}{\mathfrak{h}}
\newcommand{\Meas}{{\operatorname{Meas}}}

\newcommand{\sing}{{\operatorname{sing}}}

\newcommand{\rs}{{\operatorname{rs}}}

\newcommand{\LG}{{^LG}}

\renewcommand{\top}{{\operatorname{top}}}
\newcommand{\nilp}{{\operatorname{nilp}}}

\newcommand{\Lie}{\operatorname{Lie}}

\begin{document}
\title[Functorial transfer in rank one]{Functorial transfer between relative trace formulas in rank one.}
\author{Yiannis Sakellaridis}
\email{sakellar@jhu.edu}
\address{Department of Mathematics, Johns Hopkins University, Baltimore, MD 21218, USA.}


\begin{abstract}
According to the Langlands functoriality conjecture, broadened to the setting of spherical varieties (of which reductive groups are special cases), a map between $L$-groups of spherical varieties should give rise to a functorial transfer of their local and automorphic spectra. The ``Beyond Endoscopy'' proposal predicts that this transfer will be realized as a comparison between limiting forms of the (relative) trace formulas of these spaces. 

In this paper we establish the local transfer for the identity map between $L$-groups, for spherical affine homogeneous spaces $X=H\backslash G$ whose dual group is $\SL_2$ or $\PGL_2$ (with $G$ and $H$ split). More precisely, we construct a transfer operator between orbital integrals for the $(X\times X)/G$-relative trace formula, and orbital integrals for the Kuznetsov formula of $\PGL_2$ or $\SL_2$. Besides the $L$-group, another invariant attached to $X$ is a certain $L$-value, and the space of test measures for the Kuznetsov formula is enlarged, to accommodate the given $L$-value.

The fundamental lemma for this transfer operator is proven in a forthcoming paper of Johnstone and Krishna. The transfer operator is given explicitly in terms of Fourier convolutions, making it suitable for a global comparison of trace formulas by the Poisson summation formula, hence for a uniform proof, in rank one, of the relations between periods of automorphic forms and special values of $L$-functions.

\end{abstract}

\maketitle
\tableofcontents{}

\section{Introduction}

\subsection{Relative functoriality}\label{ssrelfunctoriality}

According to the Relative Langlands Program, the local and automorphic spectra of a spherical $G$-variety $X$ should be determined by its $L$-group $\LG_X$, which comes equipped with a distinguished morphism 
\begin{equation}\label{Lgroups}
 \LG_X\times\SL_2\to \LG,
\end{equation}
cf.\ \cite{GN, SV, KnSch}.

Roughly speaking, this means, locally, that the Plancherel formula for $L^2(X(F))$ (where $F$ is a local field)  should read:
$$ \left<\Phi_1,\overline{\Phi_2}\right>_{L^2(X)} = \int_{\varphi} J_\varphi^\Planch(\Phi_1\otimes\Phi_2) \nu_X(\varphi),$$
where the integral is over the space of Langlands parameters into $\LG_X$,
$\nu_X$ is the standard measure \cite[\S 17.3]{SV} on this space, and $J_\varphi^\Planch$ is a \emph{relative character}
$$ J_{\varphi} :\mathcal S(X\times X)\to \Pi_\varphi \otimes \widetilde{\Pi_{\varphi}} \to \CC.$$
Here, $\Pi_{\varphi}$ is the sum of irreducible representations in the Arthur packet associated to the composition of $\varphi$ with the canonical map \eqref{Lgroups}. 

Globally, an analogous decomposition (in terms of ``global Langlands parameters'') should hold for the spectral side of the \emph{relative trace formula} of $X$ --- more precisely, for its stable part ---, a distribution on the adelic points of the quotient $X\times X/G$ (with $G$ acting diagonally), whose spectral decomposition should read, roughly:
$$ \RTF_X(\Phi_1\otimes \Phi_2) = \int_{\varphi} J_\varphi^\gl(\Phi_1\otimes\Phi_2) \mu_X(\varphi).$$
Moreover, the global distributions $J_\varphi^\gl$, which can be expressed in terms of squares of periods of automorphic forms, should (under some assumptions on $X$) be equal to Euler products of the local distributions $J_{\varphi_v}^\Planch$, establishing a link between periods of automorphic forms and special values of $L$-functions; this is the generalized Ichino--Ikeda conjecture proposed in \cite[\S 17.4]{SV}. 

We currently have no general tools to address these very general, and uniform, conjectures. In this paper, I will propose a uniform approach which works in the case when $\LG_X=\SL_2$ or $\PGL_2$ and, hopefully, generalizes to higher rank (although I cannot yet propose such a generalization). The idea is to find a way to compare the relative trace formula for \emph{any} such variety, with the corresponding Kuznetsov formula, i.e., the relative trace formula for the Whittaker model $N_\psi\backslash G^*$ of the group $G^*=\PGL_2$ or $\SL_2$ (respectively). The Kuznetsov formula, not the Arthur--Selberg trace formula, seems to be the appropriate base case for such a type of functoriality, but it requires some modification, because it does not produce on the spectral side the same $L$-functions as the relative trace formula for $X$. Roughly speaking, the spectral side of the Kuznetsov formula is weighted by the factors 
$$ \frac{1}{L(\varphi, \Ad, 1)},$$
(where $\varphi$ denotes a global Langlands parameter into $\LG_X = \LG^*$), while the relative trace formula for $X$ will have an extra $L$-factor, depending on $X$, in the numerator:
$$ \frac{L_X(\varphi)}{L(\varphi, \Ad, 1)}.$$

For example, in the case of the Arthur--Selberg trace formula (when $X=H$, a reductive group), we have $L_X(\varphi) = L(\varphi, \Ad, 1)$, which is why no $L$-values appear in the end, while for $X=\Gm\backslash \PGL_2$ we have $L_X(\varphi)= L(\varphi,\Std, \frac{1}{2})^2$, corresponding to the square of the Hecke period. These $L$-functions are obtained by \emph{enlarging} the space of test measures for the Kuznetsov formula. Thus, our comparison is achieved via a \emph{transfer operator}, which is a linear isomorphism
\begin{equation}
\mathcal T: \mathcal S^-_{L_X}(N_\psi\backslash G^*/N_\psi) \xrightarrow\sim \mathcal S(X\times X/G),
\end{equation}
between the appropriately enlarged space of test measures for the Kuznetsov formula of $G^*$, and the standard space of test measures for the relative trace formula of $X$.

\subsection{Rank-one spherical varieties}

These ideas were explored, in the special cases $X = T\backslash \PGL_2$ (where $T$ is a torus) and $X=\SL_2 = \SO_3\backslash \SO_4$, in the papers \cite{SaBE1, SaBE2, SaHanoi, SaTransfer1, SaTransfer2}, in the case of $X=\SL_2$ generalizing the thesis of Rudnick \cite{Rudnick}. However, it was not clear at that point if those cases were part of a general pattern, or  just reflections of methods already known. In this paper, I demonstrate for the first time that there is a general ``operator of functoriality'' in rank one, as general and uniform as the aforementioned conjectures. 

Spherical varieties of rank one are, in some sense, the building blocks of all spherical varieties, in the same sense that the group $\SL_2$ is the building block of all reductive groups: to a general spherical variety $X$, and each simple coroot $\gamma$ of its dual group (better known as the ``spherical roots'' of $X$), there is an associated rank-one (up to center) spherical variety $X_\gamma$ which is a degeneration of $X$. Thus, the comparisons studied here should be essential in understanding cases of higher rank.

The list of spherical varieties of rank one consists of a finite number of families, classified by Akhiezer in \cite{Akhiezer} --- see also the tables of \cite{Wasserman, KnVS}. Up to the action of the ``center'' $\mathcal Z(X):=\Aut^G(X)$, the \emph{affine homogeneous spherical varieties} $X=H\backslash G$ over an algebraically closed field in characteristic zero whose dual group $\check G_X$ is either $\SL_2$ or $\PGL_2$ are listed in the following table: 

\vspace{0.5cm}~

~\hspace{-3.5cm}
$ \begin{array}{|c|c|c|c|c|c|}
\hline
&X & P(X) & \check G_X & \gamma & L_X \\
 & & & & &\\
\hline
\hline
A_1& \Gm \backslash\PGL_{2} &B &  \SL_{2} & \alpha & L(\Std,\frac{1}{2})^2 \\
\hline
A_n& \GL_n \backslash\PGL_{n+1} &P_{1,n-1,1} &  \SL_{2} & \alpha_1+\dots + \alpha_n& L(\Std,\frac{n}{2})^2 \\
\hline
B_n& \SO_{2n}\backslash\SO_{2n+1} & P_{\SO_{2n-1}} & \SL_2 & \alpha_1 + \dots + \alpha_n&  L(\Std, n-\frac{1}{2}) L(\Std,\frac{1}{2})\\
\hline
C_n& \underset{(n \ge 2)}{\Sp_{2n-2}\times \Sp_{2}\backslash \Sp_{2n}} & P_{\SL_2 \times \Sp_{2(n-2)}} & \SL_{2} & \alpha_1+2\alpha_2+\dots + 2\alpha_{n-1}+ \alpha_n & L(\Std,n-\frac{1}{2}) L(\Std,n-\frac{3}{2}) \\ 
\hline
F_4& \Spin_9\backslash F_4 & P_{\Spin_7} & \SL_2 & \alpha_1+2\alpha_2+3\alpha_3+2\alpha_4 & L(\Std, \frac{11}{2}) L(\Std, \frac{5}{2}) \\
\hline
G_2& \SL_3\backslash G_2& P_{\SL_2} & \SL_2 & 2\alpha_1+ \alpha_2  & L(\Std, \frac{5}{2}) L(\Std, \frac{1}{2}) \\
\hline
\hline
D_2& \SL_2=\SO_{3}\backslash\SO_{4} & B & \PGL_2 & \alpha_1 + \alpha_2 &  L(\Ad,1) \\
\hline
D_n& \SO_{2n-1}\backslash\SO_{2n} & P_{\SO_{2n-2}} & \PGL_2 & 2\alpha_1+ \dots + 2\alpha_{n-2} + \alpha_{n-1}+\alpha_n &  L(\Ad,n-1) \\
\hline
D_4'' & \Spin_7\backslash\Spin_8 & P_{\Spin_6} & \PGL_2 & 2\alpha_1+ 2\alpha_{2} + \alpha_{3}+\alpha_4 & L(\Ad,3)\\
\hline
B_3''& G_2\backslash\Spin_7 & P_{\SL_3} & \PGL_2 & \alpha_1 + 2\alpha_2 + 3\alpha_3 &  L(\Ad,3)\\
\hline
\end{array}
$
\begin{equation}\label{thetable}\end{equation}

The various columns of this table will be explained below. 
In this paper, I work over a local field $F$ in characteristic zero, and will only consider the case when both $G$ and $H$ are split over $F$. Under these restrictions, as we will see (Proposition \ref{uniqueclass?}), each line in the first group of the table above (from $A_1$ to $G_2$) corresponds to a unique isomorphism class of $G$-varieties, while each line in the second group (from $D_2$ to $B_3''$) corresponds to a set of isomorphism classes parametrized by square classes in $F^\times$.

For almost all of the varieties in the table above, a version of the local relative Langlands conjecture for $L^2(X)$ was established by Wee Teck Gan and Ra\'ul Gomez \cite{GG}, on a case-by-case basis using the usual and exceptional theta correspondences. Similar, and other, methods can be used to study global periods for the spaces of Table \eqref{thetable}; examples in the literature include \cite{RS,GanGur,Flicker}.

In any case, the local and global conjectures should be seen as corollaries of a deeper fact, which is encoded in the comparison of trace formulas that I propose in this paper. Moreover, the approach of the present paper is \emph{classification-free} (except for a minor result in Lemma \ref{lemmaclosedorbits}), and relies on a sophisticated theory developed by Friedrich Knop, on the geometry of the moment map 
$$ T^*X \to \g^*.$$

I now explain the various entries in the table: 
The column $\gamma$ denotes the \emph{normalized spherical root} of the spherical variety, in the language of \cite{SV}, described in the basis of simple roots labelled as in Bourbaki. This is the positive coroot for the canonical embedding \eqref{Lgroups}. (The $L$-groups can be replaced by their identity components $\check G_X$, $\check G$, here, since we take $G$ to be split.) This spherical root is either a root of $G$ or the sum of two strongly orthogonal roots; I have chosen the representative of the equivalence class up to center to be such that the dual group is $\SL_2$ in the former case, and $\PGL_2$ in the latter. Thus, we obtain two families of spherical varieties of rank one, whose prototypes are, respectively, the examples labelled $A_1$ and $D_2$ above (which, of course, are special cases of $A_n$ and $D_n$). The case $D_4''$ is obtained from $D_4$ by application of the triality automorphism of $\Spin_8$ (which does not descend to $\SO_8$). Because of the two prototypes, we say (following \cite{SV}) that the spherical root is ``of type $T$'' (for ``torus'') in the first family and ``of type $G$'' (for ``group'') in the second.

By $P(X)$ we denote the conjugacy class of parabolics stabilizing the open Borel orbit. In the table, I describe the parabolics in a way that should be self-explanatory, by indicating the semisimple part of their Levi quotient $L(X)$ or, in the case of $\GL_n$, an ordered partition of $n$. In the case of $G=G_2$, $P_{\SL_2}$ is such that its Levi contains a long root. Notice that the roots of the Levi $L(X)$ are always orthogonal to the spherical root $\gamma$.
The parabolic $P(X)$ determines the restriction of the map \eqref{Lgroups} to the ``Arthur-$\SL_2$'' factor, which has to map to a principal $\SL_2$ in the Levi subgroup of $\check G$ dual to $P(X)$. 

Finally, \label{refLX} $L_X$ stands for a $\frac{1}{2}\ZZ$-graded representation $r=\bigoplus_{n \in \frac{1}{2}\ZZ} r_n$ of $\check G_X$, which I call the ``$L$-value associated to $X$'', thinking of the $L$-value 
$$\prod_n L(r_n\circ \varphi, n)$$ attached to any Langlands parameter $\varphi$ into $\check G_X$. For this reason, I denote this graded representation by $\prod_n L(r_n, n)$.
This is the $L$-value attached to the square of the global $H$-period, according to the generalization of the Ichino--Ikeda conjecture \cite{II}  proposed in \cite[\S 17.4]{SV} and the local unramified calculation, performed for classical groups only, of \cite{SaSph}.  It would be desirable to have the same calculation generalized to all cases, including non-classical groups. In any case, the $L$-value here will be determined directly in terms of the geometry of the $G$-space $X\times X$, as follows:

\begin{itemize}
 \item When the spherical root is of type $T$, the associated $L$-value is always of the form $L(\Std, s_1)L(\Std, s_2)$, for some positive half-integers $s_1\ge s_2$. These are determined by the relations $s_1+s_2=\frac{\dim X}{2}$, and $  s_1 = \frac{\left<\check\gamma, \rho_{P(X)}\right>}{2}$.
 \item When the spherical root  is of type $G$, we have $L_X = L(\Ad, s_0)$, with $s_0= \left<\check\gamma, \rho_{P(X)}\right> = \frac{\dim X-1}{2}$, always an integer.
\end{itemize}
Here, $2\rho_{P(X)}$ is the sum of positive roots in the nilpotent radical of the Lie algebra of $P(X)$; it can be considered as a cocharacter into the canonical maximal torus of $\check G_X$, hence the value of the spherical coroot $\check \gamma$ (the positive root of $\check G_X$) makes sense on it.

\subsection{Notation and the main result} \label{ssnotation}

I introduce the notions and notation necessary in order to formulate the main result. A more complete index of notation appears in Appendix \ref{index}.

All varieties will be defined over a local, locally compact field $F$ in characteristic zero, and we write $X=X(F)$, etc, when no confusion arises. In particular, all measures or functions will be on the $F$-points of the varieties under consideration.

We denote \label{refS} by $\mathcal S(X)$ the space of ($\CC$-valued) Schwartz \emph{measures} on the $F$-points of a smooth variety $X$; these are smooth measures which, in the non-Archimedean case, are of compact support, and in the Archimedean case are of rapid decay, together with their polynomial derivatives. For the definitions in the Archimedean case, which use the semialgebraic structure of $X(F)$, I point the reader to \cite{AGSchwartz}, in particular \S A.1.1. (They are not to be confused with measures/functions in the Harish-Chandra Schwartz space, which are of slower decay.) Notice that, in the Archimedean case, when $F=\mathbb C$ we need to consider $X(F)$ as a \emph{real} semialgebraic manifold, since the complex semialgebraic topology is not fine enough. Moreover, in this case the space $\mathcal S(X)$ has a natural Fr\'echet topology; if $X$ is admits a nowhere vanishing Nash (smooth semialgebraic) density $\omega$, so that every element of $\mathcal S(X)$ can be written as $\Phi\cdot \omega$ for a Schwartz function $\Phi$, the topology is induced by the seminorms 
$$ \Vert \Phi \omega \Vert_D = \sup |D\Phi|,$$
where $D$ ranges over all polynomial Nash differential operators. The general case can be reduced to that by a Nash partition of unity \cite[Theorem 4.4.1]{AGSchwartz}.

For uniformity of language, I will often write ``smooth of rapid decay'' to describe this behavior of Schwartz measures, with the understanding that this means compact support in the non-Archimedean case, and that it includes their derivatives (the functions $D\Phi$ above) in the Archimedean case. Whenever needed, the space of Schwartz \emph{functions} will be denoted by $\mathcal F(X)$.

The notation $X\sslash G$ will stand for the affine, invariant-theoretic quotient $\Spec F[X]^G$ of a $G$-variety $X$, and if $\pi: X\to X\sslash G$ denotes the canonical quotient map, the image $\pi_* \mathcal S(X)\subset \Meas(X\sslash G)$ of the pushforward map of measures will be denoted \label{refSmodG} by $\mathcal S(X/G)$. In the Archimedean case, where $\mathcal S(X)$ is a nuclear Fr\'echet space, the space $\mathcal S(X/G)$ inherits a quotient Fr\'echet topology; in the non-Archimedean case, any reference to topology should be ignored.

Let $X$ be one of the spaces in Table \eqref{thetable}, with a reductive group $G$ acting on it, and let $G^*$ \label{refGstar} denote the group $\PGL_2$, if $\check G_X=\SL_2$, or $\SL_2$, if $\check G_X=\PGL_2$. Let $N\subset G^*$ be the upper triangular unipotent subgroup, identified with the additive group $\Ga$ in the obvious way, and let $\psi: F\to \CC^\times$ be a nontrivial character, considered also as a character of $N$. We fix throughout an additive Haar measure on $F$, which is self-dual with respect to $\psi$. We extend the notation of Schwartz spaces to the quotient that we will denote by $N_\psi\backslash G^*/N_\psi$: If $A^*\subset G^*$ is the torus of diagonal elements, and $w = \begin{pmatrix} & -1 \\ 1\end{pmatrix}$, we embed $A^*$ in the affine quotient $N\backslash G^*\sslash N$ by $A^*\ni a \mapsto [wa]$, the class of the element $wa$, and let $\mathcal S(N_\psi\backslash G^*/N_\psi)$ \label{refSKuz} denote the space of measures on $A^*$ of the form 
\begin{equation} f(a) = \pi^\psi_*(\Phi dg)(a):= \left(\int_{N\times N} \Phi(n_1 wa n_2) \psi^{-1}(n_1n_2) d(n_1,n_2)\right)\cdot \delta(a) da,
\end{equation}
where $\Phi$ is a Schwartz \emph{function} on $G^*$, $\delta$ is the modular character of the upper triangular Borel subgroup (the quotient of the right by the left Haar measure), and $da$ is a Haar measure on $A^*$. This is a twisted version of the pushforward of the measure $\Phi dg$ to $N\backslash G^*\sslash N$, which for suitable Haar measures reads:
$$ \pi_* (\Phi dg)(a) = \left(\int_{N\times N} \Phi(n_1 wa n_2) d(n_1,n_2)\right)\cdot  \delta(a) da.$$

We also fix a coordinate on $A^*$ which we denote by 
$$ \xi(a) = e^\alpha(a), \mbox{ when } G^*=\PGL_2;$$
$$ \zeta (a) = e^\frac{\alpha}{2}(a), \mbox{ when } G^*=\SL_2,$$
where $\alpha$ is the positive (upper triangular) root of $G^*$, and we use exponential notation to denote the corresponding character, since weights are written additively. Then, $N\backslash G^*\sslash N$ is identified with $\mathbbm A^1$, with coordinates $\xi$, resp.\ $\zeta$, the images of the elements
$$\begin{pmatrix} & -1 \\ \xi\end{pmatrix}, \,\,\,
\begin{pmatrix} & -\zeta^{-1} \\ \zeta\end{pmatrix},$$
respectively (when $\xi, \zeta\ne 0$).

The elements of $\mathcal S(N_\psi\backslash G^*/N_\psi)$, viewed as measures on $\mathbbm A^1=N\backslash G^*\sslash N$, are smooth of rapid decay away from zero, while in a neighborhood of $0\in \mathbbm A^1$ they have a singularity which is called the ``Kloosterman germ'', because in the non-Archimedean case they are smooth multiples of the measures
$$ \xi\mapsto \left(\int_{|u|^2=|\xi|}  \psi^{-1} \left(\frac{u}{\xi} + u^{-1}\right) du \right) d^\times \xi,$$
resp.\
$$ \zeta\mapsto \left(\int_{u\in \pm 1 + \mathfrak p} \psi^{-1} \left(\frac{u+u^{-1}}{\zeta}\right) du \right) d\zeta.$$
(See \cite[Proposition 4.9.2]{SaBE1} for the former, and the latter is similar. Notice that there are two separate germs in the case of $\SL_2$, corresponding to the choice of $\pm 1$.)

We define \label{refenlarged} enlarged spaces of test measures 
$$\mathcal S^-_{L_X} (N_\psi\backslash G^*/N_\psi) \supset  \mathcal S (N_\psi\backslash G^*/N_\psi)$$ 
for the Kuznetsov formula, associated to the $L$-values $L_X$ that appear in Table \eqref{thetable}, as follows (see also \cite[\S 2.2]{SaTransfer1}): their elements coincide with elements of $\mathcal S (N_\psi\backslash G^*/N_\psi)$ away from infinity, but in a neighborhood of infinity, instead of being of rapid decay, they are allowed to be of the following form:
\begin{itemize}
 \item When $G^*=\PGL_2$ and $L_X=L(\Std, s_1)L(\Std,s_2)$ with $s_1\ge s_2$,
 \begin{equation}\label{expansionPGL2}
  (C_1(\xi^{-1}) |\xi|^{\frac{1}{2}-s_1}+C_2(\xi^{-1}) |\xi|^{\frac{1}{2}-s_2}) d^\times \xi,
 \end{equation}
 where $C_1$ and $C_2$ are smooth functions; this should be replaced by 
 \begin{equation}
    |\xi|^{\frac{1}{2}-s_1}(C_1(\xi^{-1}) +C_2(\xi^{-1}) |\xi|^{s_1-s_2} \log|\xi|) d^\times \xi
 \end{equation}
 when $|\xi|^{s_1-s_2}$ is a smooth function --- that is, 
 \begin{equation}\label{logconditions}
  \parbox{25em}{
 \begin{itemize}
  \item in the non-Archimedean case, when $s_1=s_2$;
  \item in the real case, when $s_1-s_2\in 2\mathbb N$;
  \item in the complex case, when $s_1-s_2\in \mathbb N$.
  \end{itemize}
  }
 \end{equation}
 (We use the arithmetic normalization of absolute values, which is compatible with norms to the base field; this is the square of the usual absolute value in the complex case.)
 \item When $G^*=\SL_2$ and $L_X=L(\Ad,s_0)$, 
 \begin{equation}\label{expansionSL2}
  C(\zeta^{-1}) |\zeta|^{1-s_0} d^\times \zeta,
 \end{equation}
 where $C$ is a smooth function.
\end{itemize}

In the Archimedean case, all of these spaces have an obvious Fr\'echet topology, which by a partition of unity can be reduced to the topology of the standard space of Kuznetsov test measures $\mathcal S(N_\psi\backslash G^*/N_\psi)$ (topologized as a quotient of $\mathcal S(G^*)$), and seminorms on the Schwartz functions $C, C_1, C_2$ in the above asymptotic expressions at $\infty$.
The basic theorem proven in this paper is the following:

\begin{theorem}\label{maintheorem}
Let $\C_X=(X\times X)\sslash G$. There is an isomorphism $\C_X \simeq \mathbbm A^1$, and the map $X\times X\to \mathbbm A^1$ is smooth away from the preimage of two points of $\mathbbm A^1$, that we will call singular. We fix the isomorphisms as follows:
\begin{itemize}
 \item When $\check G_X=\SL_2$, we take the set of singular points to be $\{0, 1\}$, with $X^\diag\subset X\times X$ mapping to $1 \in \C_X\simeq \mathbbm A^1$.
 \item When $\check G_X = \PGL_2$, we take the set of singular points to be $\{-2,2\}$, with $X^\diag\subset X\times X$ mapping to $2 \in \C_X\simeq \mathbbm A^1$.
\end{itemize}

Then, there is a continuous linear isomorphism: 
 \begin{equation}
  \mathcal T: \mathcal S^-_{L_X}(N_\psi\backslash G^*/N_\psi) \xrightarrow\sim \mathcal S(X\times X/G),
 \end{equation}
 given by the following formula:
 
\begin{itemize}
 \item When $\check G_X=\SL_2$ with $L_X=L(\Std, s_1) L(\Std,s_2)$, $s_1\ge s_2$,
 \begin{equation}
    \mathcal Tf(\xi) =  |\xi|^{s_1-\frac{1}{2}}  \left( |\bullet|^{\frac{1}{2}-s_1} \psi(\bullet) d\bullet\right) \star  \left( |\bullet|^{\frac{1}{2}-s_2} \psi(\bullet) d\bullet\right) \star f(\xi).
 \end{equation}
 \item When $\check G_X=\PGL_2$ with $L_X=L(\Ad, s_0)$,
 \begin{equation}
    \mathcal Tf(\zeta) =  |\zeta|^{s_0-1}  \left( |\bullet|^{1-s_0} \psi(\bullet) d\bullet\right) \star  f(\zeta).
 \end{equation}
\end{itemize}

\end{theorem}

By $\left( |\bullet|^s \psi(\bullet) d\bullet\right) \star $ we denote the operator of \emph{multiplicative} convolution by the measure $\left( |x|^s \psi(x) dx\right)$ (in the variable $y=\xi$ or $\zeta$, respectively): 
$$\left( |\bullet|^s \psi(\bullet) d\bullet\right) \star f (y) = \int_{F^\times} |x|^s \psi(x) f(x^{-1}y) dx = |y|^{s+1} \int f(u^{-1}) |u|^s \psi(uy) du.$$
The measure $dx$ is the \emph{additive} Haar measure on $F$ that we have fixed. If these integrals do not converge, convolution should be understood as the Fourier transform of the distribution $u\mapsto f(u^{-1}) |u|^s$, followed by multiplication by $|y|^{s+1}$. 

The operator $\mathcal T$ 
is clearly the correct operator of functoriality between the relative trace formula for $X$ and the Kuznetsov formula. Indeed, it was shown in \cite{SaBE1, SaBE2, SaTransfer1} that in the basic cases $A_1, D_2$ it satisfies the appropriate fundamental lemma for the Hecke algebra, and that it pulls back relative characters to relative characters (see \cite[\S 6--7]{SaHanoi} for precise references); these statements can also be confirmed in the general $A_n$-case by ``unfolding''. In an upcoming paper \cite{JK}, Daniel Johnstone and Rahul Krishna prove the appropriate fundamental lemma for the transfer operator in all cases. There remains to prove the fundamental lemma for the full Hecke algebra, in order to be able to use this operator globally (together with the ``Hankel transforms'' for the functional equations of the standard and adjoint $L$-functions, discussed in \cite[\S 8]{SaHanoi} and \cite[\S 8]{SaTransfer2}), and obtain a uniform proof of functoriality and the relation between $X$-periods of automorphic forms and the $L$-value $L_X$. After a version of this paper appeared, Gan and Xiaolei Wan confirmed that the theta correspondence descends to the transfer operators introduced here, in the case of $X=\SO_n\backslash \SO_{n+1}$ \cite{GW}; their methods can be applied to other cases, as well. It is quite satisfactory to observe that various different methods for comparing period integrals, such as the ``unfolding'' method and the theta correspondence, all descend to the same statements at the level of relative trace formulas!

\subsection{Outline of the proof}

As mentioned, the proof of Theorem \ref{maintheorem} is classification-free, and relies on Friedrich Knop's theory of the moment map. The main issue is to analyze the quotient $X\times X/G$, and to describe, more or less explicitly, the germs of pushforward Schwartz measures for the map $X\times X \to X\times X\sslash G$. 

To every spherical variety $X$ one can attach a canonical ``universal Cartan'', that is, a torus $A_X$, and a ``little Weyl group'' $W_X$ acting on it.  The dual torus to $A_X$ is the canonical maximal torus of the dual group $\check G_X$, and $W_X$ is its Weyl group. Hence, in the rank-one case that we are considering, $A_X\simeq \Gm$ and $W_X = \ZZ/2$, acting on $\Gm$ by inversion. 

The main result about the space $\mathcal S(X\times X/G)$ of pushforward measures, for $X$ as in Table \eqref{thetable}, is the following:
\begin{theorem}\label{Xtheorem}
 There is a canonical isomorphism $\C_X:= X\times X\sslash G\simeq A_X\sslash W_X$, and the map $X\times X\to \C_X$ is smooth away from the preimages of $[\pm 1]$, where $[\pm 1]$ denote the images of $\pm 1\in A_X$ in $A_X\sslash W_X$. 
 
 In particular, there are two distinguished closed $G$-orbits $X_1 = X^\diag$ and $X_{-1}$ (over $[\pm 1]$, respectively); if $d_{\pm 1}$ denote their codimensions, then $d_1=\dim X$ and 
 $$d_{-1} = \epsilon \left<2\rho_{P(X)},\check\gamma\right>-d_1+2,$$ 
 where $\check\gamma$ is the spherical coroot, $2\rho_{P(X)}$ is the sum of roots in the unipotent radical of $P(X)$, and
$$ \epsilon = \begin{cases} 1, \mbox{ when the spherical root is of type $T$ (dual group $\SL_2$)};\\
2, \mbox{ when the spherical root is of type $G$ (dual group $\PGL_2$).}
\end{cases}$$
In the case of root of type $G$, $d_1=d_{-1}$.

 The space $\mathcal S(X\times X/G)$ consists of those measures on $\C_X(\simeq \mathbbm A^1)$ which are smooth and of rapid decay, together with their polynomial derivatives, away from neighborhoods of $[\pm 1]$ (compactly supported in the non-Archimedean case), while in the neighborhood of $[\pm 1]$ their germs coincide with germs for the twisted pushforward maps:
 $$\mathbbm A^2/(\Gm, |\bullet|^\frac{2-d_{\pm 1}}{2}),$$
 for spherical root of type $T$, and
 $$\mathfrak{sl}_2/(B_\ad, \delta_2^\frac{3-d_{\pm 1}}{2}),$$ 
 for spherical root of type $G$, where $B_\ad$ denotes the Borel subgroup of $\PGL_2$, and $\delta_2$ is its modular character.
\end{theorem}

For the precise meaning of these ``twisted pushforwards'', I point the reader to the precise formulation of Theorem \ref{Xtheorem2}. In other words, the germs for the general case are twisted versions of the germs for the ``basic cases'' $A_1$ and $D_2$. This indirect description of the germs allows us to relate these germs of pushforward measures for $X\times X/G$ with the Kloosterman germs for the Kunzetsov formula of $G^*$, based on results of \cite{SaBE1, SaTransfer2}.

Since $X=H\backslash G$ is homogeneous, we can also write $X\times X\sslash G = H\backslash G\sslash H$; when $X$ is symmetric (as is the case for most of the cases in Table \eqref{thetable}, except for those denoted by $G_2$ and $B_3''$), the identification of this with $A_X\sslash W_X$ is due to Richardson \cite{Richardson}. 

In any case, to obtain this and Theorem \ref{Xtheorem} in general, I use Knop's theory of the moment map in a somewhat paradoxical way: While the cotangent bundle together with its moment map $T^*X\to \g^*$ is classically used for microlocal analysis on $X$, here I use it to obtain an explicit resolution\footnote{The space $X\times X$ is smooth, but here we take into account the $G$-action, and use the term ``resolution'' to refer to the fibers of the quotient map $X\times X\to X\times X\sslash G$: a resolution is a blowup that turns them into normal crossings divisors.} of the space $X\times X$ under the $G$-action. The basic idea is, roughly, to study the space
$$ Z:=T^*X\times_{\g^*} T^*X,$$
which is the union of conormal bundles to the $G$-orbits on $X\times X$. Where the $G$-orbits are of codimension one, their conormal bundles are of dimension one, and the map from the projectivization:
$$ \P Z \to X\times X$$
restricts to an isomorphism. The important issue is to understand the conormal bundles where this map fails to be an isomorphism.

It eventually turns out that $\P Z$ is not quite the correct resolution, because it can be quite singular. A closely related space is a space that I denote by $\P J_X$, and which is obtained in Section \ref{sec:resolution} as follows:

Let $\a_X^*$ be the dual Lie algebra to the torus $A_X$, and $\c_X^*=\a_X^*\sslash W_X$ --- both of these spaces are isomorphic to the affine line. There is a smooth abelian group scheme $J$ over $\c_X^*$ whose general fiber is isomorphic to $A_X$, but the isomorphism is only determined up to the action of $W_X$, and whose special fiber (over $0=$ the image of $0\in \a_X^*$) is isomorphic to $\{\pm 1\} \times \Ga$. This group scheme is known, for example, as the group scheme of regular centralizers over the Kostant section of $\mathfrak{sl}_2$, and it can be abstractly defined as
$$J = \left(\Res_{\a_X^* / \c_X^*} (A_X \times \a_X^*)\right)^{W_X},$$
where $\Res_{\a_X^* / \c_X^*}$ denotes Weil restriction of scalars from $\a_X^*$ to $\c_X^*$.

Knop has shown \cite{KnAut} that, except perhaps for the non-identity component of the nilpotent fiber of $J$, this group scheme acts canonically on $T^*X$ over $\g^*$. (This action is canonical in that its differential is the Hamiltonian vector field induced from canonical isomorphisms $T^*X\sslash G\xrightarrow\sim \c_X^*$ and $\operatorname{Lie}(J) \simeq T^* \c_X^*$.)
Thus, being a bit imprecise as far as the action of the non-identity component of the nilpotent fiber goes, we have a map
$$ J\times_{\c_X^*} T^*X \to T^*X\times_{\g^*} T^*X,$$
and it is its composition with the map to $X\times X$ (after projectivization) which will give rise to the desired resolution:
$$ \P(J\times_{\c_X^*} T^*X) \to X\times X.$$

On the other hand, we have, by definition, a canonical quotient map $J\to A_X\sslash W_X$, and this can be used to prove the isomorphism $X\times X\sslash G\xrightarrow\sim A_X\sslash W_X$. 

Recall that there is a bijection between points of $X\times X\sslash G$ and closed (geometric) orbits of $G$ on $X\times X$. The diagonal $X_1=X^\diag\hookrightarrow X\times X$ corresponds to the class of $1\in A_X$, and the fiber of $J\times_{\c_X^*} T^*X$ over it is just its conormal bundle $N^*_{X_1} = T^*X$. To correct the imprecision about the non-identity component $\{-1\}\times \Ga$ of the nilpotent fiber of $J$, we replace $(\{-1\}\times \Ga) \times_{\c_X^*} T^*X$ by a copy of $\Ga \times_{\c_X^*} N_{X_{-1}}^*$, where $X_{-1}$ denotes the closed $G$-orbit over $[-1]\in
A_X\sslash W_X$, $N^*_{X_{-1}}$ denotes its conormal bundle, and $\Ga$ maps to $0\in \c_X^*$. This replacement leads to a smooth scheme $J_X\rightrightarrows T^*X$, birational to $J\times_{\c_X^*} T^*X$, such that the resulting map from its ``projectivization''
$$ \P J_X \to X\times X,$$
is isomorphic to the blowup of $X\times X$ at the closed orbits $X_1$ and $X_{-1}$.

The formula of Theorem \ref{Xtheorem} on the codimensions of orbits is obtained in Section \ref{sec:integration} by a degeneration argument --- developing the analog of the Weyl integration formula for $X\times X$ under the diagonal $G$-action, and deforming $X$ to its horospherical ``boundary degeneration'' $X_\emptyset$, where this integration formula is very explicit.

The map $\P J_X \to A_X\sslash W_X$ is easy to describe, and a standard analysis of pullbacks of Schwartz measures under resolutions shows, in Section \ref{sec:Schwartzmeasures}, that the elements of the pushforward space $\mathcal S(X\times X/G)$ are measures on $A_X\sslash W_X\simeq \mathbbm A^1$ whose singularities at $[\pm 1]$ are linear combinations of multiplicative characters of the form $x\mapsto |x|^{\frac{d_{\pm 1}}{2}-1} \eta(x)$, where $\eta$ is quadratic.

The last task, in Section \ref{sec:germs} is to understand this linear combination of these characters. By linearization, this is equivalent to understanding the pushforwards of Schwartz measures under a map
$$ V\xrightarrow{Q} \mathbbm A^1,$$
where $Q$ is a nondegenerate, split quadratic form on a vector space $V$ of dimension $d=d_{\pm 1}$. A key proposition, \ref{reductiontobasic}, identifies these pushforwards with twisted pushforwards on a two- or three-dimensional quadratic space, as in Theorem \ref{Xtheorem}.

\subsection{Acknowledgments} This work started as a joint project with Daniel Johnstone and Rahul Krishna, who eventually undertook the proof of the fundamental lemma. I thank them for many helpful conversations, including an observation of R.\ Krishna which simplified the proof of the key proposition \ref{reductiontobasic}. I am grateful to the University of Chicago and Ng\^o Bao Ch\^au for their hospitality during the winter and spring quarters of 2017, at the end of which I discovered the existence of a uniform transfer operator, and to the Institute for Advanced Study where most of the writing was done during the academic year 2017--2018. I also thank Wee Teck Gan for providing several references, and an anonymous referee for a very detailed and helpful report. Finally, I would like to acknowledge my intellectual debt to Herv\'e Jacquet and Friedrich Knop, two formidable mathematicians who, since the 80s and 90s, from different perspectives, laid the ground for the questions that I am addressing in this paper, not always properly appreciated by the mainstream of the field; I believe that we have only seen the tip of the iceberg as far as relations between their work go.

This work was supported by NSF grants DMS-1502270, DMS-1801429, DMS-1939672, and by a stipend at the IAS from the Charles Simonyi Endowment.

\section{The moment map and the structure of Borel orbits}

\subsection{Invariant theory of the cotangent bundle and its polarizations} \label{ssmomentmap}

Throughout the paper, $X$ will denote one of the homogeneous spherical varieties of rank one appearing in Table \eqref{thetable}. However, in this section I revisit (and slightly reformulate) the theory of the cotangent bundle of $X$ due to Friedrich Knop, which holds true for any homogeneous, quasi-affine spherical variety $X$ under the action of a connected reductive group $G$. 

To any such $X$, one attaches a conjugacy class of parabolics, \label{refPX} denoted by $P(X)$, characterized by the property that, if $B\subset G$ is a Borel subgroup and $\mathring X\subset X$ its open orbit, (a representative of) $P(X)$ is given by 
$$P(X):= \{g\in G| \mathring Xg =\mathring X\}.$$
The Levi quotient (and sometimes, a Levi subgroup) of $P(X)$ will be denoted by $L(X)$. 

Let $A$ denote the reductive quotient of a Borel subgroup --- viewed as a canonical torus associated to $G$, up to unique isomorphism, the so-called \emph{(universal) Cartan of $G$}. 
Let $\B$ denote the full flag variety of $G$, and let $\B_X$ denote the flag variety of parabolics in the conjugacy class of $P(X)$. 

The unipotent radical of a parabolic $P$ will be denoted by $U_P$. Having fixed a Borel subgroup $B$, we will also use the letter $N$ for $U_B$. The quotient $\mathring X\sslash N$ is a homogeneous space under the action of $A$; its action factors through the faithful action of a quotient $A\twoheadrightarrow A_X$ \label{refAX} which we will call \emph{the (universal) Cartan of $X$}. In fact, it is known that $A_X$ is a quotient of $P(X)$:
$$ P(X)\to L(X)\to A_X,$$
and that $P(X)$ acts on $\mathring X\sslash N$ through this quotient. We will denote \label{refL1} by $L_1$ the kernel of the map $L(X)\to A_X$.

The rank of $A_X$ is, by definition, the rank of $X$; thus, for all varieties of Table \eqref{thetable}, $A_X\simeq\Gm$. The character group of $A_X$ will be denoted by $\Lambda_X$, and called the weight lattice of $X$. We use similar notation for other $B$-orbits (or $B$-orbit closures) $Y$: $\Lambda_Y$ will denote \label{refLambdaY} the set of characters of nonzero rational $B$-eigenfunctions on $Y$, and $A_Y = \Spec F[\Lambda_Y]$ the torus quotient by which $A$ acts on $Y\sslash N$. The \emph{rank} of $Y$ is the rank of the group $\Lambda_Y$. It is known that $\mathring X$ has maximal rank among all $B$-orbits on $X$.

We will denote Lie algebras of algebraic groups by the same letter in Gothic lowercase, and linear duals by a star exponent. The cotangent space $T^*X$ of $X$ comes equipped with a moment map 
$$\mu: T^*X\to \g^*,$$
which, for $X=H\backslash G$, is simply the $G$-equivariant extension of the embedding $\mathfrak h^\perp\subset \g^*$ to $T^*X = \mathfrak h^\perp \times^H G$
This gives rise to a $G$-invariant map:
$$ T^*X \to \g^*\sslash G= \a^*\sslash W.$$
We let\label{refhatg} $\hat\g^*=\g^*\times_{\a^*\sslash W} \a^*$, and $\tilde\g^*=\{(Z,B)| Z\in \g^*, B\in \B, Z\in \u_B^\perp\}$; the latter is the Springer--Grothendieck resolution, and the canonical quotient $\u_B^\perp\to \a^*$ induces natural, proper maps $\tilde\g^*\to \hat\g^*\to\g^*$.

We define \label{refhatTX} the following covers of the cotangent bundle:
\begin{itemize}
 \item $\widehat{T^*X} = T^*X \times_{\a^*\sslash W} \a^* = T^*X \times_{\g^*} \hat\g^*$ (the \emph{polarized} cotangent bundle);
 \item $\widetilde{T^*X} = \{(v,B)| v \in T^*X, B\in \B, \mu(v) \in \u_B^\perp\} = T^*X \times_{\g^*} \tilde\g^*$.
\end{itemize}
Hence, we have proper maps $\widetilde{T^*X}\to \widehat{T^*X} \to T^*X$, as base changes of the corresponding maps between covers of $\g^*$.

Following Knop, we construct canonical maps, that we will call \emph{Knop's sections},  \label{refKsec}
\begin{eqnarray}
\hat\kappa_X: (X\times \B_X)^0 \times \a_X^* & \to & \widehat{T^*X}, \\
\tilde\kappa_X: (X\times \B)^0 \times \a_X^* & \to & \widetilde{T^*X}
\end{eqnarray}
over $X$, linear in the $\a_X^*$-argument,
where the exponent ``$0$'' denotes the subset of pairs $(x, P)$ with $x$ in the open $P$-orbit. The maps are given as follows: by linearity, it is enough to define them for the lattice $\Lambda_X=\Hom(A_X,\Gm)\subset \a_X^*$, where we consider characters as elements of $\a_X^*$ by identifying them with their differentials at the identity. Let $\chi\in \Lambda_X$, let $P\in \B_X$ or $P\in \B$, respectively, and let $f_\chi$ be a rational, nonzero $P$-eigenfunction on $X$ with eigencharacter $\chi$. If $x$ is in the open $P$-orbit, then 
$$ \hat\kappa_X(x,P,\chi) = (d_x \log f_\chi,\chi) \in \widehat{T^*X} = T^*X \times_{\a^*\sslash W} \a^*, $$
and
$$ \tilde\kappa_X(x, P,\chi) = (d_x \log f_\chi,P) \in \widetilde{T^*X} \subset T^*X \times \B,$$
where $d_x$ denotes the differential evaluated at $x$, and $d_x \log f_\chi= \frac{d_x f_\chi}{f_\chi(x)}$. The term ``sections'' is due to the fact that $\tilde\kappa_X$ is a partial section of the natural map $\widetilde{T^*X}\to X\times \B\times \a^*$; we also apply it to $\hat\kappa_X$ by abuse of language, despite the fact that it is a partial section of some map only over an open subset $\mathring\a_X^*$ of ``regular elements'' (to be defined below).

The following facts are known, or can easily be inferred, as we indicate below, from the work of Knop: 
\begin{enumerate}
 \item All maps $\widetilde{T^*X}\to \widehat{T^*X} \to T^*X$ are proper and dominant. 

 This is obvious from the definitions.

 \item There is, by definition, a natural map $\widetilde{T^*X}\to X\times \B$. The image of any irreducible component is an irreducible, $G$-stable subset, which therefore contains a largest $G$-orbit, giving rise to a map 
 $$ \{\mbox{irreducible components of }\widetilde{T^*X}\} \to \{G\mbox{-orbits in }X\times \B\}.$$
 If we fix $B\in \B$, $G$-orbits on $X\times \B$ are in bijection with $B$-orbits on $X$. Under this map, the irreducible components \emph{of maximal dimension} in $\widetilde{T^*X}$ are in bijection with the \emph{Borel orbits of maximal rank} in $X$. 
 
 This is \cite[Proposition 6.3]{KnOrbits}.  We will denote by $\widetilde{T^*X}^\bullet$ \label{refTXbullet} the irreducible component corresponding to the open Borel orbit; it is the closure of the image of $\tilde\kappa_X$ in $\widetilde{T^*X}$.

 \item Considering only these components of maximal dimension in $\widetilde{T^*X}$, and their images in $\widehat{T^*X}$, we obtain a canonical bijection between \emph{the irreducible components of $\widehat{T^*X}$} and the \emph{Borel orbits of maximal rank in $X$}. 

 This is \cite[Theorem 6.4]{KnOrbits}, together with the non-degeneracy statement of \cite[Lemma 3.1]{KnMotion}. 
Thus, this bijection is characterized by the fact that the component corresponding to a $B$-orbit $Y$ contains all pairs $(v\in T^*_Y X, Z\in \a^*)$ with $\mu(v)\in \mathfrak u_B^\perp$ and $Z$ its image under the canonical map $\mathfrak u_B^\perp\to \a^*$. In particular, the closure of the image of Knop's section $\hat\kappa_X$ is the irreducible component corresponding to the open $B$-orbit, to be denoted by $\widehat{T^*X}^\bullet$.

 \item The stabilizer of $\widehat{T^*X}^\bullet$ under the natural action of $W$ on $\widehat{T^*X}$ (induced from its action on $\a^*$) is a certain semidirect product $W_{L(X)}\rtimes W_X$, and the image of $\widehat{T^*X}^\bullet$ in $\a^*$ coincides with $\a_X^*$. Here, $W_{L(X)}$ is the Weyl group of the Levi quotient of $P(X)$, which is the largest subgroup of $W$ acting trivially on $\a_X^*$, and $W_X$ is the so-called \emph{little Weyl group} of the spherical variety, which acts faithfully on $\a_X^*$. For the examples of Table \eqref{thetable}, $W_X=\ZZ/2$.
 
 The fact that $W_{L(X)}$ is precisely the centralizer of $\a_X^*$ again follows from the non-degeneracy statement of \cite[Lemma 3.1]{KnMotion}; I point the reader to the proof of \cite[Theorem 6.2]{KnOrbits} for the other statements.

\end{enumerate}

The image of the moment map, followed by the Chevalley quotient: 
$$
 T^*X \xrightarrow{\mu} \mathfrak g^* \to \mathfrak a^*\sslash W
$$
is equal to the image of the map
$$
\mathfrak a_X^*\sslash W_X \to \mathfrak a^*\sslash W,
$$
induced from the inclusion $\mathfrak a_X^*\hookrightarrow \mathfrak a^*$. Knop has shown that the map to $\a^*\sslash W$ lifts canonically to a map \label{refmuinv}
\begin{equation} \label{invmoment}
 \mu_\inv: T^*X \to \c_X^*:=\a_X^*\sslash W_X,
\end{equation}
descending from the map $\widehat{T^*X}^\bullet\to \a_X^*$. Under \eqref{invmoment}, $\c_X^*$ is identified with the invariant-theoretic quotient $T^*X\sslash G$ \cite[Satz 7.1]{KnWeyl}. The map $\mu_\inv$ is 
the \emph{invariant moment map}. Thus, we have a commutative diagram
\begin{equation}\label{momentfact}\xymatrix{
 T^*X \ar[r]\ar[dr]_{\sslash G} &  \g_X^* \ar[r]\ar[d]^{\sslash G} & \g^* \ar[d]^{\sslash G}\\ 
 & \c_X^*= \a^*_X\sslash W_X \ar[r] & \a^*\sslash W,
 }
\end{equation}
where $\g_X^*$ is \label{refgX} the spectrum of the integral closure of the image of $F[\g^*]$ in $F[T^*X]$. (I point the reader to \cite[\S 6]{KnWeyl} for the definition of $\mathfrak g_X^*$, denoted there by $M_X$, and the map $\mathfrak g_X^*\to \mathfrak c_X^*$.)

We let $\mathring\a_X^*$ denote the open $W_X$-stable subset where $W_X$ acts freely and $\mathring\c_X^*$ its image --- in our rank-one cases, these are just the complements of zero. Vectors in $T^*X$ (and its various covers) or $\mathfrak g_X^*$ which live over $\mathring{\c}_X^*$ will be called \emph{regular semisimple}, \label{refrs} and denoted $T^*X^\rs$, resp.\ $\mathfrak g_X^{*,\rs}$. The reader should not confuse this notion with the property of being regular in $\g^*$; in fact, the centralizer of the image of an element of $T^*X^\rs$ in $\g^*$ is conjugate to a Levi of $P(X)$ over the algebraic closure. Hence, ``regular semisimple'' elements in $T^*X$ have an image in $\g^*$ which is semisimple and ``as regular as possible'', though not necessarily regular.

Now we restrict to the case when $X$ has rank one and $W_X=\ZZ/2$. Thus, $\a_X^*$ is a one-dimensional vector space, and $\mathring\a_X^*$ is the complement of zero. We can also identify $\c_X^*$ with $\mathbbm A^1$, always letting the point $0\in \c_X^*$ (the image of $0\in \a_X^*$) correspond to $0\in \mathbbm A^1$. Then, the invariant moment map $T^*X\to \c_X^*$ can be considered as a quadratic form on the fibers of $T^*X$.

\begin{lemma}\label{nondegenerate}
 For $X$ affine homogeneous of rank one with $W_X=\ZZ/2$, the map $T^*X\to \c_X^*$ is a nondegenerate quadratic form on every fiber of $T^*X$ over $X$, and $\g_X^{*,\rs}$ is equal to the image of $T^*X^\rs$ in $\g^*$.
\end{lemma}

\begin{proof}
 Let $x\in X$ with stabilizer $H$; the fiber of $T^*X$ over $x$ is canonically identified with $\mathfrak h^\perp$. Since $H$ is reductive, there is a nondegenerate invariant symmetric bilinear form on $\g^*$ with nondegenerate restriction to $\mathfrak h^\perp$. Indeed, recall that every invariant symmetric bilinear form on a simple (non-abelian) Lie algebra is a multiple of the Killing form. Take an invariant, nondegenerate extension $Q$ of the Killing form on the semisimple part of $\mathfrak g$; it is a standard fact that $Q$ restricts to a nondegenerate form on any Cartan subalgebra. By invariance, the restriction of $Q$ to $\mathfrak h$ is nondegenerate; therefore, it is also nondegenerate on $\mathfrak h^\perp$, when we use the form to identify $\g$ with $\g^*$.
 The corresponding quadratic form
 $$ \mathfrak h^\perp \to \mathbbm A^1$$ is $H$-invariant. It thus has to factor through the invariant-theoretic quotient 
 $$\mathfrak h^\perp\sslash H = T^*X\sslash G = \c_X^*,$$
 which therefore has to be nondegenerate.
 
 Notice that the image of $\a_X^*$ in $\a^*\sslash W$ is birational to the quotient of $\a_X^*$ by its normalizer in $W$. Since we are in rank one, and $W_X$ is non-trivial, this normalizer acts on $\a_X^* = \mathbbm A^1$ in the same way as $W_X$, namely by $\pm 1$, therefore the image of $\a_X^*$ in $\a^*\sslash W$ is birational to $\c_X^*$. (They are in fact equal, since the Killing form provides the quadratic generator of $F[\c_X^*]$, but we won't need that.) In particular, the $\Gm$-orbit $\mathring c_X^*$ embeds into $\a^*\sslash W$.

 By \cite[Satz 5.4]{KnWeyl}, the closure $\overline{\mu(T^*X)}$ of the image of the moment map is equal to the set of $G$-translates of $(\mathfrak l_1+\u_P)^\perp$, where $P$ is a parabolic in the class of $P(X)$, and $L_1 U_P$ is the kernel of the map $P\to A_X$. Choosing a linear section $\sigma$ of the natural quotient map $(\mathfrak l_1+\u_P)^\perp\to \a_X^*$, the subset 
 $$\sigma(\mathring \a_X^*)\cdot G \subset \overline{\mu(T^*X)} = (\mathfrak l_1+\u_P)^\perp\cdot G$$ 
 is open, coincides with the preimage of $\mathring \c_X^* \subset \a^*\sslash W$ in $\mu(T^*X)$, and is smooth because the action map 
$$  \sigma(\mathring\a_X^*)\times^L G \to \sigma(\mathring \a_X^*)\cdot G,$$
 where $L$ is the centralizer of $\sigma(\a_X^*)$, is \'etale. 
 
 Thus, the factorization \eqref{momentfact} of the moment map restricts to 
 $$ T^*X^\rs \to \g_X^{*,\rs} \to \sigma(\mathring \a_X^*)\cdot G\subset \g^*.$$
 By construction, the map $\g_X^*\to \g^*$ is finite, thus the same is true for its restriction over the subset $\sigma(\mathring \a_X^*)\cdot G$. If $\overline{\g_X^*}$, $\overline{\c_X^*}$ denote the images of ${\g_X^*}$, ${\c_X^*}$ under \eqref{momentfact}, it follows (in arbitrary rank) from \cite[Satz 6.4]{KnWeyl} that $\g_X^*$ is birational to the fiber product $\c_X^* \times_{\overline{\c_X^*}}\overline{\g_X^*}$. (This proposition identifies $\g_X^*$ with the invariant-theoretic quotient by $W_X$ of a variety which is birational to $  \sigma(\mathring\a_X^*)\times^L G$.) Again, $\overline{\c_X^*}$ is birational (actually equal) to $\c_X^*$, so the map $\g_X^*\to \overline{\g_X^*}$ is birational. Hence, its restriction over $\mathring \c_X^*$ is a finite, birational map of normal varieties, therefore an isomorphism.
\end{proof}

\begin{proposition}\label{rsinKnop}
 For $X$ of rank one with $W_X=\ZZ/2$,  the restrictions of Knop's sections $\hat\kappa_X$ and $\tilde\kappa_X$ to $\mathring\a_X^*$ are isomorphisms onto the subsets of regular semisimple vectors on $\widehat{T^*X}^\bullet$, resp.\ $\widetilde{T^*X}^\bullet$.
\end{proposition}

\begin{proof}
Knop's section $\tilde \kappa_X$ is an embedding, and $\hat\kappa_X$ is an embedding over $\mathring a_X^*$, so it is enough to prove surjectivity. 
We have a dominant, proper map $\widetilde{T^*X}^\bullet\to \widehat{T^*X}^\bullet$, and the image of $\tilde\kappa_X$ surjects onto the image of $\hat\kappa_X$, so it is enough to prove the proposition for $\widetilde{T^*X}$. 

The $B$-orbits of non-maximal rank in $X$ have rank zero, and therefore any cotangent vector over such an orbit which is perpendicular to $\mathfrak u_B$ maps to $0\in \a^*$ (because the stabilizer of a point, modulo $U_B$, is equal to $B/U_B$). Therefore, the regular semisimple elements $(v, B)\in \widetilde{T^*X}^\rs$ all live over $B$-orbits of maximal rank, that is, if $x$ is the image of $v$ in $X$, then $x\cdot B$ is a $B$-orbit of maximal rank. It suffices to show that those which belong to $\widetilde{T^*X}^\bullet$ live over the open orbit. 

If not, that is, if there is a regular semisimple vector $(v,B)\in \widetilde{T^*X}^\bullet$ which lives over an orbit $Y$ of maximal rank other than the open one, hence also belongs to the irreducible component of $\widetilde{T^*X}^Y\subset \widetilde{T^*X}$ indexed by $Y$, that means that the intersection of two distinct irreducible components 
$$ \widetilde{T^*X}^\bullet \cap \widetilde{T^*X}^Y$$
contains regular semisimple vectors. In particular, the same holds for the intersection of the corresponding irreducible components of $\widehat{T^*X}$, $$ \widehat{T^*X}^\bullet \cap \widehat{T^*X}^Y.$$

I claim that, in rank one with $W_X=\ZZ/2$, the map $\widehat{T^*X}\to T^*X$ is \'etale over $T^*X^\rs$; indeed, the image of $T^*X^\rs$ is the subset $\mathring\c^*_X$ of $\a^*\sslash W$, and the normalizer of $\a_X^*$ in $W$ has to coincide with $W_{L(X)}\rtimes W_X$, because $W_X=\ZZ/2$ is the group of automorphisms finite order of the lattice $\Lambda_X\subset \a_X^*\simeq\mathbbm A^1$, and $W_{L(X)}$ is its centralizer. The distinct $W$-conjugates of $\mathring\a_X^*$ have empty intersections (because these are distinct one-dimensional vector subspaces of $\a^*$ with their origins removed), thus we have
$$ \widehat{T^*X}^\rs = \bigsqcup_{w\in W_{L(X)}\rtimes W_X\backslash W} T^*X^\rs\times_{\mathring\c^*_X} (\mathring\a_X^*)^w$$
as algebraic varieties, where $w$ denotes a representative for the given coset, and $(\mathring\a_X^*)^w$ is the $w$-conjugate of $\mathring\a_X^*$ inside of $\a^*$. Hence, the map $\widehat{T^*X}^\rs\to T^*X^\rs$ is the base change of the \'etale maps $(\mathring\a_X^*)^w\to \mathring\c_X^*$, hence \'etale.

But this implies that the components of $\widehat{T^*X}^\rs$ have empty intersections.
\end{proof}

\begin{corollary}\label{corollaryGtransitive}
 In the setting of the previous proposition, $G$ acts transitively on every fiber of $\widetilde{T^*X}^{\bullet,\rs}$, $\widehat{T^*X}^{\bullet,\rs}$ or $T^*X^\rs$ over, respectively, $\mathring\a^*_X$, $\mathring\a^*_X$ or $\mathring\c^*_X$. The stabilizer of any point on $\widehat{T^*X}^{\bullet,\rs}$ is conjugate to the kernel $L_1$ of the canonical map $L(X)\to A_X$, where $L(X)$ is a Levi subgroup of $P(X)$; the stabilizer of any point on $\widetilde{T^*X}^{\bullet,\rs}$ is conjugate to a Borel subgroup of such an $L_1$; finally, the stabilizer of any point on $T^*X^\rs$ is conjugate over the algebraic closure to such a subgroup $L_1$.
\end{corollary}

\begin{proof}
 It is enough to prove transitivity of the action for $\widetilde{T^*X}^{\bullet,\rs}$. By Proposition \ref{rsinKnop}, Knop's section is an isomorphism onto regular semisimple vectors:
$$ (X\times \B)^0 \times \mathring\a_X^* \xrightarrow\sim \widetilde{T^*X}^{\bullet,\rs}.$$
The group $G$ acts transitively on $(X\times \B)^0$, hence on the fiber over any point in $\mathring\a_X^*$.

The statement on stabilizers, in the first two cases, follows by examining stabilizers on the domains of Knop's sections. In the last case, it follows from the fact that the preimage of a $G$-orbit in the $W_X$-torsor $\widehat{T^*X}^{\bullet,\rs}\to T^*X^\rs$ consists, over the algebraic closure, of $W_X$-many $G$-orbits.
\end{proof}

In particular, considering the map $T^*X\to \g_X^*$ and setting 
$$\hat\g_X^*=\g_X^*\times_{\c_X^*} \a_X^*,$$
so that we have a map
$$ \widehat{T^*X}^\bullet\to \hat\g_X^*,$$
the $G$-stabilizer of any element $\hat Z\in \hat\g_X^{*,\rs}$, resp.\ $Z\in \g_X^{*,\rs}$, acts transitively on its fiber in $\widehat{T^*X}^\bullet$, resp.\ $T^*X$. The $G$-stabilizer of such an element is a Levi subgroup of $G$ over the algebraic closure and, in the case of $\hat\g_X^{*,\rs}$, a Levi subgroup $L$ equipped with a choice of parabolic $P$ in the class of $P(X)$. (Indeed, this is the parabolic in the class $P(X)$ containing $B$, for any lift $(Z, B)$ of $\hat Z$ to $\tilde\g_X^*$ --- equivalently, the parabolic $P\in \B_X$ in the domain of Knop's section.) Corollary \ref{corollaryGtransitive} implies that $L$ acts on the fiber of $\hat Z$ precisely through the quotient $L\to P\to A_X$. This is a special case of \cite[Proposition 2.4]{KnMotion}, giving rise to a canonical action of $A_X$ on $\widehat{T^*X}^{\bullet,\rs}$, which commutes with the action of $G$: 

\begin{equation}\label{AXaction}
 A_X\times \widehat{T^*X}^{\bullet,\rs}\to \widehat{T^*X}^{\bullet,\rs}.
\end{equation}

Corollary \ref{corollaryGtransitive} implies:

\begin{corollary}\label{polarizedtorsor}
 In the setting of Proposition \ref{rsinKnop}, $\widehat{T^*X}^{\bullet,\rs}$ is canonically an $A_X$-torsor over $\hat\g_X^{*,\rs}$; namely, the centralizer $L$ of any point on $\hat\g_X^{*,\rs}$ acts transitively on its fiber via the canonical quotient $L\to A_X$.
\end{corollary}

In Section \ref{sec:resolution}, we will see (following Knop, again) how to formulate the analog of this for the map $T^*X^\rs\to \g_X^{*,\rs}$, and to extend it to an action of a group scheme over the whole space $\g_X^*$.

\subsection{Borel orbits over the algebraic closure}

From now, $X$ is always affine homogeneous spherical $G$-variety of rank one, with $W_X=\ZZ/2$. Its weight lattice $\Lambda_X$ (the character group of $A_X$) is thus isomorphic to $\ZZ$. 

In the present subsection, we work over $\bar F$, the algebraic closure of $F$. 

Recall, again, that the rank of a $B$-orbit on $X$ is the rank of the torus $B_x\backslash B/N$, where $B_x$ is the stabilizer of a point $x$ on the orbit, and $N\subset B$ the unipotent radical. For what follows, for a Borel subgroup $B$ and a simple positive root $\alpha$, we denote by $P_\alpha$ the parabolic generated by $B$ and the simple root space associated to the root $-\alpha$, and by $\mathcal R(P_\alpha)$ its radical (so that $P_\alpha/\mathcal R(P_\alpha)\simeq \PGL_2$). (We also use the notation $P_{\alpha\beta}$ for the parabolic generated by $B$ and the negative root spaces associated to two simple roots $\alpha, \beta$, etc.)

Knop has defined in \cite{KnOrbits} a rank-preserving action of the Weyl group of $G$ on the set of Borel orbits (or, equivalently, the set of $B$-orbit closures), which is transitive on the subset of orbits of maximal rank, which includes the open orbit. For the reflection $w_\alpha$ associated to a simple root $\alpha$, and a $B$-orbit $Y$, one considers the spherical $\PGL_2$-variety $YP_\alpha/\mathcal R(P_\alpha)$ which is of one of the following four types:
 \begin{equation}\label{fourtypes}
  \parbox{30em}{
\begin{enumerate}
 \item $\PGL_2\backslash \PGL_2$, i.e., a point;
  \item $T\backslash \PGL_2$, where $T$ is a torus;
  \item $\mathcal N(T)\backslash \PGL_2$, where $\mathcal N(T)$ is the normalizer of a torus;
  \item $S\backslash \PGL_2$, with $N_2\subset S \subset B_2$, where $B_2\supset N_2$ denote a Borel subgroup of $\PGL_2$, and its unipotent radical.
\end{enumerate}
  }
\end{equation}
In the first three cases, there is a single orbit of largest rank in $YP_\alpha$, and it is fixed by $w_\alpha$. In the last case, there are two such orbits, say $Y$ and $Z$, and $w_\alpha$ interchanges them; moreover, for their character groups $\Lambda_Y=\Hom(A_Y,\Gm)$, $\Lambda_Z=\Hom(A_Z,\Gm)$, we have: 
\begin{equation}\label{charsKnopaction}
 \Lambda_Z = \Lambda_Y^{w_\alpha}
\end{equation}
inside of $\Hom(A,\Gm)$.

Since, in our case, $X$ is of rank one, all $B$-orbits are of rank one or zero.

Following Brion \cite{BrOrbits}, we have by \cite[\S 3.1]{SV}:

\begin{lemma}\label{Briongenerators}
 There is a $B$-orbit $Z$ of rank one, and a simple root $\alpha$, or two orthogonal simple roots $\alpha,\beta$, such that, setting $P=P_\alpha$, resp.\ $P=P_{\alpha\beta}$, the spherical variety $ZP/\mathcal R(P)$ is isomorphic to one of the following:
 \begin{enumerate}
  \item $T\backslash \PGL_2$, where $T$ is a torus;
  \item $\mathcal N(T)\backslash \PGL_2$, where $\mathcal N(T)$ is the normalizer of a torus;
  \item $\PGL_2$ as a $\PGL_2^2$-space (when $P=P_{\alpha\beta}$). 
 \end{enumerate}
Moreover, these possibilities are mutually exclusive, as in the first case the weight lattice $\Lambda_X$ is spanned by a root $\gamma$ of $G$, in the second case it is spanned by the double $2\gamma$ of a root, and in the third case it contains the sum $\gamma=\gamma_1+\gamma_2$ of two strongly orthogonal roots.
\end{lemma}

The weight $\gamma$ of the lemma above called the (normalized) \emph{spherical root} of $X$, by \cite[\S 3.1]{SV}.
Correspondingly to the three cases, we will say that the spherical root is of type $T$, $N$ or $G$. 

I caution the reader that this is not the standard normalization of spherical roots in the theory of spherical varieties (e.g., as in \cite{Luna}), and it also differs from yet another normalization that appears in \cite{KnAut}; however, in this paper I will call $\gamma$ the spherical root.

\begin{corollary} \label{allrankone}
 The spherical root is of type $G$ if and only if all $B$-orbits are of rank one.
\end{corollary}

\begin{proof}
 Indeed, the spherical root were of type $G$, but there were a $B$-orbit $Z_0$ of rank zero, since the parabolics of type $P_\alpha$ generate the group, and $ZG= X$, there is a sequence of spherical roots $\alpha_1, \dots, \alpha_m$ such that $Z_i:= Z_{i-1} P_{\alpha_i}$ is of dimension $>\dim Z_{i-1}$, and $\overline{Z_m}=X$. By the analysis of the four types of \eqref{fourtypes}, the open $B$-orbit on $Z_i$ either has rank zero or one, and since $\mathring X$ has rank one, there is a $B$-orbit $Z$ of rank zero, and a simple root $\alpha$, such that $ZP_\alpha$ contains a $B$-orbit of rank one. But this can happen only if $ZP_\alpha/\mathcal R(P_\alpha)$ is isomorphic to $T\backslash \PGL_2$ or $\mathcal N(T)\backslash \PGL_2$, which is impossible by the above lemma.
 
 Vice versa, a spherical root of type $T$ or $N$ requires, by definition (Lemma \ref{Briongenerators}), that there be an orbit $Z$ of rank one and a parabolic of type $P_\alpha$ such that $ZP_\alpha/\mathcal R(P_\alpha)$ is isomorphic to $T\backslash \PGL_2$ or $\mathcal N(T)\backslash \PGL_2$; in both of those cases, the closed orbit in $ZP_\alpha$ has smaller rank than $Z$.
\end{proof}

We now study closed $B$-orbits, first over the algebraic closure:

\begin{lemma}
 Let $Z\subset X$ be a $B$-orbit, and $H\subset G$ the stabilizer of a point on $Z$. The following are equivalent:
 \begin{enumerate}
  \item $H\cap B$ is a Borel subgroup of $H$;
  \item $Z$ is closed.
 \end{enumerate}
\end{lemma}

Here, for non-connected groups, by slight abuse of language we use ``Borel'' for any solvable subgroup such that the quotient is projective, whether it is connected or not.

\begin{proof}
$B$-orbits on $H\backslash G$ are in natural bijection with $H$-orbits on the flag variety $G/B$. Since the latter is projective, a closed $H$-orbit on $G/B$ is a projective homogeneous $H$-variety with solvable stabilizers, hence of the form $H/B_H$ for a Borel subgroup of $H$; vice versa, if it is of this form, it is projective and hence closed.
\end{proof}

\begin{lemma}\label{lemmaclosedorbits}
Let $Z\subset X$ be a closed $B$-orbit. Then, one of the following two holds:
\begin{enumerate}
 \item $Z$ is of rank zero, and for every simple root $\alpha$ such that $Y:=ZP_\alpha\ne Z$, we have $Y/\mathcal R(P_\alpha) \simeq T\backslash \PGL_2$ or $\mathcal N(T)\backslash \PGL_2$ (notation as above);
 \item or, $Z$ is of rank one, and for all simple roots $\alpha$ we have $ZP_\alpha=Z$, except for two orthogonal simple roots $\alpha, \beta$ for which, setting $Y:=ZP_{\alpha\beta}$, we have $Y/\mathcal R(P_{\alpha\beta}) \simeq 
\PGL_2^\diag\backslash \PGL_2^2$.
\end{enumerate}
\end{lemma}

\begin{proof}
 In the first case, we only need to exclude the possibility that $ZP_\alpha/\mathcal R(P_\alpha)=S\backslash \PGL_2$ with $N_2\subset S\subset B_2$. Let $H\subset G$ be the stabilizer of a point on $Z$. Since $Z$ is of rank zero, $H\cap B$ and $G\cap B$ have the same rank and, in particular, contain a common maximal torus $T$. Decomposing the Lie algebras of $G$ and $H$ into $T$-eigenspaces, we see that each root space for $H$ is also a root space for $G$. This means that the statement $ZP_\alpha/\mathcal R(P_\alpha)=S\backslash \PGL_2$ lifts to the statement that $\mathfrak h\cap \mathfrak p_\alpha$ contains the root space $\mathfrak n_\alpha$ corresponding to the root $\alpha$ (but not its opposite). Thus, if $L_\alpha'$ denotes the commutator of the standard Levi (with respect to $T$) of $P_\alpha$, then $Y=ZP_\alpha$ contains the space $(T\cap L_\alpha')N_\alpha\backslash L_\alpha'$ as a subvariety. But this is nontrivial projective, a contradiction, since $X$ is assumed affine. 
 
 In the second case, observe first that the spherical root is necessarily of type $G$. Indeed, if $H$ is the stabilizer of a point in the closed orbit, it follows from the previous lemma that $\rk(H)=\rk(G)-1$. Since $X$ is homogeneous, this holds for all stabilizers, and therefore there cannot be a Borel orbit of rank zero, so by Corollary \ref{allrankone}, the spherical root is of type $G$.
  
 I will now rely on the classification of Table \eqref{thetable}, since I currently do not have a proof which avoids any kind of classification. Since all $B$-orbits are of rank one, they form a partially ordered set (by dimension) which can be identified, using Knop's action, with the homogeneous set $(W_{L(X)}\rtimes W_X)\backslash W$ for the Weyl group, with the minimal length of a representative of a coset corresponding to the codimension of the orbit. The little Weyl group $W_X=\ZZ/2$, for spherical roots of type $G$ (i.e., the cases $D_n$--$B_3''$ of Table \ref{thetable}), is generated by the element $w_\gamma=w_{\gamma_1} w_{\gamma_2}$, where $\gamma_1,\gamma_2$ are the two strongly orthogonal roots such that $\gamma=\gamma_1+\gamma_2$. These partially ordered sets, together with the graph of Knop's action, have been depicted in \cite[6.19, 6.16]{SaSph}, and one observes that there is actually a unique minimal element (Borel orbit) in this partially ordered set, and two unique, mutually orthogonal simple roots raising it to the same Borel orbit. 

\end{proof}

Notice that a subset $Y\subset X$ as in the lemma is closed, since this is the case for $Z$ and the action map $X\times^B P\to X$ is proper; thus, the image of the closed subset $Z\times^B P$ is closed. 

\begin{definition}
\label{boppair}
A pair $(Y,P)$ as in Lemma \ref{lemmaclosedorbits}, where $P=P_\alpha$ in the first case and $P=P_{\alpha\beta}$ in the second, will be called \emph{a basic orbit-parabolic pair}. In other words, a basic orbit-parabolic pair consists of a parabolic $P$ of type $P_\alpha$ or $P_{\alpha\beta}$ (where $(\alpha,\beta)$  is a strongly orthogonal pair of roots), and a closed $P$-homogeneous subvariety $Y\subset X$ such that $Y/\mathcal R(P)$ is isomorphic to $T\backslash \PGL_2$ or $\mathcal N(T)\backslash \PGL_2$ or, respectively, $\PGL_2^\diag\backslash\PGL_2^2$. 
\end{definition}

Such a pair will play an important role in various arguments in this paper, since by the above lemma it allows us to reduce many arguments to the basic rank-one cases, labelled $A_1$ and $D_2$ in Table \eqref{thetable}. In reality, the precise choice of parabolic will never matter; only its class matters, and if $\B_P$ denotes the flag variety of parabolics in this class, $Y$ can be replaced by the $G$-orbit on $X\times \B_P$, whose fiber over $P$ is $Y$.

\begin{example}
\label{examplebopp}
 Consider the space $X=\SO(M)\backslash \SO(M \oplus F)$, where $M$ is a non-degenerate quadratic space, and $M\oplus F$ denotes a non-degenerate quadratic space of one dimension larger. If $n$ is even, the stabilizer $P$ of an isotropic flag $M_1 \subset M_2\subset \dots \subset M_{\frac{n}{2}-1}=:M'$ of dimension $\frac{n}{2}-1$ in $M$ is a parabolic of type $P_\alpha$, and any of the two points on the hyperboloid $X$ which are represented by vectors perpendicular to $M$ is a closed $P$-orbit $Y$. The quotient $Y/\mathcal R(P)$ is isomorphic to the hyperboloid $\SO_2\backslash \SO_3$, where $\SO_3$ stands for $\SO(M'^\perp/M')$ (where the orthogonal complement is taken in $M\oplus F$).
 
 If $n$ is odd, we can similarly choose an isotropic flag $M_1 \subset M_2\subset \dots \subset M_{\frac{n-3}{2}}=:M'$ of dimension $\frac{n-3}{2}$ in $M$, and this defines a parabolic $P$ of type $P_{\alpha\beta}$. Again, the points on the hyperboloid $X$ perpendicular to $M$ represent closed $P$-orbits $Y$, and the quotient $Y/\mathcal R(P)$ is isomorphic to $\SO_3\backslash\SO_4$.
 \end{example}

In this paper we do not consider arbitrary homogeneous spherical varieties of rank one. If we divide affine homogeneous spherical varieties of rank one into equivalence classes, with $X\sim X'$ if their quotients by the groups of their $G$-automorphisms are isomorphic, then it turns out that only one representative in each equivalence class is appropriate for the relative trace formula comparison that we are performing. This representative is the one listed in Table \eqref{thetable}, and is described by the following:

\begin{proposition}\label{correctrepresentative}
 If the spherical root $\gamma$ is of type $N$ --- equivalently, if $\Lambda_X$ is spanned by $2\gamma$, then there is an equivariant two-fold cover $X'\to X$ with $\Lambda_{X'}$ spanned by $\gamma$. 
 
 If the spherical root $\gamma$ is of type $G$, then there is an equivariant finite cover $X'\to X$ (possibly $X'=X$) such that $\Lambda_{X'}$ is spanned by $\frac{\gamma}{2}$. Moreover, $\Aut^G(X')=\mathbb Z/2$.
 
 Moreover, in both cases, the stabilizers of points on $X'$ are connected.
\end{proposition}

\begin{proof}
The first statement is \cite[Lemme 6.4.1]{Luna}.

For the second, if $X=H\backslash G$, replace $G$ by the simply connected cover of its derived group; it necessarily acts transitively on $X$, because if the connected center of $G$ was not acting trivially, the rank of $X$ would be greater than one. 
So, we can without loss of generality denote that by $G$. I claim that $\frac{\gamma}{2} \in \Hom(A,\Gm)$. To show this, it is enough to show that its pairing with every simple coroot is an integer. Without loss of generality, we may replace $\gamma$ by its Weyl group conjugate $\gamma'=\alpha-\beta$ which belongs to the character group of a closed $B$-orbit $Z$, where $\alpha, \beta$ are two orthogonal simple roots, as in Lemma \ref{lemmaclosedorbits}. Clearly, the pairing of $\frac{\gamma'}{2}$ with $\check\alpha, \check\beta$ is integral. On the other hand, by the same lemma, for every simple root $\delta\ne \alpha, \beta$, we have $YP_\delta = Y$; in that case, the image of the cocharacter $\check\delta$ in $A=B/N$ belongs to the image of the stabilizer of a point on $Y$, hence $\left<\gamma', \check\delta\right>=0$, proving the claim.  The result on the existence of $X'$ now follows as a simple special case of \cite[Th\'eor\`eme 2]{Luna}; it also follows from the same that there exists an $X_0$ whose weight lattice is spanned by $\gamma$. The cover $X'\to X_0$ gives rise to a nontrivial involution $\tau$ of $X'$, hence to an embedding $\ZZ/2\hookrightarrow \Aut^G(X')$. (This embedding can also be established as a special case of \cite[Theorem 1.2]{KnAut}.) On the other hand, we also have an embedding 
$$\iota:\Aut^G(X') \hookrightarrow \Hom(\Lambda_{X'}, \bar F^\times),$$ by considering the action of an automorphism on the lines of $B$-eigenvectors in the function field $\bar F(X')$. To show that $\ZZ/2$ is the whole automorphism group, it suffices to show that $\iota(\tau)$ is trivial on the (normalized) spherical root $\gamma$. 

The set of $B$-orbits of maximal rank is acted upon transitively by the Weyl group action of Knop.
Any $G$-automorphism $\sigma$ of $X'$ preserves the open Borel orbits and commutes with 
Knop's action on the set of Borel orbits, hence preserves all orbits of maximal rank. If $Y, Z$ are two $B$-orbits of maximal rank, and $\alpha$ is a simple root such that $Z\cdot P_\alpha  = Z\cup Y$, so that $Z^{w_\alpha} = Y$ under Knop's action, the quotient $Z\cdot P_\alpha/U_{P_\alpha}$ is a quotient of $N\backslash P_\alpha$ (where $N$ is the unipotent radical of $B$), and $\sigma$ descends to an automorphism of $Z\cdot P_\alpha/U_{P_\alpha}$, which is induced by the ``left'' action of an element of $A=B/N$. On the character groups of the two $B$-orbits, therefore, which are related by \eqref{charsKnopaction}, $\sigma$ induces conjugate automorphisms, i.e., if $\iota_Z$ denotes the analogs of the map $\iota$:
$$ \iota_Z: \Aut^G(X') \to \Hom(\Lambda_Z, \bar F^\times),$$
then those are compatible with the $W$-action: 
$$\iota_{Z^w}(\sigma) = \iota_Z(\sigma)^w.$$

Let $(Y',P_{\alpha\beta})$, now, be a basic orbit-parabolic pair for $X'$, as in Lemma \ref{lemmaclosedorbits}. Since $Y'$ is preserved by $G$-automorphisms, a $G$-automorphism $\sigma$ of $X'$ restricts to a  $P_{\alpha\beta}$-automorphism of $Y'$. The quotient $Y'/\mathcal R(P_{\alpha\beta})$ is isomorphic to $\PGL_2$, and the weight lattice of its closed $B$-orbit is spanned by the weight $\alpha-\beta$. Since $\PGL_2$ has no $\PGL_2\times\PGL_2$-automorphisms, we see that $\iota_{Y'}(\sigma)$ has to be trivial on $\alpha-\beta$, and therefore $\iota(\sigma)$ has to be trivial on its $W$-conjugate $\gamma$.

If the stabilizers were not connected, then there would exist a further cover $X''\to X'$, giving an embedding of character groups of rank one with nontrivial cokernel: $\Lambda_{X'}\subset \Lambda_{X''}$. But this is impossible in both cases for the spherical root $\gamma$ of $X'$: If it is of type $T$, the same should hold for $X''$, and their character groups are both spanned by $\gamma$; if it is of type $G$, the same should hold for $X''$, and the character group in such a case is spanned by either $\gamma$ or $\frac{\gamma}{2}$ (in this case, by the latter), because no smaller fraction of $\gamma$ can belong to the weight lattice of $A$. Thus, stabilizers are connected.

\end{proof}

\begin{definition}\label{defcorrectrep}
The variety $X'$ of Proposition \ref{correctrepresentative} will be called the ``correct representative'' of its equivalence class modulo $G$-automorphisms.
\end{definition}

\subsection{Borel orbits and the moment map over $F$} \label{ssrationality}

We now return to our non-algebraically closed field $F$ in characteristic zero, with $G$ split over $F$. We maintain the assumptions of the previous subsection for $X$, and, moreover, we assume that it is the correct representative given by Proposition \ref{correctrepresentative}.

Under these assumptions, our main goal for the rest of this section is to prove the following:

\begin{proposition}\label{rationality}
 The following are equivalent:
\begin{enumerate}
 \item The stabilizer $H$ of one, equivalently any, $F$-point on $X$ is a split reductive group.
 \item One, equivalently every, closed $B_{\bar F}$-orbit on $X_{\bar F}$ is defined over $F$.
 \item The invariant moment map $T^*X\to \c_X^*$, viewed as a quadratic form on the fibers of $T^*X$, is split (maximally isotropic) on one, equivalently every, fiber.
\end{enumerate}
\end{proposition}
Notice that $H$ is connected, by the fact that we are working with correct representatives (Definition \ref{defcorrectrep}), and Proposition \ref{correctrepresentative}. The proposition is not true without the assumption on ``correct representatives'', e.g., for the variety $\mathcal N(T)\backslash \PGL_2$.

We begin with some preliminary lemmas and constructions.
If $B\subset G$ is a Borel subgroup, since $B$ is split solvable, by standard Galois cohomology, every $B$-orbit which is defined over $F$ has an $F$-point. This holds, in particular, for the open $B$-orbit, which is unique hence defined over $F$, therefore $X$ has $F$-points.

\begin{lemma}\label{maxrankdefined}
 All $B$-orbits of maximal rank (over $\bar F$) are defined over $F$. Moreover, if $Y$ is a $B$-orbit of maximal rank, then $Y(F)$ meets any $G(F)$-orbit on $X(F)$ nontrivially.
\end{lemma}

\begin{proof}
The open $B$-orbit is defined over $F$ and has rank one. By the definition of Knop's action of the Weyl group on the set of Borel orbits, and the fact that $G$ is split, one sees from the definition that the action is defined over $F$, therefore every $F$-orbit of maximal rank is defined over $F$. Any $G(F)$-orbit is open in $X(F)$ (in the Hausdorff topology induced by the topology of $F$), and therefore has to contain $F$-points of the open $B$-orbit. Finally, we use Knop's action to deduce the same result for any $B$-orbit of maximal rank: if $Y^{w_\alpha} = Z$, with $Z$ open in $YP_\alpha$, and $Z(F)$ contains a point $z$ in a given $G(F)$-orbit, then there is a $g\in P_\alpha(F)$ such that the $P_\alpha(F)$-stabilizer of $zg$ is contained in $B(F)$ (by the fact that the map $P_\alpha(F)\to B\backslash P_\alpha(F)$ is surjective), therefore the same $G(F)$-orbit also contains a point of $Y(F)$.
\end{proof}

Now let $(Y,P)$ be a basic orbit-parabolic pair, see Definition \ref{boppair}. 

\begin{lemma}\label{boppdefined}
The variety $Y$ is defined over $F$, and the action map $Y\times^P G\to X$ is surjective on $F$-points, i.e., $Y(F)\times^{P(F)} G(F)\twoheadrightarrow X(F)$.  
\end{lemma}

\begin{proof}
By definition, the variety $Y$ contains a $B$-orbit of maximal rank (rank one). By Lemma \ref{maxrankdefined}, this $B$-orbit is defined over $F$, and its $F$-points meet every $G(F)$-orbit; hence, the map $Y(F)\times G(F)\to X(F)$ is surjective.
\end{proof}

Now we define a proper cover of the cotangent bundle $T^*X$ which is intermediate between this and a component of the cover $\widetilde{T^*X}$ defined in \S \ref{ssmomentmap}. Let $(Y,P)$ be a basic orbit-parabolic pair. Similarly to the definition of $\widetilde{T^*X}$, we let
$$\widetilde{T^*X}^P = \{(v,P')| v\in T^*X, P'\sim P, \mu(v) \in \u_{P'}^\perp\}.$$
Here, $P'\sim P$ means that $P'$ is conjugate to $P$. (Really, $P$ denotes here a class of parabolics, and $Y$ can be thought of as a $G$-orbit on $P\backslash G\times X$; nothing depends on a choice in this class.)

We let $\widetilde{T^*X}^{P,Y} = T_Y^*X^{\u_P^\perp} \times^P G$, where the exponent $~^{\u_P^\perp}$ means that the image under the moment map belongs to $\u_P^\perp$ --- hence, $T_Y^*X^{\u_P^\perp}$ is the pullback to $Y$ of the cotangent bundle of the quotient $Y/U_P$. 

\begin{lemma}\label{partialcover}
 $\widetilde{T^*X}^{P,Y}$ is an irreducible component of $\widetilde{T^*X}^P$ and, therefore, is proper and dominant over $T^*X$.
\end{lemma}

\begin{proof}
 It is irreducible, since $T_Y^*X^{\u_P^\perp}$ is a vector bundle over the irreducible variety $Y$, and a closed subset of $\widetilde{T^*X}^P$. Moreover, considering the natural map $\widetilde{T^*X}\to \widetilde{T^*X}^P$, it is clear from the definition that $\widetilde{T^*X}^{P,Y}$ contains the image of the irreducible component of maximal dimension corresponding to the open $B$-orbit in $Y$; therefore, it is an irreducible component of $\widetilde{T^*X}^P$, and dominant over $T^*X$.
\end{proof}
Now let $Y_2=Y/U_P$. Since $Y$ is homogeneous under $P$, this is homogeneous under the Levi quotient $L$ of $P$, and its quotient by $\mathcal Z(L)$ (the center of $L$) is the rank-one spherical variety $Y/\mathcal R(P)$ described in Lemma \ref{lemmaclosedorbits}, hence isomorphic to $T\backslash \PGL_2$ or $\PGL_2$. There are natural quotient maps of (the total spaces of) vector bundles
$$ T_Y^*X^{\u_P^\perp} \to T^*Y^{\u_P^\perp} \to T^* Y_2,$$
with the former an isomorphism on the base and the latter an isomorphism on the fiber.

\begin{lemma}\label{lemmaY2}
 The connected center of $L$ acts trivially on $Y_2$, and $Y_2$ is isomorphic either to $T\backslash \PGL_2$, where $T$ is a torus, or to $\SL_2$. Moreover, we have $T^* Y_2\sslash L \xrightarrow\sim \c_X^*$, fitting into a natural commutative diagram:
 \begin{equation}\label{moddiagram}
  \xymatrix{
  \widetilde{T^*X}^{P,Y} \ar[r] \ar[d] & T^*X\ar[d] \\
  T^* Y_2 \times^P G \ar[r] & \c_X^*.
  }
 \end{equation}
\end{lemma}

\begin{proof}
 The variety $Y_{2,\ad}:= Y_2/\mathcal Z(L) = Y/\mathcal R(P)$ is, by Lemma \ref{lemmaclosedorbits}, isomorphic to $T\backslash \PGL_2$ or $\PGL_2$. If the identity component of $\mathcal Z(L)$ was acting nontrivially on $Y_2$, the rank of $Y_2$, hence of $Y$, would be larger than one, a contradiction. Thus, $Y_2 \to Y_{2,\ad}$ is a finite cover. Since $X$ is taken to be the correct representative in its class, the character group $\Lambda_X$ is generated by the spherical root $\gamma$, if that is of type $T$, or by $\frac{\gamma}{2}$, if it is of type $G$. Correspondingly, by \eqref{charsKnopaction}, the character group of $Y$ (equivalently, of $Y_2$) is generated by $\alpha$, in the first case, and by $\frac{\alpha+\beta}{2}$, in the second, in the notation of Lemma \ref{lemmaclosedorbits}. But this means that in the first case $Y_2=Y_{2,\ad}$, while in the second it has to be a two-fold connected cover of it, hence isomorphic to $\SL_2$. 
 
 We have $A_{Y_2} = A_Y = A_X^w$ for some element $w$ of the Weyl group of $G$, and $W_{Y_2} = W_X = \ZZ/2$. This gives rise to a canonical isomorphism
 $$ \c_{Y_2}^* = \c_X^*,$$
 and the construction of the invariant moment maps $T^* Y_2\to \c_{Y_2}^*$ and $T^*X\to \c_X^*$ is clearly compatible with this isomorphism, showing the commutativity of the diagram.
\end{proof}

\begin{proposition}\label{propVP}
 Let $(Y,P)$ be a basic orbit-parabolic pair; then the map $\widetilde{T^*X}^{P,Y} \to T^*X$ is surjective on $F$-points.
 
 Let $y\in Y$ with image $y_2\in Y_2$, and let $V\supset V_P$ and $V_2$ be, respectively, the fibers of $T^*X$, $\widetilde{T^*X}^{P,Y}$, and $T^*Y_2$ over $y$, $(y,P)$, and $y_2$, respectively. The kernel of the map $V_P\to V_2$ is an isotropic subspace of $V$ (with respect to the quadratic map $V\to \c_X^*$) of dimension 
\begin{equation}\label{dimensions}
 \dim\ker(V_P\to V_2) = \frac{\dim V - \dim V_2}{2}.
\end{equation}

 The quadratic space $V$ is split\footnote{We will be using ``split'' in a slightly non-standard way: for quadratic spaces of odd dimension $d$, ``split'' will mean maximally isotropic, i.e., containing a $\frac{d-1}{2}$-dimensional isotropic subspace.}  (maximally isotropic) if and only if $V_2$ is.
\end{proposition}

\begin{example}
 In the setting of Example \ref{examplebopp} (in particular, $G=\SO(M\oplus F)$ for a quadratic space $M$, $X = \SO(M)\backslash \SO(M\oplus F)$, and we have a basic orbit-parabolic pair $(Y,P)$ as described in that example), for the point $y$ represented by a vector perpendicular to $M$ in the hyperboloid $X$, the fiber of $T^*X$ can be identified with the dual of $M$, or with $M$ itself (using the quadratic form). The space $V_P$ is the orthogonal complement of the isotropic flag $M'$ described in that example, and the space $V_2$ is the quotient $M'^\perp/M'$. If $\dim M$ is even, the group $G$ can be split even if $M$ is not, and this happens precisely when the two-dimensional quadratic space $(M'^\perp\cap M)/M'$ is nonsplit.
\end{example}

Notice that $V_2$ is isomorphic to $\mathbbm A^2$, for spherical roots of type $T$, and to $\mathfrak{sl}_2$, for spherical roots of type $G$, with $P\cap H$ acting through a $1$-dimensional torus quotient in the first case, and through a quotient isomorphic to $\PGL_2$, in the second case.

\begin{proof}
 Let $L$ denote the Levi quotient of the class of parabolics $P$, considered as 
 an abstract algebraic group depending only on the class of $P$, defined uniquely up to conjugacy. We have a canonical map of coadjoint quotients 
 \begin{equation}
  \mathfrak l^*\sslash L \to \g^*\sslash G.
 \end{equation}
 Let 
 $$\widehat{\g^*}^P = \g^* \times_{\g^*\sslash G} \mathfrak l^*\sslash L,$$ 
 $$\widetilde{\g^*}^P = \{(Z,P')| P'\sim P, Z \in \u_{P'}^\perp\},$$ 
 so we have natural, proper maps 
 $$\widetilde{\g^*}^P \to \widehat{\g^*}^P \to \g^*.$$ 
 At the level of $F$-points, the map $\widehat{\g^*}^P \to \g^*$ is not surjective, but the map $\widetilde{\g^*}^P \to \widehat{\g^*}^P$ is. 
 
 Now, considering $\c_X^*$ as a subset of $\g^*\sslash G$ and, similarly, $\c_{Y_2}^*$ as a subset of $ \mathfrak l^*\sslash L$, the canonical isomorphism $\c_{Y_2}^*\simeq \c_X^*$ that we saw in Lemma \ref{lemmaY2} gives rise to a lift:
 \begin{equation}\label{liftcX} \c_X^* \to \mathfrak l^*\sslash L.
 \end{equation}
 The composition of this with the invariant moment map $ T^*X \to \c_X^*$ gives a lift:
 \begin{equation}\label{liftTX} T^*X \to \widehat{T^*X}^P:= T^*X \times_{\g^*} \widehat{\g^*}^P,
 \end{equation}
 which fits into a commutative diagram
 \begin{equation} \xymatrix{  \ar[d]\widetilde{T^*X}^Y \ar[r] &  \widehat{T^*X}^Y \ar[d]\ar[r] & \a^*_Y \ar[d]  \\
 \widetilde{T^*X}^{P,Y} \ar[dr]\ar[r] & \widehat{T^*X}^P \ar[r]\ar@/_/[d] & \c_X^*\\
 & T^*X, \ar@/_/[u]\ar[ur] & }\end{equation}
 where $\widetilde{T^*X}^Y$, $\widehat{T^*X}^Y$ denote the irreducible components corresponding to the open $B$-orbit in $Y$ (see \S \ref{ssmomentmap}). The scheme-theoretic image of $T^*X$ in $\widehat{T^*X}^P$ is the same as the image of $\widehat{T^*X}^Y$, hence an irreducible component $\widehat{T^*X}^{P,Y}$ of $\widehat{T^*X}^P$. The fact that the map  $\widetilde{\g^*}^P \to \widehat{\g^*}^P$, hence its base change $\widetilde{T^*X}^{P,Y} \to \widehat{T^*X}^{P,Y}$, is surjective on $F$-points, proves the surjectivity statement.
 
 Let us now count dimensions. The map $\widetilde{T^*X}^Y\to \widetilde{T^*X}^{P,Y}$ is finite, hence 
 $$\dim \widetilde{T^*X}^{P,Y} = \dim \widetilde{T^*X}^Y = \dim \widehat{T^*X} + \dim \B_{L(X)} = \dim T^*X + \dim \B_{L(X)},$$
 where $\B_{L(X)}$ denotes the full flag variety of the Levi quotient of $P(X)$. In terms of the above vector spaces, setting $P_H= P\cap H$, where $H$ is the stabilizer of $y\in Y$, we have 
 $$ T^*X = V\times^H G,$$
 $$ \widetilde{T^*X}^{P,Y} = V_P \times^{P_H} G,$$
 hence the relation above translates to 
 \begin{equation}\label{firstdim}
  \dim V_P + \dim (P_H\backslash H) = \dim V + \dim \B_{L(X)}.
 \end{equation}

 On the other hand, by Corollary \ref{polarizedtorsor}, the general fiber of the map $T^*X\to \c_X^*$ is a single $G$-orbit, which is of the form $L_1\backslash G$ over the algebraic closure, where $L_1$ is (non-canonically) isomorphic to  $\ker(L(X)\to A_X)$, with $L(X)$ a Levi of $P(X)$; the general $G$-orbit on $\widetilde{T^*X}$ and $\widetilde{T^*X}^{P,Y}$ is isomorphic over the algebraic closure to  $B_{L_1}\backslash G$, where $B_{L_1}$ is a Borel subgroup of $L_1$. 

 Similarly, the general $L$-orbit on $T^*Y_2$ is isomorphic to $T_1\backslash L$, where $T_1$ is (non-canonically) isomorphic to  $\ker(A\to A_X)$ over the algebraic closure; notice that $P(Y_2)$ is the class of Borel subgroups of $L$. If we consider $T^*Y_2$ as a $P$-space, we have to include $U_P$, which acts trivially, in the point stabilizers. 
 
 Hence, the relative dimension of the map $T_Y^*X^{\u_P^\perp} \to T^*Y_2$, by counting dimensions of stabilizers in $G$, is 
 $$ \dim U_P - \dim U_{L(X)} = \dim U_P - \dim \B_{L(X)}.$$
 The analogous relation for the corresponding vector spaces, by counting dimensions of stabilizers in $H$, is
 \begin{equation}\label{seconddim} \dim V_P - \dim V_2  = \dim U_{P_H} - \dim \B_{L(X)} = \dim (P_H\backslash H) - \dim \B_{L(X)}.\end{equation}
 
 Combining \eqref{firstdim} and \eqref{seconddim}, we obtain \eqref{dimensions}.
 
 Since the kernel of the map $V_P\to V_2$ is isotropic (it maps to $0\in \c_X^*$ by the commutativity of \eqref{moddiagram}), the quadratic space $V$ is maximally isotropic if and only if $V_2$ is. 
\end{proof}

We are now ready to prove Proposition \ref{rationality}.

\begin{proof}[Proof of Proposition \ref{rationality}]
 Suppose that the stabilizer $H$ of a point on $X$ is split. Since $H$ is connected (Proposition \ref{correctrepresentative}), a Borel subgroup of $H$ is split, solvable, connected, and is contained in a Borel subgroup $B$ of $G$ (over $F$). Then, the $H$-orbit represented by $1$ on $G/B$ is projective, hence closed, and the corresponding $B$-orbit on $X=H\backslash G$ is defined over $F$. 
 
 Vice versa, if a closed $B_{\bar F}$-orbit on $X_{\bar F}$ is defined over $F$ --- equivalently, has a point with stabilizer $H$, then $B\cap H$ is a split Borel subgroup of $H$, hence $H$ is split.
 
 Now let $(Y,P)$ be a basic orbit-parabolic pair, as in Lemma \ref{lemmaclosedorbits}. By Lemma \ref{lemmaY2}, the variety $Y_2 = Y/U_P$ is isomorphic to $T\backslash \PGL_2$ for some torus $T$ or to $\SL_2$, hence the closed $B$-orbit(s) on $Y_2$, and on $Y$, are defined over $F$, unless we are in the case of $T\backslash \PGL_2$ with $T$ not belonging to any Borel $F$-subgroup, i.e., with $T$ non-split. This is equivalent to each, equivalently one, fiber of $T^*Y_2\to \c_X^*$ being isotropic, as one can see by direct calculation: in the case of $T\backslash \PGL_2$, this is the map $\mathfrak t^\perp \to \mathfrak t^\perp\sslash T$ which is isomorphic to the norm map for the quadratic splitting ring of $T$, and in the case of $\SL_2$ it is simply the determinant on $\mathfrak{sl}_2$. By Proposition \ref{propVP}, the fibers of $T^*Y_2\to \c_X^*$ being isotropic is equivalent to each, equivalently one, fiber of $T^*X \to \c_X^*$ over a point in $Y$ being maximally isotropic. The surjectivity of the map $\widetilde{T^*X}^{P,Y}\to T^*X$ on $F$-points, by the same proposition, implies that the same is true for every fiber of $T^*X$ hence, applying this argument in the reverse direction, one closed $B_{\bar F}$-orbit being defined over $F$ implies that every closed $B_{\bar F}$-orbit is defined over $F$.
\end{proof}

Recall that we have been assuming that $G$ is split, and that $X$ is a ``correct representative'' of its class modulo $G$-isomorphisms, see Definition \ref{defcorrectrep}; the statement is not true for varieties such as $\mathcal N(T)\backslash \PGL_2$. 

\begin{remark}\label{typeGalwayssplit}
For spherical roots of type $G$, all $B$-orbits are of rank one (maximal), hence defined over $F$ by Lemma \ref{maxrankdefined}. Therefore, by the proposition, in that case all stabilizers are split.
\end{remark}

Finally, the next proposition explains how many isomorphism classes of varieties over $F$ correspond to each line of Table \eqref{thetable}, to which our results apply:

\begin{proposition}\label{uniqueclass?}
 If $X_{\bar F}$ is the correct representative of a class of rank one affine  homogeneous spherical $G_{\bar F}$-varieties (with $W_X=\ZZ/2$) over $\bar F$, then a form of $X$ as a $G$-variety over $F$, where $G$ denotes the split form of $G_{\bar F}$, always exists. Moreover, if the spherical root is of type $T$, there is a unique such form with split stabilizers, and if the spherical root is of type $G$ the isomorphism classes of such forms are naturally a torsor for the group $F^\times/(F^\times)^2$, and stabilizers are always split.
\end{proposition}

\begin{proof} 
The existence of a model over $F$ follows from the general Theorem 0.2 of \cite{Borovoi}. To apply it, we first replace the spherical subgroup $H_{\bar F}$ (over the algebraic closure) by its \emph{spherical closure} $\bar H_{\bar F}\subset \mathcal N(H_{\bar F})$; I refer the reader to the aforementioned reference for the definition of spherical closure. By the theorem, the resulting variety $\bar X_{\bar F} = \bar H_{\bar F}\backslash G_{\bar F}$ has a form $\bar X$ as a $G$-variety over $F$. As in the proof of Proposition \ref{correctrepresentative}, if we replace $G$ by the simply connected cover of its derived group and $\bar H$ by the identity component of its preimage in this simply connected cover, we will obtain a variety $X$ of rank one whose weight lattice is ``as large as possible'', that is, spanned by the spherical root $\gamma$ if that is of type $T$, and by $\frac{\gamma}{2}$ if it is of type $G$; thus, $X$ has to be a form of $X_{\bar F}$ (and, in particular, the action of the simply connected cover factors through $G$).

Given such a form $X$, the set of $G$-forms of $X$ is parametrized by the first Galois cohomology group $H^1(F,\Aut^G(X))$. 

If the spherical root is of type $G$, then, by Remark \ref{typeGalwayssplit}, stabilizers are always split. Moreover, by Proposition \ref{correctrepresentative}, $\Aut^G(X)=\ZZ/2$, so the forms are a torsor for $H^1(F,\Aut^G(X)) = F^\times/(F^\times)^2$. 

If the spherical root is of type $T$, then there is a Borel orbit (over the algebraic closure) of rank zero, and therefore the rank of $H_{\bar F}$ is equal to the rank of $G$, where $H_{\bar F}$ denotes the stabilizer of any point on $X_{\bar F}$. 
Let $T_{\bar F}\subset H_{\bar F}$ be a maximal torus. The Lie algebra $\g_{\bar F}$ splits into a direct sum of $T_{\bar F}$-eigenspaces, and the subalgebra $\h_{\bar F}$ is a subsum of that. If, now, $G$ is defined and split over $F$, we may assume that $T_{\bar F}$ is the extension to $\bar F$ of a maximal split torus $T\subset G$, hence the eigenspaces are defined over $F$, and the subalgebra $\h_{\bar F}$ is the extension to $\bar F$ of a subalgebra over $F$. In other words, there is a form of $X$ such the stabilizer of a point (hence every point, by Proposition \ref{rationality}) is split.

Assume that $X$ is such a form. 
Let $(Y,P_\alpha)$ be a basic orbit-parabolic pair, so that $Y$ contains a $B$-orbit of rank one, and two closed $B$-orbits of rank zero. These closed $B$-orbits now are defined over $F$, by Proposition \ref{rationality}. Arguing as in the proof of Proposition \ref{correctrepresentative}, we get an injection
$$ \Aut^G(X) \hookrightarrow \Aut^{L_\alpha}(Y/U_{P_\alpha}),$$
where $L_\alpha$ is the Levi quotient of $P_\alpha$. The weight lattice of $Y/U_{P_\alpha}$ is spanned by $\alpha$, hence the $L_\alpha$-variety $Y/U_{P_\alpha}$ is isomorphic to the quotient of $L_\alpha$ by a maximal torus $T$, which now has to be split. The group $\Aut^{L_\alpha}(Y/U_{P_\alpha})$ is isomorphic to $\ZZ/2$, with the nontrivial automorphism interchanging the two closed Borel orbits. Hence, if $\sigma$ is a nontrivial $G$-automorphism of $X$, it interchanges the two closed $B$-orbits in $Y$. As a result, for any nontrivial element of $H^1(F,\Aut^G(X))$, defining a form $X'$ of $X$, these $B$-orbits are not defined over $F$ in $X'$, and, by Proposition \ref{rationality} again, stabilizers of points on $X'$ are not split. Therefore, $X$ is the unique form with split stabilizers.

\end{proof}

\begin{example}
 Consider the $G=\SO_4=(\SL_2\times\SL_2)/\{\pm 1\}^\diag$-action on $\GL_2$. All orbits are isomorphic to $H\backslash G=\SO_3\backslash\SO_4=\PGL_2\backslash \SO_4$, but the $G(F)$-conjugacy class of the embedding of $H$ in $G$ depends on the square class of the determinant. 
\end{example}

\section{Structure and resolution of $X\times X/G$} \label{sec:resolution}

In this section we study the diagonal action of $G$ on $X\times X$, and the morphism to the invariant-theoretic quotient:
$$ X\times X\to \C_X:=X\times X\sslash G.$$

In the group case $X=H$, $G=H\times H$, the invariant theoretic quotient $X\times X\sslash G$ coincides with the invariant theoretic quotient of $H$ by $H$-conjugacy, and is naturally identified with $A_H\sslash W_H$, by the Chevalley isomorphism. Moreover, the quotient map $H\to A_H\sslash W_H$ is smooth over the points where the quotient map $A_H\to A_H\sslash W_H$ is smooth.

A Chevalley isomorphism for $X\times X\sslash G = H\backslash G\sslash H$ was proven by Richardson \cite{Richardson} in the case where $X=H\backslash G$ is a symmetric variety. A Chevalley isomorphism for general reductive group actions on affine varieties, in terms of fixed point sets of generic stabilizers, was proven by Luna and Richardson \cite{LR}. Almost all spherical varieties of Table \eqref{thetable} are symmetric, but this is not the case for the examples denoted by $G_2$ and $B_3''$; those two are symmetric under the action of a bigger group, but \emph{a priori} the invariant-theoretic quotient $X\times X\sslash G$ could be different. In any case, I do not know how to deduce all the results that are needed directly from the aforementioned references, even for symmetric spaces (and, in fact, do not all hold without choosing the ``correct representatives'' in each equivalence class, see Definition \ref{defcorrectrep}). 

Thus, in this and the following section, we introduce a more conceptual way to analyze this quotient, for varieties of rank one, and prove the results that we need to analyze orbital integrals. Hopefully, this approach will also be useful for higher-rank cases, where only a small fraction of spherical varieties are symmetric. The main results that we need are the following:
\begin{enumerate}
 \item The invariant-theoretic quotient $\C_X$ is canonically isomorphic to $A_X\sslash W_X$ --- which is just an affine line $\mathbbm A^1$ (Proposition \ref{GITquotients}).
 \item The quotient map $X\times X\to \C_X$ is smooth away from two points $[\pm 1]$ on $\C_X$ (Proposition \ref{genericisomorphism}).
 \item Each point of $\C_X$ corresponds to a closed $G$-orbit in $X\times X$; if $X_{\pm 1}$ denote the closed orbits over the singular points $[\pm 1]$, in an \'etale neighborhood of those the map $X\times X\to \C_X$ is modeled on a non-degenerate quadratic form (see \eqref{Luna} and Proposition \ref{genericisomorphism}).
 \item The codimensions $d_{\pm 1}$ of the singular closed orbits $X_{\pm 1}$ satisfy a formula of the form:
 $$d_1+ d_{-1}=2\epsilon \left<\rho_{P(X)},\check\gamma\right>+2.$$
 (I point the reader to Theorem \ref{thmintegrationformula} for the result and notation.)
 
\end{enumerate}

This general approach will be provided to us by Knop's theory of the cotangent bundle of $X$ and the invariant collective motion --- the infinitesimal action of the $G$-invariant Hamiltonian vector fields obtained from the invariant moment map \eqref{invmoment}. It will turn out, in rank one, that even for the non-symmetric cases the quotient $\C_X$ is (canonically) isomorphic to $A_X\sslash W_X$, with singularities of the quotient map $X\times X\to A_X\sslash W_X$ only over the singularities of $A_X\to A_X\sslash W_X$; this fact does not generalize to higher rank, although the constructions of Knop do.

By the end of this section, we will have constructed a resolution of the $G$-space $X\times X$, denoted
$$ \R : \P J_X \to X\times X.$$
As mentioned in the introduction, the term ``resolution'' refers here to the fibers of the map $X\times X\to \C_X$, which under the resolution become normal crossings divisors.

\subsection{Knop's abelian group scheme; the $[\pm 1]$ and nilpotent divisors}\label{ssKnopsgroupscheme}

Consider $A_X\times \a_X^*$ as a constant group scheme over $\a_X^*$, with the simultaneous action of $W_X$ on $A_X$ and on $\a_X^*$.  
Let  \label{refJ}
\begin{equation}J = (\Res_{\a_X^*/\c_X^*} (A_X\times \a_X^*))^{W_X},\end{equation}
where $\Res_{\a_X^*/\c_X^*}$ denotes Weil's restriction of scalars from $\a_X^*$ to $\c_X^*$.
It is a group scheme over $\c_X^*$, whose sections over any $U\subset \c_X^*$ (as a special case of the universal property of Weil restriction) are the $W_X$-invariant sections of $A_X\times \pi^{-1}(U)$, where $\pi:\a_X^*\to \c_X^*$ is the canonical map. It comes with a canonical birational morphism
\begin{equation}\label{birational} J\to (A_X\times \a_X^*)\sslash W_X,\end{equation}
obtained functorially by identifying the right hand side with 
$$\Res_{\a_X^*/\c_X^*} \left((A_X\times \a_X^*)\sslash W_X\times_{\c_X^*} \a_X^*\right)^{W_X}.$$
This map is an isomorphism over $\mathring{\c}_X^* = \c_X^*\smallsetminus\{0\}$. It is an easy exercise in restriction of scalars to see that, in rank one with $W_X=\ZZ/2$, this scheme is given in coordinates by 
\begin{equation}\label{Jcoordinates} J \simeq \Spec F[t_0,t_1,\xi]/(t_0^2-\xi t_1^2 - 1),\end{equation}
where we have identified $\c_X^*=\mathbbm A^1$, with a coordinate $\xi$ such that $\xi=0$ corresponds to $0\in\a_X^*$. 
This group scheme is also familiar as the group scheme of regular centralizers over the Kostant section of the Lie algebra $\mathfrak{sl}_2$, under the adjoint action of $\SL_2$ \cite[\S 2.4]{Ngo-FL}.

We have a canonical identification
\begin{equation}\label{Jbc} J\times_{\c_X^*} \mathring{\a}_X^* \simeq A_X \times \mathring{\a}_X^*,
\end{equation}
compatible with pullback of sections from $\c_X^*$ to $\a_X^*$, and the identification of sections of $J$ with $W_X$-equivariant sections of $A_X\times \a_X^*$.

On the other hand, the fiber of $J$ over $0\in \c_X^*$ is isomorphic to $\Ga \times \ZZ/2$. We let $J^0\subset J$  \label{refJ0} denote the open group subscheme whose fiber over any point of $\c_X^*$ is the identity component of the fiber of $J$. The group scheme $J$ has a canonical action of $\Gm$, induced from the action of $\Gm$ on $A_X \times \a_X^*$ (on the second factor). In the coordinates above:
\begin{equation}
 a\cdot (\xi, t_0, t_1) = (a^2\xi, t_0, a^{-1} t_1).
\end{equation}

In what follows, for any variety $Y$ equipped with a map to $\c_X^*$ we will denote 
$$J\bullet Y:= J\times_{\c_X^*} Y$$
(and similarly for $J^0$).

We can distinguish three divisors on $J\bullet Y$, two of them to be denoted $[\pm 1]_Y$, and another to be denoted $(J\bullet Y)^\nilp$. The first two are the images of $(\pm 1) \cdot Y$ (where $(\pm 1)$ are understood as $W_X$-invariant sections of $A_X$ over $\a_X^*$, hence as sections of $J$ over $\c_X^*$), and the third is the preimage of $0\in \c_X^*$. In the coordinates used above, the sum of the divisors $[\pm 1]_Y$ is given by the equation $t_1=0$ (with the value of $t_0$ distinguishing the irreducible components), and $(J\bullet Y)^\nilp$ is given by $\xi =0$. Notice that the union of these three divisors is precisely the preimage of the corresponding points $[\pm 1]\in A_X\sslash W_X$ under the map $J\bullet Y\to A_X\sslash W_X$ descending from \eqref{birational}; in the coordinates above, this union corresponds to the equations $t_0=\pm 1$.

Notice that $J\bullet Y$ is smooth over $Y$ (because it is obtained by base change from the smooth group scheme $J\to \c_X^*$), and the divisors $[\pm 1]_Y$ are isomorphic to $Y$. Hence, if $Y$ is smooth, this is also the case for these divisors and the scheme $J\bullet Y$. On the other hand, the nilpotent divisor $(J\bullet Y)^\nilp$ is smooth on its intersection with the smooth locus of the morphism $Y\to \c_X^*$. 

\begin{lemma}\label{lemmatransverse}
Assume that $Y\to \c_X^*$ is smooth. Then, the divisors $[\pm 1]_Y$ and $(J\bullet Y)^\nilp$ intersect transversely, and the morphism $J\bullet Y \to A_X\sslash W_X$ is smooth away from $[\pm 1]_{Y}$.
\end{lemma}

\begin{proof}
All these properties are stable under smooth base change, and since the morphism $Y\to \c_X^*$ is smooth, the problem reduces to the case $Y=\c_X^*$, $J\bullet Y=J$. 

It is then immediate to check from the equation $t_0^2-\xi t_1^2 = 1$ that the intersections of the divisors $\xi=0$ and $t_1=0$ are transverse.
 
For the second statement, we note that we can identify $A_X\sslash W_X\simeq \mathbbm A^1$ so that the map $J\to A_X\sslash W_X$ corresponds to the coordinate $t_0$. The cotangent space of $J$ is generated by $dt_0, dt_1, d\xi$ subject to the equation 
 $$ 2t_0 dt_0 -2\xi t_1 dt_1 -t_1^2 d\xi =0.$$
 Thus, $dt_0 = 0$ only when $t_1=0$.
\end{proof}

We now specialize to the scheme $Y=T^*X$, endowed with the invariant moment map $\mu_\inv: T^*X\to \c_X^*$.
Lemma \ref{nondegenerate} implies that this map is smooth away from the zero section (whose complement we will denote by $T^*X_{\ne 0}$), hence, by \ref{lemmatransverse}, we get:

\begin{corollary}\label{corollarydivisors}
 The divisors $[\pm 1]_{T^*X_{\ne 0}}$ and $(J\bullet T^*X_{\ne 0})^\nilp$ in $J\bullet T^*X_{\ne 0}$ intersect transversely, and the morphism $J\bullet T^*X_{\ne 0} \to A_X\sslash W_X$ is smooth away from $[\pm 1]_{T^*X_{\ne 0}}$.
\end{corollary}

\subsection{Integration of the invariant collective motion}

The relative Lie algebra $\Lie(J)$ of $J$ over $\c_X^*$ can be canonically identified with the cotangent space of $\c_X^*$. Indeed, sections of both are canonically identified with $W_X$-invariant sections of the cotangent bundle of $\a_X^*$. Thus, a section of the cotangent bundle of $\c_X^*$ can be viewed as a section of $\Lie(J)$, and at the same time induces a Hamiltonian vector field on $T^*X$, by pullback and the symplectic structure. The flow along this vector field is known as the \emph{invariant collective motion}. Knop has shown \cite[Theorem 4.1]{KnAut}
that there is an action of the group scheme $J^0$ on $T^*X$ over $\g_X^*$
\begin{equation}\label{actionJ0T}\xymatrix{\mathcal C: & J^0 \bullet T^*X \ar[dr] \ar[rr] && T^*X \ar[dl]\\ && \g_X^*}\end{equation}
that integrates this vector field. Over $\mathring\c_X^*=\c_X^*\smallsetminus\{0\}$, this action lifts through the isomorphism \eqref{Jbc} to the canonical action of $A_X$ that we discussed \eqref{AXaction} on the regular semisimple part of the polarization $\widehat{T^*X}$. In particular, on the regular semisimple part this action is induced from the action of the centralizers of coadjoint vectors: if $v\in T^*X^\rs$ with image $\mu(v)=Z\in \g^*$, the centralizer of $Z$ is a twisted Levi $L$ (conjugate over the algebraic closure to a Levi of $P(X)$), and it acts on $v$ through a quotient which, over the algebraic closure, is isomorphic (up to the $W_X$-action) to $A_X$.

This action may, but does not always, extend to $J$, as the following examples show:

\begin{example}
Let $X=\SL_2$ under the $G=\SO_4=\SL_2\times \SL_2/\{\pm 1\}^\diag$-action. Then $\c_X^*=\mathfrak{sl}_2\sslash \SL_2$ (under the adjoint action), and $J$ can be identified with the group scheme of regular centralizers over $\c_X^*$, i.e., the group scheme of centralizers in $\SL_2$ over a Kostant section $$\c_X^*\hookrightarrow \mathfrak{sl}_2^{\reg}.$$
Hence, it acts (faithfully) on $(T^*\SL_2)_{\ne 0}$, by either left or right translation; it is easy to see that this extends to the trivial action on the zero section. In particular, the action of $J^0$ extends to $J$.
\end{example}

\begin{example}
Let $X=\Gm\backslash \PGL_2$. The group scheme of regular centralizers for $\PGL_2$ again acts faithfully on $T^*X$, but in this case it is isomorphic to $J^0$. However, one can easily see that the entire group scheme $J$ acts, with the action of the $(-1)$-section induced from the nontrivial $G$-automorphism of $X$.
\end{example}

\begin{example} \label{pgl3}
Let $X=H\backslash G=\GL_2\backslash \PGL_3$, the variety of direct sum decompositions $\Ga^3=V_2\oplus V_1$ of a based three-dimensional vector space into the sum of a two- and a one-dimensional subspace. The fiber $\mathfrak h^\perp$ of $T^*X$ over the point $x_0$ corresponding to the decomposition $\left<e_1, e_2\right> \oplus \left< e_3\right>$ is isomorphic to $\Std\oplus \Std^*$, the direct sum of the standard representation and its dual, as a (right) representation of $H=\GL_2$. The quadratic form of the invariant moment map is 
$$(v, v^*)\mapsto \left<v, v^*\right>.$$
Hence, $T_{x_0}^*X^\nilp=$ the variety of mutually orthogonal pairs $(v, v^*)$. 

Represent $H$ as the upper left copy of $\GL_2$ in $\PGL_3$, identify $\g^*= \g = \mathfrak{sl}_3$ through the trace pairing, and consider representatives 
$$Z_{\epsilon, y}= \begin{pmatrix}
  && 0 \\ && 1 \\ y & \epsilon^2 & 0
 \end{pmatrix},
 $$
with $\epsilon\ne 0$, for the split regular semisimple $H$-orbits on $\mathfrak h^\perp$ (i.e., those over the image of $\mathring\a_X^*(F)\to \c_X^*(F)$). The variable $y$ is redundant, at the moment, but will play a role as $\epsilon \to 0$. The centralizer of $Z_{\epsilon, y}$ under the right coadjoint representation of $G$ is the torus of matrices (modulo center) of the form 
$$ \begin{pmatrix}
    a\\
    \frac{(-2a+b+c)y}{2\epsilon^2} & \frac{b+c}{2} & \frac{b-c}{2 \epsilon} \\
    \frac{(b-c)y}{2\epsilon} & \frac{(b-c)\epsilon}{2} & \frac{b+c}{2}
   \end{pmatrix}.$$
It acts on the point $x_0\in X$ through the quotient $(a,b,c)\mapsto \frac{b}{c}$, which is isomorphic to $A_X\simeq \Gm$ (up to inversion). 

Let us examine whether the action of $-1\in A_X$ on $T^*X^\rs$ extends to the nilpotent limit $\epsilon\to 0$. Representing $-1$ by a matrix as above, corresponding to $(a,b,c)=(1,1,-1)$: 
$$g_{\epsilon, y}= \begin{pmatrix}
    1\\
    \frac{-y}{\epsilon^2} & 0 & \frac{1}{ \epsilon} \\
    \frac{y}{\epsilon} & \epsilon & 0
   \end{pmatrix} \sim \begin{pmatrix}
    1\\
    -y & 0 & \epsilon \\
    \frac{y}{\epsilon^2} & 1 & 0
   \end{pmatrix}$$
(where $\sim$ means same left $H$-coset), we easily see that $x_0\cdot g_{\epsilon, y}$ does not have a limit in $X$ as $\epsilon \to 0$, unless $y=0$. On the other hand, if $y$ also tends to zero in such a way that $\frac{y}{\epsilon^2}$ has a limit, the action extends to the limit. Geometrically, this means that if we blow up $\mathfrak h^\perp= \Std\oplus \Std^*$ over the divisor $\Std \cup \Std^*$, and remove the strict transform of the nilpotent divisor, then the action of $-1$ extends to the blowup --- this will be relevant when discussing the ``second'' singular orbit on $X\times X/G$, below, see Example \ref{pgl3-second}.
\end{example}

\begin{lemma}\label{torsor}
 For $X$ of rank one with $W_X=\ZZ/2$, the restriction of $T^*X\to \g_X^*$ over $\g_X^{*,\rs}$ (that is, over $\mathring\c_X^*$ under the invariant moment map) is a $J\bullet \g_X^{*,\rs}$-torsor.
 \end{lemma}

\begin{proof}
 This is a consequence of Corollary \ref{polarizedtorsor}: Since $\widehat{T^*X}^\rs$ is an $A_X$-torsor over $\hat\g_X^{*,\rs}$, and the action is the base change of the action of the group scheme $J\bullet \g_X^{*,\rs}$ on $T^*X^\rs$, the latter is also a torsor.
\end{proof}

\subsection{Resolution in a neighborhood of the diagonal}

We denote by $N^*_A B$ the conormal bundle in a smooth variety $B$ of a subvariety $A$. In the case that $B=X\times X$, we will be using the shorthand notation $N^*_A$ for $N^*_A B$; we denote by $N_A$ the normal bundle.

The fiber product $J\bullet T^*X$ carries a natural $\Gm$-action, induced from the action on $T^*X$ and on $J$. The quotient $\Gm\backslash (J\bullet T^*X_{\ne 0})$ will be denoted by the symbol of projectivization, $\P (J\bullet T^*X)$. We use analogous notation for all similar spaces with a $\Gm$-action.

Consider the combination of the projection and action maps:
\begin{equation}\label{Jtoconormal}
 J^0\bullet T^*X \to T^*X \times_{\g^*} T^*X.
\end{equation}

An immediate corollary of Lemma \ref{torsor} is:
\begin{corollary}\label{cortorsor}
The map \eqref{Jtoconormal} is  an isomorphism over $\mathring\c_X^*$ (i.e., on the sets of regular semisimple vectors).
\end{corollary}

\begin{proof}
By Lemma \ref{nondegenerate}, over $\mathring\c_X^*$ we have $T^*X^\rs \times_{\g^*} T^*X^\rs =  T^*X^\rs \times_{\g_X^{*,\rs}} T^*X^\rs$, and by Lemma \ref{torsor} the right hand side is isomorphic to $J^0\bullet T^*X^\rs = J\bullet T^*X^\rs$.
 
\end{proof}

The space on the right hand side of \eqref{Jtoconormal} can be thought of as the union of all conormal bundles to all $G$-orbits on $X\times X$. Indeed, the conormal bundle on a $G$-orbit for the diagonal action on $X\times X$ is the subbundle of the cotangent bundle $T^*X \times T^*X$ determined by the vanishing of the diagonal moment map $(v_1, v_2)\mapsto \mu(v_1)+\mu(v_2)$; multiplying $v_2$ by $-1$, this becomes the fiber product over $\g^*$.
This union of conormal bundles comes with its own map to $\g^*$, which is not the moment map for the $G^\diag$-action on $X\times X$ (which is zero), but ``remembers'' the fact that $X\times X$ had a $G\times G$-action. This is the microlocal analog of the spectral decomposition of the relative trace formula.

The rough idea behind the resolution of $X\times X$ that we are about to construct is that, generically, the projectivization of this union of conormal bundles is isomorphic to $X\times X$ (because we are in rank one, and generic $G$-orbits will be of codimension one), hence the projectivization of $T^*X \times_{\g^*} T^*X$ is, roughly, a ``resolution'' of $X\times X$. However, this space is quite singular, and we will use an extension of the space on the left as a smooth replacement.

\begin{proposition}\label{arounddiagonal}
Consider the map 
\begin{equation} \R: \P(J^0\bullet T^*X) \to X\times X,
\end{equation}
 descending from \eqref{Jtoconormal}. We regard $\P T^*X$ as a divisor in $\P (J^0\bullet T^*X)$,  descending from the divisor $[1]_{T^*X}$ (see \S \ref{ssKnopsgroupscheme}).

The map $\R$ factors through the blowup $\widetilde{X\times X}^1\to X\times X$ at the diagonal $X_1:=X^\diag$, and is an isomorphism from a $G$-stable neighborhood of $\P T^*X$ to a neighborhood of the exceptional divisor in $\widetilde{X\times X}^1$.
\end{proposition}

\begin{proof}

By the universal property of blowups, the map $\R$ factors through a morphism to the blowup
\begin{equation}\label{toblowup}
 \tilde\R : \P J_X\to  \widetilde{X\times X}^1,
\end{equation}
sending the divisor $\P T^*X$ to the exceptional divisor of the blowup, which is isomorphic to $\P TX$. 

To show that this map is an isomorphism in a $G$-stable neighborhood of the divisor $\P T^*X$, it is enough to show that the induced map $d\tilde\R$ from the normal bundle of the divisor $\P T^*X$ to the normal bundle of the exceptional divisor is an isomorphism. 

The normal bundle of $\P T^*X$ can be identified with $\P(\Lie J \bullet T^*X)$,  and the normal bundle of the exceptional divisor can be identified with the blowup of the tangent bundle $TX$ at the zero section, i.e., of the normal bundle to $X_1 = X^\diag$. The map $d\tilde\R$ is lifted from the analogous map
\begin{equation}\label{mapLie} d\R: \P(\Lie J \bullet T^*X_{\ne 0}) \to TX= N_{X_1},\end{equation}
the partial differential of the map $\R$. We compute this map:

Recall that $\Lie J$ is canonically isomorphic to the cotangent space of $\c_X^*$. 
A section $\sigma$ of $\Lie J$ gives, by pullback of differential forms via the invariant moment map, a section $\mu_\inv^*\sigma$ of the cotangent bundle of $T^*X$, hence a vector field $v_\sigma$ on $T^*X$, by the symplectic structure. If $\pi: T^*X\to X$ denotes the canonical projection, and $\pi^*(TX)$ is the pullback of the tangent bundle, we have a canonical projection of vector bundles on $T^*X$ 
$$\pr: T(T^*X)\to \pi^*(TX),$$
corresponding to the projection of vector fields on $T^*X$ ``to the $X$-direction''.

Let $T^*_X (T^*X) = T^*(T^*X)/\pi^*(T^*X)$ denote the relative cotangent bundle of $T^*X$ over $X$; we similarly have a projection of vector bundles
$$\pr': T^*(T^*X) \to T^*_X(T^*X).$$
Moreover, we have canonical identifications
$$ T^*_X (T^*X) \xrightarrow\sim T^*X\times_X TX \overset{\sim}\leftarrow \pi^*(TX),$$
and the image $\pr(v_\sigma)$ of the aforementioned vector field, as a section of $\pi^*(TX)$, coincides under this identification with ``the restriction of $\mu_\inv^*\sigma$ to the fiber direction'', that is, with the section $\pr'(\mu_\inv^*\sigma)$ of $T^*_X(T^*X)$.

Consider, for example, an identification $\c_X^*\xrightarrow{\xi}\mathbbm A^1$ (with $0\in\c_X^*$ mapping to $0\in \mathbbm A^1$), and let $\sigma=d\xi$. Then, $\xi\circ\mu_\inv$ can be viewed as a quadratic form on $T^*X$, and it gives rise to a map 
$$\iota_\xi:T^*X\to TX$$ 
over $X$. As we have seen in Lemma \ref{nondegenerate}, 
the quadratic form is nondegenerate; hence, $\iota_\xi$ is an isomorphism. 

The differential of the quadratic form, restricted to each fiber, is the graph of $\iota_\xi$, considered as a subset of $T^*X\times_X TX = \pi^*(TX)$. 
Therefore, the section of $\Lie J$ corresponding to $d\xi$ defines a vector field on $T^*X$, whose projection to the $X$-direction is the graph of $\iota_\xi$. Hence, the map $d\R$ descends from 
$$ v^*\in T^*X_{\ne 0}  \mapsto (d\xi, v^*) \in \Lie J \bullet T^*X_{\ne 0} \mapsto \iota_\xi(v^*) \in TX.$$

In particular, the map $d\R$ is fiberwise an isomorphism, hence its lift $d\R$ to the blowup  is an isomorphism, and $\tilde R$ is an isomorphism from a $G$-stable neighborhood of the divisor $\P J_X$ to a $G$-stable neighborhood of the exceptional divisor.

\end{proof}

\begin{remark}
 Notice that, under the isomorphism $\Lie J = T^*\c_X^*$, the action of $\Gm$ on $\Lie J$ is the following: it acts in the canonical way on the base ($a\in \Gm$ acts by multiplying a coordinate $\xi$ on $\c_X^*$ by $a^2$); this induces an inverse pullback isomorphism:
$$ (a^*)^{-1}: T^*\c_X^* \to T^*\c_X^*,$$
and we multiply this by $a$, fiberwise. In terms of the coordinates $t_0^2-\xi t_1^2=1$, the action of $\Gm$ on $J$ is given by $a\cdot (\xi, t_0, t_1) = (a^2\xi, t_0, a^{-1} t_1)$. One can now directly see that the map $\Lie J \bullet T^*X_{\ne 0} \to TX$ that was described above  is, indeed, $\Gm$-equivariant.
\end{remark}

The proposition implies:

\begin{corollary}\label{corregulardense}
 There is an open, dense, $G$-stable subset $(X\times X)^\circ \subset X\times X$ on which every $G$-orbit is of codimension one, and the nonzero vectors of its conormal bundle are regular semisimple. Moreover, every $G$-orbit on this subset is isomorphic, over the algebraic closure, to $L_1\backslash G$, where $L_1$ is (non-canonically) isomorphic to  $\ker(L(X)\to A_X)$, with $L(X)$ a Levi of $P(X)$. 
\end{corollary}

\begin{proof}
 Indeed, consider the space $T^*X\times_{\g^*} T^*X$ as the union of the conormal bundles of all $G$-orbits on $X\times X$; the map
 from $J\bullet T^*X^\rs$ is an isomorphism onto the regular semisimple subset 
 $T^*X^\rs\times_{\g^{*,\rs}} T^*X^\rs$ by Corollary \ref{cortorsor}, and the space $J\bullet T^*X^\rs$ is a union of $\Gm\times G$-orbits of codimension one. By Proposition \ref{arounddiagonal}, a dense open $G$-stable subset of its projectivization is isomorphic to an open subset of a $G$-stable neighborhood of the diagonal $X_1=X^\diag$. The description of $G$-orbits follows from Corollary \ref{polarizedtorsor}.
\end{proof}

\subsection{Closed orbits and invariant-theoretic quotients}

The following lemma will be very basic in our analysis of the space $X\times X$:

\begin{lemma}\label{lemmaregular}
 Every $G$-orbit in the image of the regular semisimple set $\P(J^0\bullet T^*X^\rs)$ under the map $\R$ is closed. 
 
 In particular, non-closed $G$-orbits on $X\times X$ do not contain regular semisimple vectors in their conormal bundles.
\end{lemma}

\begin{proof}
 The second statement follows from the first, by Corollary \ref{cortorsor}.

To prove the first, we need some preliminary results. Choose a basic orbit-parabolic pair $(Y,P)$ as in Lemma \ref{lemmaclosedorbits}. The subvariety $Y\times Y\subset X\times X$ is closed, hence the map $(Y\times Y)\times^P G\to X\times X$ is proper.

Consider the space $\widetilde{T^*X}^{P,Y} = T_Y^*X^{\u_P^\perp} \times^P G$, defined in \S \ref{ssrationality}. By Lemma \ref{partialcover}, it surjects onto $T^*X$, hence $J^0\bullet \widetilde{T^*X}^{P,Y}$ surjects onto $J^0\bullet T^*X$ (and same for the subsets of regular semisimple vectors). But the latter has dense image in $X\times X$, by Proposition \ref{arounddiagonal}, hence so does the former. The subspace $J^0\bullet T_Y^*X^{\u_P^\perp}$ maps to $Y\times Y$, and we have a commutative diagram:
 \begin{equation}\label{Ydiagram} \xymatrix{
 J^0\bullet \widetilde{T^*X}^{P,Y} \ar[r] \ar[d]^{\R_Y} & J^0\bullet T^*X \ar[d] \\
 \tilde Y\times^P G \ar[r] & X\times X,}\end{equation}
where $\tilde Y\subset Y\times Y$ denotes the closure of the image of $J^0\bullet T_Y^*X^{\u_P^\perp}$ in $Y\times Y$ under the map to $Y\times Y$. Hence, the map 
$$ \tilde Y\times^P G\to X\times X$$ 
is surjective and proper. It is thus enough to show that every $P$-orbit in the image of $J\bullet T_Y^*X^{\u_P^\perp,\rs}$ in $\tilde Y$ is closed.

Let $H\subset G$ be the stabilizer of a point on $Y$, and $H_P=H\cap P$. 
Consider the quotient $Y\to Y_2 = Y/U_P$, and remember that $Y_2$ is isomorphic to $T\backslash \PGL_2$ or $\SL_2$ (Lemma \ref{lemmaY2}), with the connected center of the Levi quotient of $P$ acting trivially; in particular, stabilizers of points in $Y_2$ are reductive.
Thus, we can fix a Levi decomposition $H_P = H_L \cdot H_U$, and a Levi subgroup $L\subset P$ which contains $H_L$. Notice that $Y_2 = H_L\backslash L$, hence we can also consider $Y_2$ as a \emph{subvariety} of $Y$ (depending on the choices that we have made). 

Use an invariant bilinear form on $\g$ to identify $\u_P^\perp = \mathfrak p$. By semisimplicity, the image of every element of $T_Y^*X^{\u_P^\perp,\rs}$ under the moment map is $U_P$-conjugate to an element of $\mathfrak l$. Thus, it is enough to show that the $P$-orbit of $\R_Y(a\cdot v)$ is closed, for any $v\in T_Y^*X^{\u_P^\perp,\rs}$ with $\mu(v)\in \mathfrak l$, and any $a\in J$ over the image of $v$ in $\c_X^*$. 

Recall that the $J$-action on $v$ is induced from the action of the centralizer of $\mu(v)$ in $G$; since this acts by a one-dimensional quotient, the same action is induced from the centralizer of  $\mu(v)$ in $L$, which is a torus acting nontrivially. Hence, given $a\in J$ over $\mu_\inv(v)$, there is an $l\in L$ with $a\cdot v = v\cdot l$. Hence, considering $Y_2$ as a subvariety of $Y$, $\R_Y(a\cdot v)\in Y_2\times Y_2$. To avoid confusion, we will be denoting by $y$ a point in $Y_2\times Y_2$ considered as a subset of $Y\times Y$, and by $\bar y$ its image in $Y_2\times Y_2$ considered as a quotient of $Y\times Y$.
It is immediate to confirm that the $L^\diag$-orbit of $\R_Y(a\cdot v)$ is closed in $Y_2\times Y_2$; for every point $y$ on that orbit, the preimage of $\bar y$ under the \emph{quotient} $\tilde Y\to Y_2\times Y_2$ is the closure of the $U_P$-orbit of $y$, considered as a point in $\tilde Y$. But $\tilde Y$ is affine, hence the $U_P$-orbit of $y$ is closed. Hence, the $P$-orbit of $\R_Y(a\cdot v)$ is closed, completing the proof of the lemma.
\end{proof}

Now recall the birational map \eqref{birational}: $J\to (A_X\times \a_X^*)\sslash W_X$.   
This induces a map $$J^0\bullet T^*X\to (A_X\times \a_X^*)\sslash W_X,$$ which is a smooth geometric quotient\footnote{A smooth geometric quotient $X\to Y$ is a smooth surjective morphism of $G$-varieties, with $G$ acting trivially on $Y$, such that geometric fibers are $G$-orbits. \label{geomquotient}} by the $G$-action when restricted to $J\bullet T^*X^\rs$ (i.e., over $\mathring\c_X^*$), by Lemma \ref{torsor}.

\begin{proposition}\label{GITquotients}
 There is a commutative diagram
 $$ \xymatrix{
 J^0\bullet T^*X_{\ne 0} \ar[r] \ar[d] & T^*X_{\ne 0} \times_{\g_X^*} T^*X_{\ne 0} \ar[r] & X\times X \ar[d] \\
 (A_X\times \a_X^*)\sslash W_X \ar[r] & A_X\sslash W_X\ar[r] & \C_X:=X\times X\sslash G }$$
 which identifies: 
 \begin{itemize}
  \item $\C_X$ with $A_X\sslash W_X$;
  \item for every $c\in \C_X$ with corresponding closed $G$-orbit $C\subset X\times X$, the fiber of $J\bullet T^*X^{\rs}$ over $c$ with the set $N_C^{*,\rs}$ of regular semisimple vectors in the conormal bundle to $C$;
  \item the quotient $N^*_C\sslash G$ with the fiber of $(A_X\times \a_X)^*\sslash W_X$ over $c$.
 \end{itemize}
\end{proposition}

\begin{proof}
Since $T^*X\sslash G=\c_X^*$, we have $\P (J^0 \bullet T^*X)\sslash G = A_X\sslash W_X$, so the composition $\P (J^0 \bullet T^*X)\to X\times X\to \C_X$ indeed factors through a map $A_X\sslash G\to \C_X$. On the other hand, with notation as in the proof of Lemma \ref{lemmaregular}, and by the diagram \eqref{Ydiagram}, the composition 
$$ \P(J^0\bullet \widetilde{T^*X}^{P,Y}) \twoheadrightarrow \P (J^0 \bullet T^*X) \to X\times X$$
factors through a proper, surjective map $\tilde Y\times^P G\to X\times X$. I claim that the composition 
$$ \P(J^0\bullet \widetilde{T^*X}^{P,Y}) \to \P (J^0 \bullet T^*X) \to A_X\sslash W_X$$ 
factors through a $G$-invariant map 
$$\tilde Y \times^P G \to A_X\sslash W_X.$$ Indeed, recall that $\tilde Y$ is the closure of the image of $J^0\bullet T_Y^*X^{\u_P^\perp}$ in $Y\times Y$; the map $J^0\bullet T_Y^*X^{\u_P^\perp}\to \tilde Y\sslash P$ factors through $A_X\sslash W_X$ for the same reasons, but on the other hand we have a quotient map $\tilde Y \to Y_2\times Y_2$ (again, in the notation of the proof of Lemma \ref{lemmaregular}), and $Y_2\times Y_2\sslash P = A_X\sslash W_X$, so the map to $A_X\sslash W_X \to \tilde Y\sslash P$ is an isomorphism. 

Therefore, the map $A_X\sslash W_X\to \C_X$ is surjective.
Proposition \ref{arounddiagonal} implies that it is also birational. Since both $A_X\sslash W_X$ and $\C_X$ are normal, the map has to be an isomorphism.

Recall that geometric points of $\C_X$ correspond bijectively to closed geometric orbits of $G$ on $X\times X$. By Lemma \ref{lemmaregular}, these have to be \emph{precisely} the images of $\P(J^0\bullet T^*X^\rs)$ in $X\times X$. Notice that, by Lemma \ref{torsor}, the set $\P (J^0\bullet T^*X^\rs)$ contains a unique $G$-orbit over any point of $\C_X$.

Let $c\in \C_X$ be such a point, $C\subset X\times X$ the closed $G$-orbit over $c$, and denote by an index $~_c$ various fibers over $c$. By Corollary \ref{cortorsor}, the fiber $J^0\bullet T^*X^\rs$ over $c$ coincides with the regular semisimple part $N^{*,\rs}_C$ of its cotangent bundle.  Hence, 
$$N^{*,\rs}_C\sslash G= (J\bullet T^*X^\rs)_c\sslash G = ((A_X\times\mathring\a_X^*)/W_X)_c.$$

On the other hand, the invariant moment map gives rise to a $G$-invariant map $N^*_C\to \c_X^*$ which, by considering the regular semisimple and the zero vectors, has to be surjective. Since $N^*_C\sslash G$ is normal, it has to coincide with the spectrum of the integral closure of $F[\c_X^*]$ in the function field of $((A_X\times\mathring\a_X^*)/W_X)_c$, which coincides with $((A_X\times \a_X^*)/W_X)_c$; that is, with $\a_X^*$ (up to $\pm 1$), if $c\ne [\pm 1]$, and with $\c_X^*$, if $c=[\pm 1]$.
\end{proof}

\begin{corollary}\label{containsnormal}
 The closure of the image of the map \eqref{Jtoconormal} contains ${N^*_C}$, for any \emph{closed} $G$-orbit $C\subset X\times X$. Every closed $G$-orbit contains regular semisimple vectors in its conormal bundle.
\end{corollary}

\begin{proof}
The second statement was already explained in the proof of the previous lemma, and the first follows by taking the closure of regular semisimple vectors,  and using the fact that $N^*_C$ is irreducible.
\end{proof}

The statement is not true for non-closed $G$-orbits, which could contribute smaller irreducible components to $T^*X\times_{\g^*} T^*X$. 
\begin{example}
 Consider the case of $X=\GL_2\backslash\PGL_3$, discussed in Example \ref{pgl3}, where we identified $T^*_{x_0}X=\h^\perp$ as the representation $\Std\oplus \Std^*$ of $H=\GL_2$. Let $v^*\in h^\perp$ be a nonzero \emph{irregular} nilpotent vector, i.e., either in $\Std$ or in $\Std^*$. The orbit of $x_0$ under its centralizer $G_{v^*}$ is two-dimensional. Hence, the fiber of 
 $$T^*X\times_{\g^*} T^*X$$
 over $v^*\in T_{x_0}^*X$ under the first projection is (at least) two-dimensional, while the fiber of 
 $$J^0\bullet T^*X$$
 over it is one-dimensional. Thus, there are conormal vectors to $G$-orbits on $X\times X$ which are not contained in the closure of the image of \eqref{Jtoconormal}.
\end{example}

\subsection{Blow-up of $X\times X$ at the closed orbits}

We have already seen (Proposition \ref{arounddiagonal}) that the map $\P (J^0 \bullet T^*X)\to X\times X$ is an isomorphism, generically; in particular, generic fibers over $\C_X=A_X\sslash W_X$ (Proposition \ref{GITquotients}) are single $G$-orbits. We can now determine which ones:

\begin{proposition}\label{genericisomorphism}
 \begin{enumerate}
  \item Let $c\in \C_X$ with corresponding closed $G$-orbit $C\subset X\times X$. If $c\ne [\pm 1]$, the linear map $N^*_C \to N^*_C\sslash G$ (which is identified with $\a_X^*$, up to $\pm 1$, by Proposition \ref{GITquotients})
 is an \emph{isomorphism} on each fiber over $C$. If $c=[\pm 1]$, the quadratic map $N^*_C\to N^*_C\sslash G=\c_X^*$ is \emph{nondegenerate} on each fiber over $C$.
  \item The map $\P (J^0 \bullet T^*X)\to X\times X$ is an isomorphism over $\mathring\C_X := \C_X\smallsetminus\{[\pm 1]\}$. In particular, in the setting of Corollary \ref{corregulardense}, we can take $(X\times X)^\circ = (X\times X)\times_{\C_X} \mathring\C_X$. 
  \item The map $X\times X\to \C_X$ is a smooth geometric quotient by the $G$-action away from $[\pm 1]\in \C_X$.
 \end{enumerate}
\end{proposition}

\begin{proof}
Let $V$ be the fiber of $N^*_C$ over a point of $C$, and $H_0$ the stabilizer of that point. If the linear map (when $c\ne [\pm 1]$)
$$ V\to \a_X^*,$$
where $\a_X^*$ is identified, up to $\pm 1$, with the fiber of $A_X\times \a_X^*\sslash W_X$ over $c$,
resp.\ if the quadratic map (when $c=[\pm 1]$)
$$ V\to \c_X^*$$
were trivial on a nonzero, necessarily $H_0$-stable subspace $V_0\subset V$, that space would have an $H_0$-stable complement 
$$ V= V_0 \oplus V_0',$$
identifying the invariant-theoretic quotients
$$N^*_C\sslash G = V\sslash H_0 = V_0'\sslash H_0 = \c_X^*.$$ 
In particular, the $G$-orbit of a generic point of $N^*_C$ (corresponding to an $H_0$-orbit on $V$ not belonging to $V_0'$) is not closed, a contradiction, since by Proposition \ref{GITquotients} and Corollary \ref{corollaryGtransitive} the fibers over all points of $\mathring \c_X^*$ are $G$-homogeneous.

This proves the first claim, and it implies that closed $G$-orbits in $X\times X$ over $c\ne [\pm 1]$ are of codimension one, hence coincide with the whole fiber. Hence, by Corollary \ref{cortorsor}, the map from $J\bullet T^*X^\rs$ to the union $T^*X_{\ne 0} \times_{\g^*} T^*X_{\ne 0}$ of nonzero vectors in the union of conormal bundles is an isomorphism over $\mathring \C_X$, therefore the projectivization $\P(J\bullet T^*X^\rs)$ (or, equivalently, $\P(J^0\bullet T^*X)$) is isomorphic to $X\times X$ over this subset. 

The last claim follows from the analogous claim for $\P (J \bullet T^*X^\rs)$: Since $T^*X^\rs$ is a $J$-torsor over $\g_X^{*,\rs}$ by Lemma \ref{torsor}, we have isomorphisms of geometric quotients
$$ ((X\times X)\times_{\C_X} \mathring\C_X)/G = \left(\P (J \bullet T^*X^\rs)\times_{\C_X}\mathring \C_X\right)/G = \P J^\rs  \times_{\C_X}\mathring \C_X = \mathring \C_X.$$
\end{proof}

Finally, we are ready to construct a resolution of $X\times X$. We now know that this space contains two closed $G$-orbits, $X_1 = X^\diag$ and $X_{-1}$, around which the map to $A_X\sslash W_X$ may fail to be a geometric quotient, namely, the ones over the points $[\pm 1]\in\C_X$. The resolution that we will construct will eventually turn out to be, simply, the blowup at those two subsets. However, generalizing Proposition \ref{arounddiagonal}, we will construct this resolution by a slight modification of the space $J\bullet T^*X$, that we will denote by $J_X$.

The scheme $J_X$ \label{refJX} will be glued from two open subsets: the first is $U_1:=J^0\bullet T^*X$. The second is $U_{-1} := J^0\bullet N^*_{X_{-1}}$, where $X_{-1}$ is the closed $G$-orbit corresponding to $[-1]\in \C_X$. We define $J_X= U_1\cup U_{-1}$, glued over their subsets of regular semisimple vectors as follows: Notice that $U_1^\rs=J\bullet T^*X^\rs$. Moreover, by Proposition \ref{GITquotients}, we have an identification of $N^{*,\rs}_{X_{-1}}$ with the subset $(-1)\cdot T^*X^{\rs}$ of $J\bullet T^*X^\rs$, hence  $U_{-1}^\rs = J\bullet (-1)\cdot T^*X^{\rs}$. This defines the isomorphism
$$ U_{-1}^\rs = J\bullet (-1)\cdot T^*X^{\rs} \ni (j,(-1)\cdot v^*)\mapsto ((-1)\cdot j, v^*)\in J\bullet T^*X^\rs = U_1^\rs,$$
hence the scheme $J_X$. This scheme retains a map 
\begin{equation}\label{birationalX} J_X \to (A_X\times \a_X^*)\sslash W_X,\end{equation}
whose restriction to $\mathring\c_X^*$ is equal to $J\bullet T^*X^\rs$.

The map \eqref{Jtoconormal} extends to $J_X$:
\begin{equation}\label{JtoconormalX} J_X \to T^*X\times_{\g^*} T^*X,
\end{equation}
and we define distinguished divisors $[1]_{J_X}$, $[-1]_{J_X}$ and $J_X^\nilp$, similarly as in \S \ref{ssKnopsgroupscheme}.

Extending Corollary \ref{corollarydivisors}, 
\begin{corollary}\label{corollarydivisorsJX} 
 The divisors $[\pm 1]_{J_{X,\ne 0}}$ and $J_{X,\ne 0}^\nilp$ in $J_{X,\ne 0}$ intersect transversely, and the morphism $J_{X,\ne 0} \to A_X\sslash W_X$ is smooth away from $[\pm 1]_{J_{X,\ne 0}}$.
\end{corollary}

Here, $J_{X,\ne 0}$ denotes the complement of the zero section $J^0\bullet X\subset J^0\bullet T^*X$ in $U_1$, and of the zero section $J^0\bullet X_{-1}\subset J^0\bullet N^*_{X_{-1}}$ in $U_{-1}$.

\begin{proof}
 On the open subset $U_1$, this is contained in Corollary \ref{corollarydivisors}. For $U_{-1}$, the same proof, based on Lemma \ref{lemmatransverse}, works, because of the non-degeneracy statement of the first part of Proposition \ref{genericisomorphism}.
\end{proof}

Now consider the composition $J_X \to T^*X\times_{\g^*} T^*X \to X\times X$. On $J_{X,\ne 0}$, it clearly factors through the projectivization $\P J_X$.

 \begin{proposition} \label{isblowup}
 The morphism 
 $$\R: \P J_X \to X\times X$$ 
 is isomorphic to the blowup of $X\times X$ at the closed $G$-orbits $X_1$ and $X_{-1}$. The preimage of any point of $\C_X=X\times X\sslash G$ under the composition of the maps $\P J_X\to X\times X\to \C_X$ is a normal crossings divisor.
 \end{proposition}

 \begin{proof}
 The statement identifying $\P J_X$ as a blowup has already been proven away from $[-1]\in \C_X$, by a combination of Propositions \ref{arounddiagonal} and \ref{genericisomorphism}. On a $G$-stable neighborhood of $[-1]_{J_X}$, it can be proven by exactly the same arguments as in Proposition \ref{arounddiagonal}. Notice that, now, the map $N^*_{X_{-1}}\to \c_X^*$, viewed as a quadratic form on the fibers by fixing a coordinate $\xi$ on $\c_X^*$, gives rise to a map from the conormal to the normal bundle: 
 $$N^*_{X_{-1}}\to N_{X_{-1}},$$
 and the non-degeneracy of this quadratic form (Proposition \ref{genericisomorphism}) implies that this map is an isomorphism. This fact implies, as in Proposition \ref{arounddiagonal}, that the map from $\P J_X$ to the blowup of $X\times X$ at $X_{-1}$ is an isomorphism around $[-1]_{J_X}$.
 
 The preimage of any point on $\C_X$ is either a unique (smooth) $G$-orbit of codimension one in $X\times X$, or is contained in the divisors $[\pm 1]_{\P J_X}$ and $\P J_X^\nilp$; by Corollary \ref{corollarydivisorsJX}, these have normal crossings.
 \end{proof}

I finish this section by stating rationality properties of the map $N^*_{X_{-1}}\to \c_X^*$, analogous to those of the invariant moment map that were proven in \S \ref{ssrationality}. Notice that, up to this point in this section, we have not used the fact that we are working with the ``correct representative'' of a variety in its class modulo $G$-automorphisms (see Proposition \ref{correctrepresentative}), but now we will.

Let $(Y,P)$ be a basic orbit-parabolic pair, as in Proposition \ref{propVP}; hence, $Y$ is a closed $P$-orbit, and the quotient $Y/\mathcal R(P)$ is isomorphic to $T\backslash \PGL_2$ or to $\PGL_2$. We have defined a cover  $\widetilde{T^*X}^P\to T^*X$, and an irreducible component $\widetilde{T^*X}^P$ thereof; by base change, we get analogous covers for $J\bullet T^*X$. We let
$$ \widetilde{J_X}^P = \{(v,P')| v\in J_X,  P'\sim P, \mu(v) \in \u_{P'}^\perp\}$$
(where $\mu$ also denotes the moment map for $J_X$), and we let $ \widetilde{J_X}^{P,Y}$ be the closure of $J\bullet \widetilde{T^*X}^{P,Y,\rs}$ in $\widetilde{J_X}^P$. (Recall that $J_X^\rs = J\bullet T^*X^\rs$.) Finally, recalling that $N^*_{X_{-1}}$ is the closure of $(-1)\cdot T^*X^\rs$ in $J_X$, let $\widetilde{N^*_{X_{-1}}}^{P,Y}$ be the closure of $(-1)\cdot \widetilde{T^*X}^{P,Y,\rs}$ in $\widetilde{J_X}^{P,Y}$. 

Explicitly, fix the pair $(Y,P)$, let $x_1\in Y\subset X$, let $V = T^*_{x_1} X$, and $V_P = V\cap \mu^{-1}(\u_P^\perp)$, as in Proposition \ref{propVP}, $v\in V_P^\rs$. Then $(-1)\cdot v \in T^*X$, the translate of $v$ under the action of the $(-1)$-section of $J$, lives over a point $x_2$ which also belongs to $Y$, because this action can be induced by the centralizer of $v$ in a Levi subgroup of $P$ (as in the proof of Lemma \ref{lemmaregular}). The point $x=(x_1,x_2)\in X\times X$ belongs to the closed orbit $X_{-1}$, by Lemma \ref{lemmaregular}, and setting $H_i$ for the stabilizer of $x_i$, and $V' = (\h_1+\h_2)^\perp\subset \g^*$, $V'_P=V'\cap \u_P^\perp$, we have 
$$ \widetilde{N^*_{X_{-1}}}^{P,Y} = V'_P\times^{(P\cap H_1\cap H_2)} G  \longrightarrow  V'\times^{(H_1\cap H_2)} G = N^*_{X_{-1}}.$$ 

Let $Y_2 = Y/U_P$, as before, and let $Y_{2,-1}$ be analog of $X_{-1}$ for $Y_2$ --- that is, the closed $L$-orbit on $Y_2\times Y_2$ which contains the projections of cotangent pairs $(v, (-1)v)$, where $v\in T^*Y_2^\rs$. If $\bar x$ is the image of $x$ in $Y_{2,-1}$, and $V_2'$ the fiber of $N^*_{Y_{2,-1}}(Y_2\times Y_2)$ over $\bar x$, then we have a quotient map $V'_P\to V'_2$.

\begin{proposition}\label{propVP2}
 In the setting above, the map $\widetilde{N^*_{X_{-1}}}^{P,Y} \to N^*_{X_{-1}}$ is surjective on $F$-points.
 
 The kernel of the map $V_P'\to V_2'$ is an isotropic subspace of $V'$ (with respect to the quadratic map $V'\to \c_X^*$) of dimension 
\begin{equation}\label{dimensions2}
 \dim\ker(V'_P\to V'_2) = \frac{\dim V' - \dim V'_2}{2}.
\end{equation}

 The quadratic space $V'$ is split (maximally isotropic) if and only if $V'_2$ is, which happens if and only if the quadratic space $V$ is.
\end{proposition}

Hence, by the equivalences of Proposition \ref{rationality}, the fibers of $N^*_{X_{-1}}$ are split quadratic spaces if and only if the stabilizer of one point on $X$ is split.

\begin{proof}
 The proof is identical to that of Proposition \ref{propVP2}; we just need to add how the property of $V'_2$ being split relates to the property of $V_2$ being split. But this is clear from considering the space $Y_2$ or, equivalently, $Y_{2,\ad}= Y_2/\mathcal Z(L)$, the latter being isomorphic to $T\backslash \PGL_2$ or to $\PGL_2$. In the both cases, there is an automorphism of rank $2$ of $Y_2$, which, applied to one copy of $Y_2$, interchanges the orbits $Y_2^\diag$ and $Y_{2,-1}\subset Y_2\times Y_2$; so, one, equivalently all fibers of the conormal bundle of the former are split if and only if one, equivalently all, fibers of the conormal bundle of the latter are.
\end{proof}

Finally:

\begin{lemma}\label{surjectiveonpoints}
 Assume that one, equivalently all, stabilizers of points on $X$ are split. Then the map $X\times X\to \C_X = A_X\sslash W_X$ is surjective on $F$-points.
\end{lemma}

\begin{proof}
 By Proposition \ref{rationality}, the stabilizers being split is equivalent to the map $T^*X\to \c_X^*$ being totally isotropic, which implies that it is surjective. This means that the projection $J^\rs\bullet T^*X^\rs\to J^\rs$ is surjective on $F$-points. The map $J^\rs\to A_X\sslash W_X$ is also surjective on $F$-points. Hence, the composition $\P J_X^\rs = \P(J^\rs\bullet T^*X^\rs)\to X\times X\to \C_X$ is surjective on $F$-points. 
\end{proof}

\begin{example}\label{pgl3-second} 
 Let us consider the case of the variety $X=\GL_2\backslash \PGL_3 =$ the variety of decompositions $\Ga^3=V_2\oplus V_1$, that we already saw in Example \ref{pgl3}. Letting $x_1$ be the decomposition $\left<e_1,e_2\right>\oplus\left<e_3\right>$, and $P=$ the stabilizer of the plane $\left<e_2,e_3\right>$, we will get $x_2=$ a decomposition $\left<e_1 + ce_2, e_3\right> \oplus \left<e_2\right>$, with the scalar $c$ depending on the chosen cotangent vector. Then we see that $H_1\cap H_2  = P\cap H_1\cap H_2 \simeq\Gm$, $V'=V_P'=(\h_1+\h_2)^\perp = \Std_1 \oplus \Std_1^*$, where $\Std_1$ denotes the standard one-dimensional representation of $\Gm$, and $V'$ intersects the nilpotent cone in $\g^*$ along \emph{irregular} orbits only. The reader should compare this with Example \ref{pgl3}, where we saw that the nilpotent limit of the action of the $(-1)$-section of $J$ on regular semisimple vectors does not exist, but it does exist at the exceptional divisor of the blowup along the (irregular nilpotent) divisor $\Std\oplus \Std^*$. This blowup, with the strict transform of the nilpotent divisor removed, is isomorphic to the conormal bundle $N^*_{X_{-1}}$. 
\end{example}

\section{Integration formula} \label{sec:integration}

The goal of this section is to prove the Theorem \ref{thmintegrationformula} below. To formulate it, we will need a way to fix measures compatibly on orbits which are isomorphic over the algebraic closure:

\begin{definition}\label{defmeasures}
Let $H_1, H_2\subset G$ be two subgroups that are conjugate over the algebraic closure, and such that the normalizer of $H_i$ acts trivially on the top exterior power $\bigwedge^\top \mathfrak h_i$ (in particular, $H_i$ is unimodular). Two $G$-invariant measures $\mu_1, \mu_2$ on the spaces $H_1\backslash G$, $H_2\backslash G$, respectively, will be said to be \emph{compatible}, if they can be presented as $\mu_i=c\cdot |\omega_i|$ for some invariant volume forms $\omega_i$, $i=1,2$, and the same scalar $c$, such that, over the algebraic closure, $\omega_2$ is conjugate to $\epsilon\omega_1$, for some $\epsilon \in \bar F$ with $|\epsilon|=1$. 
\end{definition}

Recall that we have fixed a Haar measure on $F$, so that the absolute value of a volume form is a well-defined measure. We will apply this to the following setting: Consider the diagonal action of $G$ on $X\times X$. By Proposition \ref{genericisomorphism} and Corollary \ref{corregulardense}, the $G$-stabilizers of points over any $c \in \mathring\C_X$ (the complement of $[\pm 1]$ in $\C_X=A_X\sslash W_X$) are all conjugate, over the algebraic closure, to the kernel of the map $L(X)\twoheadrightarrow A_X$, where $L(X)$ is a Levi subgroup of a parabolic of type $P(X)$.

\begin{theorem}\label{thmintegrationformula}
Let $\omega$ be a nonzero $G\times G$-invariant volume form on $X\times X$ defining an invariant measure $|\omega|$.  Fix compatible $G$-invariant measures $d\dot{g}$ (Definition \ref{defmeasures}) on all $G$-orbits over $\mathring\C_X$.

Then, identifying $\C_X\simeq \mathbbm A^1$ and letting $c_{\pm 1}$ be the coordinates of the points $[\pm 1]$, there is an additive Haar measure $dc$ on $\mathbbm A^1$ such that the following integration formula holds:
 \begin{equation}\label{integrationformula} \int_{X\times X} \Phi(x) |\omega|(x) = \int_{\C_X} |c-c_1|^{\frac{d_1}{2}-1} |c-c_{-1}|^{\frac{d_{-1}}{2}-1}  \left(\int_{(X\times X)_c} \Phi(\dot{g}) d\dot{g} \right) dc.
\end{equation}

Here, $d_1=\dim X=\codim X_1$, and 
\begin{equation}\label{codimensionseq}
d_{-1}=\codim X_{-1} =  2\epsilon \left<\rho_{P(X)},\check\gamma\right>-d_1+2,
\end{equation}
where $\check\gamma$ is the spherical coroot, and
$$ \epsilon = \begin{cases} 1, \mbox{ for roots of type $T$ (dual group $\SL_2$)};\\
2, \mbox{ for roots of type $G$ (dual group $\PGL_2$).}
\end{cases}$$
Moreover, in the case of type $G$ we have 
\begin{equation}\label{equalcodim}d_1=d_{-1} = \left<2\rho_{P(X)},\check\gamma\right>+1.
\end{equation}

\end{theorem}

Notice that the formula for $\codim X_{-1}$ is new, and will be proven as a corollary of the integration formula.

There are also analogous integration formulas for the normal/conormal bundles of the orbits $X_1$ and $X_{-1}$. Notice that, fixing the isomorphism $\c_X^*=\mathbbm A^1$, by the nondegenerate quadratic forms obtained by the invariant-theoretic quotients 
$$N_{X_1}^*\to  \c_X^* \leftarrow N_{X_{-1}}^* $$ (see Proposition \ref{genericisomorphism}), the normal and conormal bundles are $G$-equivariantly isomorphic. Again, $G$-stabilizers of points on $X\times X$ over any $\xi\ne 0\in \c_X^*$ are all conjugate, over the algebraic closure, to the kernel $L_1$ of the map $L(X)\twoheadrightarrow A_X$.

\begin{theorem}\label{thmintegrationformula-linear}
Let $\omega$ be a nonzero $G$-invariant volume form on $N_{X_{\pm 1}}^*$ which restricts to Haar measures on the fibers. Fix compatible $G$-invariant measures $d\dot{g}$ on all $G$-orbits over $\mathring\c_X^*$.

Then, there is an additive Haar measure $d\xi$ on $\c_X^*\simeq \mathbbm A^1$ such that the following integration formula holds:
 \begin{equation} \int_{N^*_{X_{\pm 1}}} \Phi(x) |\omega|(x) = \int_{\c_X^*} |\xi|^{\frac{d_{\pm 1}}{2}-1} \left(\int_{(N^*_{X_{\pm 1}})_\xi} \Phi(\dot{g}) d\dot{g} \right) d\xi.
\end{equation}
\end{theorem}

The proof of Theorem \ref{thmintegrationformula-linear} is completely analogous to that of Theorem \ref{thmintegrationformula}, and therefore I will only present that of Theorem \ref{thmintegrationformula}, leaving the reformulation for the other to the reader.

\subsection{Pullback to the polarization}

Consider the map $\P J_X\xrightarrow{\R} X\times X \to A_X\sslash W_X$. Recall that $J_X^\rs = J\bullet T^*X^\rs$; we let 
$$\widehat{J_X}^{\bullet,\rs} = J\bullet \widehat{T^*X}^{\bullet,\rs} = A_X \times \widehat{T^*X}^{\bullet,\rs} \subset J_X\times_{\c_X^*} \a_X^*,$$
where $\widehat{T^*X}^\bullet$ is the distinguished irreducible component of the polarized cotangent bundle that was defined in \S \ref{ssmomentmap}.

We have a commutative diagram
$$ \xymatrix{
\P \widehat{J_X}^{\bullet,\rs} \ar[r]^{p}\ar[d] &  \P J_X \ar[r]^{\R}\ar[d] & X\times X \ar[dl]\\
A_X \ar[r] & A_X\sslash W_X.}$$

Recall that we denote by $[\pm 1]_{\P J_X}$ the exceptional divisors of the blowup $\R$, and by $d_{\pm 1}$ the codimensions of the orbits $X_{\pm 1}$.
\begin{lemma}\label{canonicaldivisor}
 Let $K_{X\times X}$ be the canonical bundle on $X\times X$. Then 
 $$\R^* K_{X\times X} = K_{\P J_X} ((d_1-1)[1]_{\P J_X} + (d_{-1}-1)[-1]_{\P J_X}).$$ 
\end{lemma}

\begin{proof}
 This is immediate from the characterization of $\P J_X$ as the blowup of $X\times X$ at the two divisors $X_1=X^\diag$ and $X_{-1}$ (Proposition \ref{isblowup}).
\end{proof}

Hence, if $\omega$ is a nonzero, $G\times G$-invariant volume form on $X\times X$, the divisor of its pullback to $\P J_X$ is 
$$[\R^*\omega]=(d_1-1) \left[ [1]_{\P J_X}\right] + (d_{-1}-1)\left[[-1]_{\P J_X}\right],$$ 
where $d_1 = \codim (X_1) = \dim X$, and $d_{-1} = \codim (X_{-1})$. 

Set $Y=\P \widehat{J_X}^{\bullet,\rs}$, for notational simplicity. The map $p:Y \to \P J_X $ is an \'etale $\ZZ/2$-cover onto its image, and notice that $\widehat{J_X}^{\bullet,\rs} = A_X \times \widehat{T^*X}^{\bullet,\rs}$, canonically. Thus, setting $\hat\R=p\circ\R$, we have 
\begin{equation}\label{omegatopolar}[\hat\R^*\omega] = (d_1-1) [Y_1] + (d_{-1}-1)[Y_{-1}].
\end{equation}
where $Y_{\pm 1} = [\pm 1]_Y\subset Y$.

Fix a pair $(x,B)\in (X\times \B)^\circ$, defining an embedding $\mathring\a_X^*\to \widehat{T^*X}^{\bullet,\rs}$ by Knop's section $\hat\kappa_X$ (\S \ref{ssmomentmap}), and let $L_1$ be the stabilizer of the points in the image. If $L$ denotes the centralizer of the image of such a point under the polarized moment map to $\hat\g_X^*$, identified with the Levi quotient $L(X)$ of $P(X)$, then $L_1\supset [L,L]$, and $L/L_1\simeq A_X$, canonically because of the polarization. The action map identifies $$\widehat{T^*X}^{\bullet,\rs} \simeq \mathring\a_X^* \times L_1\backslash G.$$ 
Hence: 
\begin{equation}\label{polarisom}\P\widehat{J_X}^{\bullet,\rs}  = A_X \times L_1\backslash G.
\end{equation}

Fix an invariant volume form $\omega_{L_1\backslash G}$ on $L_1\backslash G$. Then, by \eqref{omegatopolar} and the fact that the only nowhere vanishing regular functions on a torus are characters, there is a Haar volume form $\omega_{A_X}$ on $A_X$ and an $m\in \ZZ$ such that
\begin{equation}\label{formonpolarization} \hat\R^*\omega = (a-1)^{d_1-1} (a+1)^{d_{-1}-1} a^m \cdot \omega_{A_X}\wedge \omega_{L_1\backslash G}
\end{equation}
 under \eqref{polarisom}, where we have identified $A_X\simeq \Gm$ to fix a coordinate $a$. On the other hand, this has to be invariant under the $W_X$-Galois action $a\mapsto a^{-1}$, hence
$$ m = 1- \frac{d_1+d_{-1}}{2}.$$

\subsection{Descent to $X\times X$}

The integration formula \eqref{integrationformula}, now, follows from \eqref{formonpolarization} by descending to $A_X\sslash W_X \simeq \mathbbm A^1$: Fix a coordinate $c$ on that space, with $c_{\pm 1}$ the coordinates of the points $[\pm 1]$. In a sufficiently small neighborhood $U$ of any point of $A_X\sslash W_X\smallsetminus\{[\pm 1]\}$ the stabilizers of all points are conjugate to a group $L_1'$ which is conjugate over the algebraic closure to $L_1$, and the preimage of $U$ in $X\times X$ is, $G$-equivariantly isomorphic to $U\times L_1'\backslash G$. Thus, fixing the compatible measures $d\dot{g}$ on the $G$-orbits as in the statement of the theorem, there is an integration formula of the form:
$$ \int_{X\times X} |\omega| = \int_{A_X\sslash W_X} \varphi(c) \int_{(X\times X)_c} d\dot{g} dc$$
for some nonnegative measurable function $\varphi$ on $A_X\sslash W_X\smallsetminus\{[\pm 1]\}$, and some additive Haar measure $dc = |\omega_{\mathbbm A^1}|$ on $\mathbbm A^1$. On the other hand, writing, in such a neighborhood $c$, the measure $d\dot{g}$ on $L_1'\backslash G$ as $|\omega'|$, for some invariant volume form $\omega'$, we see by applying \eqref{formonpolarization} over a suitable algebraic extension of $F$ that the pullback of $\omega_{\mathbbm A^1}\wedge \omega'$ to $\P \widehat{J_X}^{\bullet,\rs}$ has to be a multiple of $(a-1)^{d_1-1} (a+1)^{d_{-1}-1} a^m \cdot \omega_{A_X}\wedge \omega_{L_1\backslash G}$ by a rational function $f(a)^{-1}$ with $|f(a)| = \varphi(c(a))$. Without loss of generality, the pullback of $\omega'$ is equal to $\omega_{L_1\backslash G}$, and an elementary calculation, choosing, for example, the isomorphisms and the map $A_X=\Gm \to A_X\sslash W_X=\mathbbm A^1$ given by $c=a+a^{-1}$,
shows that, up to a scalar that we can take to be $1$ by scaling the forms, 
$$|\hat\R^*\omega_{\mathbbm A^1}| = |c-c_1|^\frac{1}{2}|c-c_{-1}|^\frac{1}{2} |\omega_{A_X}|,$$ 
and $|a-(\pm 1)|^2\cdot {|a|^{-1}} = |c - c_{\pm 1}|$ (up to a fixed scalar), hence:
$$ \varphi(c) = |c - c_1|^{\frac{d_1}{2}-1} |c - c_{-1}|^{\frac{d_{-1}}{2}-1}$$
for a suitable $dc$.

\subsection{Degeneration} \label{ssdegeneration}

We have proven the integration formula of Theorem \ref{thmintegrationformula}, except for the determination of the codimension $d_{-1}$ of the orbit $X_{-1}$. In this subsection we will prove the codimension formula \eqref{codimensionseq}:
$$ (d_1-1) + (d_{-1}-1) = 2\epsilon \left<\rho_{P(X)},\check\gamma\right>,$$
where $\check\gamma$ is the spherical coroot, and
$$ \epsilon = \begin{cases} 1, \mbox{ for roots of type $T$ (dual group $\SL_2$)};\\
2, \mbox{ for roots of type $G$ (dual group $\PGL_2$).}
\end{cases}$$
with  $d_1=d_{-1} = \left<2\rho_{P(X)},\check\gamma\right>+1$ in the case of type $G$.

To prove this, we degenerate $X$ to its \emph{boundary degeneration} $X_\emptyset$, which is a horospherical variety. The codimension formula will follow by comparing the integration formula \eqref{integrationformula} to the corresponding formula for $X_\emptyset\times X_\emptyset$, which is very easy to compute.

More precisely, consider the decomposition of the coordinate ring of $X$ as a $G$-module into irreducibles:
$$ F[X] = \bigoplus_{\lambda \in \Lambda_X^+} F[X]_\lambda.$$
The indices $\lambda$ denote, here, the highest weight of the representation, and $\Lambda_X^+\subset \Lambda_X:=\Hom(A_X,\Gm)\simeq \ZZ$ is the monoid of weights such that the corresponding rational $B$-eigenfunction is regular; from now on we identify it with $\mathbb N$. The above decomposition is not an algebra grading, a fact that is known to be equivalent to the fact that $W_X\ne 1$ (see \cite{KnAut}). Instead, it corresponds to an algebra filtration
$$ \mathcal F_\lambda = \bigoplus_{\nu\le \lambda} F[X]_\nu,$$
i.e., $\mathcal F_\lambda \cdot \mathcal F_\mu\subset \mathcal F_{\lambda+\mu}$.
The Rees family \cite[\S 4]{Popov}:
$$ F[\mathscr X] := \bigoplus_\lambda \mathcal F_\lambda \cdot t^\lambda \subset F[X][t]$$
defines an affine $G$-variety $\mathscr X$ over $\mathbbm A^1=\Spec F[t]$, together with an action of $\Gm$ that extends to a morphism $\tau: X\times \mathbbm A^1\to \mathscr X$; more canonically, $\Gm = \Spec F[\Lambda_X] = A_X$, and the base $\mathbbm A^1$ of this family is the affine embedding $\overline{A_X}\supset A_X$ on which elements of $\Lambda_X^+$ extend to regular functions, so the defining morphism is
\begin{equation}\label{tau}
\xymatrix{ X\times \overline{A_X}\ar[dr]\ar[rr]^\tau && \mathscr X\ar[dl]^\pi\\ & \overline{A_X},}
\end{equation}
and extends to a canonical action of $A_X\times G$ on $\mathscr X$ over $\overline{A_X}$.

By \cite{Popov} (see also \cite[\S 2.5]{SV}), the special fiber \label{refboundary} $X_\emptyset:=\mathscr X_0$ is an affine horospherical $G$-variety\footnote{In other references, $X_\emptyset$ denotes just the open $G$-orbit in $X_\emptyset$. Here, I have found it more convenient to use $X_\emptyset$ for the affine degeneration, and $X_\emptyset^\bullet$ for its open $G$-orbit.} with $P(X_\emptyset) = P(X)$; its open $G$-orbit $X_\emptyset^\bullet$ is isomorphic to $U^-L_1\backslash G$, where $U^-$ is the unipotent radical of a parabolic opposite to $P(X)$, and $L_1$ is the kernel of the map $L(X)\twoheadrightarrow A_X$, as above. Notice that the image of the defining morphism $\tau$ lies in the complement of the open orbit over $0\in \overline{A_X}$.

It is known that there is a family of $G$-invariant volume forms on the homogeneous parts of the fibers of $\mathscr X\to \mathbbm A^1$ that is everywhere nonvanishing, cf.\ \cite[\S 4.2]{SV}. More precisely, let us fix a parabolic $P$ in the class of $P(X)$, and let $\mathring X$ be the open $P$- (and Borel) orbit. Restricting to $U$-invariants, the above decomposition becomes a grading
$$ F[X]^U = \bigoplus_\lambda F[X]^U_\lambda,$$
with $F[X]^U_\lambda=$ the (one-dimensional) highest weight subspace of $F[X]_\lambda$, on which $P(X)$ acts by the character $\lambda$ of the quotient $A_X$. Correspondingly, the family $\mathscr X\sslash U$ becomes constant:
\begin{equation}\label{constantmodU}
 \sigma:\mathscr X\sslash U \xrightarrow\sim X\sslash U \times \overline{A_X}, 
\end{equation}
but it can be seen from the definitions that this isomorphism is related to the one (call it $\tau_U$) that we obtain by descending the defining map $\tau$ of \eqref{tau} by:
\begin{equation}\label{comparison} \sigma^{-1}(\bar x, t) = \tau_U(\bar x\cdot t^{-1}, t),\end{equation}
where we have used the canonical action of $A_X$ on $X\sslash U$; that is, the action of $A_X$ on $\mathscr X$ descends to the action
\begin{equation}\label{comparison2} a\cdot (\bar x, t) = (\bar x\cdot a, at)\end{equation}
on $X\sslash U \times \overline{A_X}$.

If $\mathring{\mathscr X}$ is the union of open $P$-orbits on the various fibers, the restriction of $\sigma$ defines an isomorphism
$$ \mathring{\mathscr X}/U \xrightarrow\sim \mathring X/U \times \overline{A_X},$$
which, by the local structure theorem \cite[Th\'eor\`eme 1.4]{BLV}, \cite[Theorem 2.3]{KnMotion}, can be lifted $P=L\cdot U$-equivariantly:
\begin{equation}\label{tildesigma}\tilde\sigma: \mathring{\mathscr X} \simeq  S \times U \times \overline{A_X},
\end{equation}
for some Levi subgroup $L\subset P$, with $S\simeq \mathring X/U$ an $L$-stable subvariety of $X$, acting by conjugation on $U$ and via the quotient $A_X$, simply transitively, on $S$. More precisely, $L$ is the centralizer of the image of an element $v\in T^*\mathring X^{\u^\perp,\rs}$ under the moment map, and $S$ is a ``flat'', that is, the $L$-orbit of the image of $v$ on $X$.

From now on, by abuse of notation, any reference to $\mathscr X$ should be taken to refer to the smooth locus of the map $\pi:\mathscr X\to \overline{A_X}$.
Consider $\Omega:=\Omega_{\mathscr X/ \mathbbm A^1}=$ the relative cotangent sheaf of $\mathscr X\xrightarrow\pi \mathbbm A^1$. Restricted to any fiber, it is canonically identified with the cotangent bundle of that fiber. Its top exterior power, $\bigwedge^\top \Omega$, restricts to the bundles of volume forms on the fibers (of the smooth locus). 

\begin{lemma}
There is a $G$-invariant section $\omega$ of $\bigwedge^\top \Omega$ which restricts to a nonvanishing volume form on each fiber. Moreover, such a form is an $A_X$-eigenform satisfying, for every $a\in A_X$,
\begin{equation}\label{omegaeigen} a^*\omega_{\mathscr X} = e^{2\rho_{P(X)}}(a) \cdot \omega_{\mathscr X}.
\end{equation}
\end{lemma}

Notice that we use exponential notation for the character $2\rho_{P(X)}=$the sum of roots in the unipotent radical of $P(X)$, since we use additive notation for roots.

\begin{proof}
A nonzero, $P$-invariant volume form on $X/U \times U$ pulls back by \eqref{tildesigma} to a $P$-invariant section $\omega_{\mathring{\mathscr X}}$  of the line bundle $\bigwedge^\top \Omega$ on $\mathring{\mathscr X}$. If $X$ admits a $G$-invariant measure (as is the case for affine homogeneous spaces), so does $X_\emptyset$ \cite[\S 4.2]{SV}, and the $g$-pullback of $\omega_{\mathring{\mathscr X}}$, for every $g\in G$, is a $gPg^{-1}$-invariant section of $\Omega$ on $\mathring{\mathscr X}g^{-1}$, which coincides with $\omega_{\mathring{\mathscr X}}$ on the intersection with $\mathring{\mathscr X}$; thus, these translates glue to a global section $\omega_{\mathscr X}$ of $\Omega$ over the union $\mathscr X^\bullet$ of all open $G$-orbits on the fibers.

This section restricts, by construction, to a nonvanishing volume form on each fiber. Notice, also, that any other such form is a multiple of $\omega_{\mathscr X}$ by a nowhere vanishing regular function on $\overline{A_X}\simeq \mathbbm A^1$, hence by a scalar.

Regarding the action of $A_X$, it is enough to prove \eqref{omegaeigen} for the restriction of $\omega_{\mathscr X}$ to $\mathring{\mathscr X}$. In terms of the isomorphism \eqref{tildesigma}, the $P$-invariant form is given by 
$$\omega_{\mathring{\mathscr X}}(s,u,t) = e^{2\rho_{P}}(s) \cdot \omega_S(s) \wedge \omega_U(u), $$
where $\omega_S$ is an $A_X$-invariant volume form on $S$, $\omega_U$ is a $U$-invariant volume form on $U$, and we have identified $S\simeq A_X$ by choosing a base point. The action of $A_X$ on $\mathring{\mathscr X}$ is given by \eqref{comparison2} on $S\times \overline{A_X}$, and trivial action on $U$, therefore this form is $e^{2\rho_{P(X)}}$-equivariant.
\end{proof}

Now we move to the space $\widetilde{{\mathscr X}}:=\mathscr X\times_{\overline{A_X}}\mathscr X$. Again, we only work over the smooth locus of the morphism to $\overline{A_X}$. The tensor product of $\omega_{\mathscr X}$ with itself gives rise to a section $\omega_{\widetilde{\mathscr X}}$ of the top exterior power of the relative cotangent bundle of $\widetilde{\mathscr X}\to \overline{A_X}$, which restricts to an invariant, nonvanishing, $G\times G$-invariant volume form on the open orbit of each fiber. 

\begin{proposition}\label{GITfamily}
 There is an isomorphism $\widetilde{{\mathscr X}}\sslash G \simeq \mathbbm A^1\times \overline{A_X}$ over $\overline{A_X}$.
\end{proposition}

\begin{proof}
 Recall the heighest weight decomposition $F[X] = \bigoplus_\lambda F[X]_\lambda$, where $\lambda$ ranges in a monoid $\Lambda_X^+\simeq \mathbb N$ of weights of $A_X$. Notice that the highest-weight modules $F[X]_\lambda$ are necessarily self-dual; indeed, twisting the action of $G$ on $X=H\backslash G$ by a Chevalley involution does not change its isomorphism class as a $G$-variety (because $H$ is reductive), hence preserves the monoid $\Lambda_X^+$; because $\Lambda_X^+$, it acts trivially on it. Thus, $(F[X]_\lambda\otimes F[X]_\lambda)^G=F$, where $G$ here acts diagonally.
 
 This gives the structure of a graded vector space to
 $$ F[X\times X]^G = \bigoplus_{\lambda, \mu} (F[X]_\lambda\otimes F[X]_\mu)^G = \bigoplus_{\lambda} (F[X]_\lambda\otimes F[X]_\lambda)^G =\bigoplus_{\lambda} F,$$
 which also corresponds to a filtration of rings, with associated graded $\gr F[X\times X]^G = F[\overline{T_1}]$, where $\overline{T_1}$ is the image, in the grading, of a nonzero element $T_1$ of $F$ in the copy labelled by the first nontrivial element of $\Lambda_X^+$. 
 
 But this shows that $T_1\in F[X\times X]^G$ generates the ring freely, thus, $F[X\times X]^G\simeq F[T_1]$.
 
 Moving now to the coordinate ring $F[\mathscr X]\subset F[X][t]$ of the Rees family, this argument shows that 
 $$ F[\widetilde{\mathscr X}]^G = F[T_1 t^2, t],$$
 hence $\widetilde{\mathscr X}\sslash G \simeq \mathbbm A^2 = \mathbbm A^1 \times \overline{A_X}$.
\end{proof}

Notice that at $t=0$, this specializes to an isomorphism $X_\emptyset\times X_\emptyset \sslash G \simeq \mathbbm A^1$ which is $A_X$-equivariant when $A_X$ acts by the \emph{square} of the generator $\lambda_1$ of $\Lambda_X^+$ on $\mathbbm A^1$ (because the action of $A_X$ on $\widetilde{\mathscr X}$ restricts to its diagonal action on the two copies of the special fiber $X_\emptyset\times X_\emptyset$). More generally, the action of $A_X$ on $\mathbbm A^1 \times \overline{A_X}=\Spec F[y, t]$, where $y=T_1t^2$ as in the proof above, is given by 
$$ a\cdot (y,t) = (\lambda_1^2(a) y, \lambda_1(a) t).$$

As mentioned, the restriction of the form $\omega_{\widetilde{\mathscr X}}$ to the fiber over any $t\in \overline{A_X}$ is a $G\times G$-invariant, nonzero volume form $\omega_t$; on the special fiber, it satisfies the integration formula:
\begin{equation}\label{integration-degen}\int_{X_\emptyset\times X_\emptyset} |\omega_0| = \int_{\mathbbm A^1} |c|^{2\epsilon \left<\rho_{P(X)},\check\gamma\right>}  \left(\int_{L_1\backslash G} |\omega_{L_1\backslash G}|\right) dc,
\end{equation}
where $\epsilon$ is as in Theorem \ref{thmintegrationformula}. Indeed, the special fiber contains an open dense subset which is $A_X\times G$-equivariantly isomorphic to $A_X\times L_1\backslash G$, and which corresponds to the open Bruhat cell under the isomorphism $X_\emptyset^\bullet\times X_\emptyset^\bullet/G = H_\emptyset\backslash G/H_\emptyset$, where $H_\emptyset \simeq \ker(P(X)^-\to A_X)$, with $P(X)^-$ opposite to $P(X)$. The parabolic $P(X)^-$ is actually conjugate to $P(X)$: indeed, if $X=H\backslash G$, since $H$ is reductive there is a Chevalley involution of $G$ which fixes $H$, and hence preserves the isomorphism class of $X$ --- but this turns the class of $P(X)$ to the class of $P(X)^-$. Thus, the integration formula for the open Bruhat cell with respect to $P(X)^-$ reads:
$$ \int_G \Phi(g) dg = \int_{A_X} \int_{H_\emptyset\times U_{P(X)}^-} \Phi(u_1 aw u_2) d(u_1, u_2) \cdot |e^{2\rho_{P(X)}}(a)| da,$$
where $w$ is the longest element of the Weyl group, and this easily translates to \eqref{integration-degen}. Here, we need to take into account that there is an isomorphism 
$$H_\emptyset\backslash G\sslash H_\emptyset \xrightarrow\sim \mathbbm A^1$$ 
which pulls back to the character $\frac{\gamma}{\epsilon}$ (a generator for $\Lambda_X=\Hom(A_X,\Gm)$)
under the sequence of maps
$$  A_X\to H_\emptyset A_X w H_\emptyset \hookrightarrow G\to H_\emptyset\backslash G\sslash H_\emptyset \xrightarrow\sim \mathbbm A^1.$$
Hence, the inverse of this sequence of maps (restricted to $\Gm$) is given by the cocharacter $\epsilon\check\gamma$.

On the other hand, consider the integration formula \eqref{integrationformula},  taking into account that the points $c_1,c_{-1}$ on $X\times X\sslash G$, expressed now in the coordinate $T_1$ as above, when we vary the parameter $t\ne 0$ become $c_{\pm 1} t^2$ in the coordinate $T_1 t^2$. The limit as $t\to 0$ must coincide with the integration formula \eqref{integration-degen} on $X_\emptyset\times X_\emptyset$, proving the codimension formula \eqref{codimensionseq}. 

Finally, for spherical roots of type $G$ we have, by Proposition \ref{correctrepresentative}, a nontrivial $G$-automorphism of $X$ of order $2$.
Applied to the first copy of $X$ in $X\times X$, this automorphism does not preserve the diagonal $X_1$, hence has to interchange it with the unique other semisimple $G$-orbit which can have codimension larger than one, that is, with $X_{-1}$.
This completes the proof of Theorem \ref{thmintegrationformula}.

\section{Schwartz measures} \label{sec:Schwartzmeasures}

We are ready to consider the pushforward of Schwartz measures:
\begin{equation}\label{pushforward} \mathcal S(X\times X)\to \Meas(\C_X),
\end{equation}
whose image we have denoted by $\mathcal S(X\times X/G)$. From now on, we assume that $X$ is not only a ``correct representative'' in its class modulo $G$-automorphisms (Definition \ref{defcorrectrep}), but also that stabilizers of points on $X$ are split; thus, by Lemma \ref{surjectiveonpoints}, the map $X\times X\to \C_X$ is surjective on $F$-points.

In this section we will obtain as much information as possible from abstract principles about the space $\mathcal S(X\times X/G)$, using the blowup $\P J_X$. We use the blowup in the way that it is used in Igusa integrals: as a resolution of the map $X\times X\to \C_X$, in the sense that preimages of points are normal crossings divisors, see Proposition \ref{isblowup}. 

We will actually be working mainly with the linearizations of this $G$-space. 
The main result of the section is Theorem \ref{pushfthm}, leaving us only a certain linear combination of scalars to compute in the next section. Strictly speaking, the techniques of the next section are sufficient to obtain the main results, Theorems \ref{Xtheorem} and \ref{maintheorem}, but using the resolution puts the results into a conceptual context, up to the computation of a linear combination of coefficients.

\subsection{Generalities on Schwartz measures}\label{ssgeneralities-Schwartz}

Before we proceed, I recall some concepts, and introduce some notation, related to cosheaves of Schwartz measures; more details can be found in \cite{AGSchwartz}.

Let $Z$ be a smooth variety, and $D\subset Z$ a divisor. We let $C^\infty(\bullet,D)$ denote the sheaf of functions on the $F$-points of $Z$ which, locally, are of the form $\Phi(z) |\epsilon_D(z)|$, where $\Phi$ is a smooth function and $\epsilon_D$ is a local generator for the divisor $D$. Informally, we consider such functions as ``smooth sections of the complex line bundle $|\mathcal L_D|$ associated to $D$''.

Consider the restricted topology of semialgebraic sets on the $F$-points of $Z$. (The $F$-points of Zariski open subsets will be enough, for our purposes.) We can define a cosheaf $\mathcal S(\bullet, D)$ of Schwartz measures valued in $|\mathcal L_D|$. In the Archimedean case, its sections are linear combinations of measures on semialgebraic open subsets which can be written as products $\Phi\cdot \omega$, where $\omega$ is a nowhere vanishing Nash (=smooth semialgebraic) density on such an open subset, and $\Phi$ is \emph{a section of $C^\infty(\bullet,D)$ of rapid decay} (together with its derivatives). In other words, sections of $\mathcal S(\bullet, D)$ are generated by linear combinations of measures supported on open subsets which admit (semialgebraic) coordinates $(x_1, \dots, x_n)$, and where the divisor $D$ is represented by a polynomial function $f_D$, and in such coordinates are of the form 
\begin{equation}\label{Schwartzonchart}
|f_D|\cdot \Phi\cdot dx_1\cdots dx_n,
\end{equation}
 where $\Phi$ is a Schwartz function on the given chart. For every open $U\subset Z$, the space $\mathcal S(U, D)$ is a Fr\'echet space; if $U$ admits such a chart, the topology is defined by the seminorms 
$$|f_D|\cdot \Phi\cdot dx_1\cdots dx_n \mapsto \sup |T\Phi|,$$
where $T$ ranges over all smooth semialgebraic differential operators defined on that chart. (In the case of affine algebraic sets, one can consider just algebraic differential operators.)
In the non-Archimedean case, sections of $\mathcal S(\bullet, D)$ are simply linear combinations of measures which, locally on a chart, can be written as \eqref{Schwartzonchart}, with $\Phi$ compactly supported. Note that they are not necessarily smooth as measures on $Z$, because of the factor $|f_D|$.

For a closed subset $Y\subset Z$ (``closed'' means semialgebraic, again, but the reader can restrict their attention to Zariski closed) we define the \emph{stalk} $\mathcal S_Y(\bullet,D)$ \label{refstalkfiber} as the cosheaf on $Z$, supported on $Y$, whose sections over an open $U\subset Z$ are the quotient 
\begin{equation}\mathcal S_Y(U,D)= \mathcal S(U,D)/\mathcal S(U\smallsetminus Y,D).
\end{equation}
 The \emph{fiber} $\overline{\mathcal S_Y(\bullet,D)}$ is the cosheaf whose sections over $U\subset Z$ are the quotient 
 \begin{equation} \overline{\mathcal S_Y(U,D)} = \mathcal S(U,D)/C^\infty_\temp(U,[Y])\mathcal S(U,D),
 \end{equation}
 where $C^\infty_\temp(U,[Y])$ denotes the ideal of those tempered (i.e., of polynomial growth together with their polynomial derivatives) smooth functions that vanish on $Y$. In the non-Archimedean case, the natural map $\mathcal S_Y(\bullet,D)\to \overline{\mathcal S_Y(\bullet,D)}$ is, of course, an isomorphism.

Our analysis of the pushforward \eqref{pushforward} starts from the following:

\begin{lemma}
 Let $Z\to Y$ be a smooth map of smooth varieties which is surjective on $F$-points. Then the pushforward of $\mathcal S(Z)$ is equal to $\mathcal S(Y)$.
\end{lemma}

\begin{proof}
 This is standard, see, e.g., \cite[Proposition 3.1.2]{SaStacks}.
\end{proof}

\begin{corollary}
 Let $U\subset X\times X$ be the preimage of $\mathring\C_X:=\C_X\smallsetminus\{[1], [-1]\}$. Then the pushforward of $\mathcal S(U)$ is the space of Schwartz measures $\mathcal S(\mathring\C_X)$. 
\end{corollary}

\begin{proof}
 Indeed, the map $X\times X\to \C_X$ is smooth over $\mathring\C_X$ by Proposition \ref{genericisomorphism}, and the map is surjective on $F$-points by Lemma \ref{surjectiveonpoints}.
\end{proof}

Thus, our remaining task is to determine the behavior of the elements of $\mathcal S(X\times X/G)$ close to the points $[\pm 1]\in \C_X$. To this end, we can linearize the problem:  Let $x\in X_{\pm 1}$, with $H_{\pm 1}$ its stabilizer in $G$ and $V_{\pm 1}$ its fiber in the conormal bundle $N^*_{X_\pm 1}$, so that $V_{\pm 1}\sslash H_{\pm 1} = N^*_{X_{\pm 1}}\sslash G = \c_X^*$, by Proposition \ref{GITquotients}. We let $\mathcal S(V_{\pm 1}/H_{\pm 1})$ be the pushforward of $\mathcal S(V_{\pm 1})$ under the map $V_{\pm 1}\to \c_X^*$.

\begin{proposition}\label{linearization}
There is an $F$-analytic isomorphism between a neighborhood $U_1$ of $0\in \c_X^*(F)$ (in the Hausdorff topology on $F$-points) and a neighborhood $U_2$ of $[\pm 1]\in \C_X(F)$, such that the space of restrictions to $U_1$ of elements of $\mathcal S(V_{\pm 1}/H_{\pm 1})$ is equal, under this isomorphism, to the space of restrictions to $U_2$ of the pushforwards of elements of $\mathcal S(X\times X)$ supported in a certain $G(F)$-stable neighborhood of $x$.
\end{proposition}

The restriction to a $G(F)$-stable neighborhood of $x$ is because the map $G(F)\to X_{\pm 1}(F)=H_{\pm 1}\backslash G(F)$ sending $g$ to $x\cdot g$ may not be surjective on $F$-points. Eventually, as we will see, the normal fibers $V_{\pm 1}$ of all points on $X_{\pm 1}$ contribute the same germs of pushforward measures, so this detail will not matter.

\begin{proof}
 The pushforward map $\mathcal S(V_{\pm 1})\to \mathcal S(V_{\pm 1}/H_{\pm 1})$ factors through the $H_{\pm 1}$-coinvariants of $\mathcal S(V_{\pm 1})$, and similarly the pushforward map $\mathcal S(X\times X)\to \mathcal S(X\times X/G)$ factors through the $G$-coinvariants of $\mathcal S(X\times X)$.
 
 Fix an isomorphism $\c_X^*\simeq \mathbbm A^1$ and use the resulting nondegenerate quadratic form (Proposition \ref{genericisomorphism}) $V_{\pm 1}\to \c_X^*\simeq\mathbbm A^1$ to identify $V_{\pm 1}$ with its linear dual, the fiber of the \emph{normal} bundle. 
 
 By Luna's \'etale slice theorem \cite{Lunaslice}, there is an $H_{\pm 1}$-stable subvariety $W\subset X\times X$ containing $x$, and a Cartesian diagram of pointed spaces with \'etale diagonal maps:
\begin{equation}\label{Luna}
\xymatrix{ 
 &  (W\times^{H_{\pm 1}}G,x) \ar[dl]\ar[dr]\ar[dd]
 \\
 (V_{\pm 1}\times^{H_{\pm 1}} G,0) \ar[dd] && (X\times X,x) \ar[dd] \\ 
& W\sslash H_{\pm 1}\ar[dl]\ar[dr]
 \\
 V_{\pm 1}\sslash H_{\pm 1} = \c_X^* && \C_X= X\times X\sslash G.}
\end{equation}
  
The \'etale diagonals induce isomorphisms between neighborhoods $U_1$ of $0\in \c_X^*(F)$ and $U_2$ of $[\pm 1]\in \C_X(F)$, and \cite[Corollary 4.2.1]{SaStacks} implies that such a diagram induces an isomorphism between the coinvariant spaces over these neighborhoods; more precisely (since we are not treating $V_{\pm 1}/H_{\pm 1}$ as a stack), between the $H_{\pm 1}(F)$-coinvariants of elements of $\mathcal S(V_{\pm 1})$ supported in the preimage of $U_1$, and the $G(F)$-coinvariants of elements of $\mathcal S(X\times X)$ supported in the intersection of the preimage of $U_2$ with the $G(F)$-orbit of the Luna slice $W(F)$. In particular, the pushforwards of those measures to $U_1\simeq U_2$ coincide. 
\end{proof}

From now on, we denote $V_{\pm 1}$ simply by $V$, and $H_{\pm 1}$ simply by $H$. The reader should not confuse that, in the case of $X_{-1}$, with the representation $X=H\backslash G$ used elsewhere in this paper. The dimension $d_{\pm 1}$ of $V$ will be denoted simply by $d$.

\subsection{Pullback to the blowup}

Let $\R_V:\tilde V\to V$ be the blowup of $V$ at the origin, and $E$ the preimage of $0$ (the exceptional divisor); it is the linear analog of the resolution $\R: \P J_X \to X\times X$.

\begin{lemma}
Pullback of Schwartz measures under the blowup $\R_V$ gives rise to a closed embedding:
 \begin{equation}\label{Schwartzembedding}
\mathcal S(V)\overset{\R_V^*}{\hookrightarrow}   \mathcal S(\tilde V, (d-1)[E]).
 \end{equation}
\end{lemma}

The space on the right is the space of Schwartz measures valued in the complex line bundle defined by the divisor $(d-1)[E]$, introduced in \S \ref{ssgeneralities-Schwartz}.
 
\begin{proof}
 This follows from writing any Schwartz measure, locally, as $f=\Phi\cdot |\omega|$, where $\Phi$ is a Schwartz function and $\omega$ a Haar volume form on $V$, and taking into account that the divisor of $\R_V^*\omega$ is $(d-1)[E]$.
\end{proof}

The blowup $\tilde V$ is canonically the total space of the tautological line bundle over $\P V=$ the exceptional divisor; let $\pi:\tilde V\to \P V$ be the projection to the zero section. 
Any element of $\mathcal S(\tilde V, (d-1)[E])$ can be written as a product $\Phi \R_V^* dv$, where $\Phi$ is a Schwartz function on $\tilde V$ and $dv$ is a Haar measure on $V$. 
The map $\Phi(v) \R_V^* dv \mapsto \Phi(\pi(v)) \R_V^* dv$
gives rise to a canonical identification
\begin{equation}
\label{stalks}
\overline{\mathcal S_E(\tilde V, (d-1)[E])} \simeq \Meas^\infty(V\smallsetminus\{0\})^{\Gm, |\bullet|^d}
\end{equation}
between the fiber of the Schwartz cosheaf $\mathcal S(\tilde V,(d-1)[E])$ over the exceptional divisor $E$, and the space of smooth measures on $V\smallsetminus\{0\}$ which are eigenmeasures for the multiplicative group of dilations with eigencharacter $|\bullet|^d$. \footnote{The action of $F^\times$ on measures is defined in duality with its (unnormalized) action on functions: $\left<a\cdot \mu, \Phi\right>=\left<\mu, a^{-1}\cdot \Phi\right>$, where $a^{-1}\cdot \Phi(x)= \Phi(a^{-1}x)$. In particular, Haar measure is $(\Gm,|\bullet|^d)$-equivariant.}

In particular, $\overline{\mathcal S_E(\tilde V, (d-1)[E])}$ contains a canonical line $\overline{\mathcal S_E(\tilde V, (d-1)[E])}_\Haar$, that will be called the line of ``Haar'' elements, corresponding to Haar measures on the right hand side of \eqref{stalks}. 

In the non-Archimedean case, the fiber $\overline{\mathcal S_E(\tilde V, (d-1)[E])}$ and the stalk $\mathcal S_E(\tilde V, (d-1)[E])$ coincide. In the Archimedean case, the image of a measure $\Phi(v) \R_V^* dv$ in the stalk at $E$ will be determined not only by the values, but also by the transversal derivatives of $\Phi$ along the exceptional divisor. Again, there is a distinguished $\mathcal C^\infty_0(V)$-submodule 
$$\mathcal S_E(\tilde V, (d-1)[E])_\Haar \subset \mathcal S_E(\tilde V, (d-1)[E])$$ 
(where $\mathcal C^\infty_0(V)$ denotes the stalk at $0\in V$ of the ring of smooth functions), generated by measures which are of the form $\R_V^* dv$ in a neighborhood of the exceptional divisor, and we have the following, almost tautological, lemma:

\begin{lemma}\label{lemmastalks}
In the non-Archimedean case, the image of $\mathcal S(V)$ under the embedding \eqref{Schwartzembedding} coincides with the space of those elements of $\mathcal S(\tilde V, (d-1)[E])$ whose image in the fiber over the exceptional divisor lies in $\overline{\mathcal S_E(\tilde V, (d-1)[E])}_\Haar$. 
 
In the Archimedean case, the image of $\mathcal S(V)$ coincides with the space of those elements of $\mathcal S(\tilde V, (d-1)[E])$ whose image in the stalk over the exceptional divisor lies in $\mathcal S_E(\tilde V, (d-1)[E])_\Haar$. 
\end{lemma}

\begin{proof}
Among measures of the form $\Phi(v) \R_V^* dv$, the image of $\mathcal S(V)$ consists of those where the function $\Phi$ is the pullback of a Schwartz function on $V$, i.e., precisely those which in a neighborhood of the exceptional divisor can be written as the product of an element of $\mathcal C^\infty_0(V)$ with the pullback $\R_V^* dv$ of a Haar measure. In the non-Archimedean case, the stalk and the fiber coincide, and the claim follows.
\end{proof}

\subsection{Pushforward to $\c_X^*$}

Now we consider pushforwards of Schwartz measures to $\c_X^*$.

\begin{proposition}\label{propbigimage}
 The image of the pushforward map 
 $$ \mathcal S(\tilde V, (d-1)[E]) \to \Meas(\c_X^*)$$ consists precisely of those measures which are smooth away from $0$, of rapid decay (together with their polynomial derivatives) at infinity (compactly supported, in the non-Archimedean case), and in a neighborhood of $0$ have the form:
 \begin{equation}\label{pushforwards-all} C_0(\xi) + |\xi|^{\frac{d}{2}-1} \sum_{\eta\in \widehat{F^\times/(F^\times)^2}} C_\eta(\xi) \cdot \eta(\xi),
 \end{equation}
 where $\eta$ runs over all quadratic characters of $F^\times$, $C_0$ and the $C_\eta$'s are smooth measures, and $\xi$ is a coordinate on $\c_X^*\simeq \mathbbm A^1$, except when $|\xi|^{\frac{d}{2}-1}\eta(\xi)$ is smooth for some $\eta$, that is:
 \begin{itemize}
  \item when $\frac{d}{2}-1 = 0$ and $\eta=1$, or 
  \item when $F=\mathbb R$ and $\frac{d}{2}-1$ is an even integer and $\eta$ is trivial, or an odd integer and $\eta$ is the sign character, or 
  \item when $F=\mathbb C$ (so, $\eta=1$) and $d$ is even,\footnote{We use the arithmetic normalization of absolute values, which is compatible with norms to the base field; this is the square of the usual absolute value in the complex case.} 
  \end{itemize}
 in which case  the term $|\xi|^{\frac{d}{2}-1} \eta(\xi)\cdot C_\eta(\xi)$ should be replaced by $|\xi|^{\frac{d}{2}-1} \eta(\xi) \log|\xi|\cdot C_\eta(\xi)$. In the Archimedean case, this map is continuous with respect to the obvious Fr\'echet topology on these measures, determined by Schwartz seminorms away from zero, and by absolute values of the derivatives of the functions $\frac{C_0}{d\xi}$, $\frac{C_\eta}{d\xi}$ at zero. 
\end{proposition}

\begin{proof}
I claim that, locally around any point of the exceptional divisor, there is a coordinate chart $(\epsilon_E, x_1, \dots, x_{d-1})$, where $\epsilon_E=0$ is a local equation for $E$, such that the map $\xi:\tilde V\to \c_X^*\simeq \mathbbm A^1$ is given by $\xi=\epsilon_E^2 x_1$. Indeed, this is seen immediately by writing the split quadratic form $V\to \c_X^*$ in coordinates: $\xi=y_1^2 - y_2^2 +\dots \pm y_d^2$, and setting, e.g., $x_1 = \frac{\xi}{y_d^2}$, $x_i = \frac{y_i}{y_d}$, for $2\le i\le d-1$, $\epsilon_E = y_d$.\footnote{In terms of the map $\P J_X\to \C_X$, of which the map $\tilde V\to \c_X^*$ is the ``linearization'', and given that $J\bullet N^*_{X_{\pm 1},\ne 0}$ is smooth over $J$ by non-degeneracy of the quadratic forms, the first of the maps $\P J_X\to \P J \xrightarrow{Q} \C_X$ is smooth, and the second is given, in  coordinates $t_0^2-\xi t_1^2=1$ for $J$, and a suitable identification $A_X\sslash W_X = \mathbbm A^1$, by $Q=(t_0 \pm 1)^{-1} \cdot \xi t_1^2$. In a neighborhood of $t_0=\pm 1$, where  the function $(t_0 \pm 1)^{-1}$ is a nonvanishing, smooth semialgebraic function of $(\xi, t_1)$, we can set $\xi'=(t_0 \pm 1)^{-1}\xi$, and we get that the map is given by $Q=\xi' t_1^2$.}

The asserted form of the pushforward of Schwartz measures under such a map is quite a standard result. One way to prove it is using Mellin transforms: The Mellin transform of a pushforward measure $\xi_* f$ with respect to the variable $\xi$, with $f = |\epsilon_E|^{d-1}\cdot C(\epsilon_E, x_1, \dots, x_{d-1})$, where $C$ is a Schwartz measure in $d$ variables, is 
$$ \widecheck{\xi_* f}(\chi) := \int_F \xi_* f(\xi) \chi^{-1}(\xi) = \int \bar C(\epsilon_E, x_1) |\epsilon_E|^{d-1} \chi^{-1}(\epsilon_E^2 x_1),$$
where $\bar C$ is the pushforward of $C$ with respect to the map $(\epsilon_E, x_1, \dots, x_{d-1})\mapsto (\epsilon_E, x_1)$.

This is the Tate zeta integral  of a Schwartz measure in two variables, in one of the variables against the character $\chi^{-1}$ and in the other against the character $|\bullet|^{d-1} \chi^{-2}$. In the non-Archimedean case, it has poles at $\chi=|\bullet|^{-1}$ and at the points $\chi=\eta \cdot |\bullet|^{\frac{d}{2}-1}$ (double if any of these points coincide, simple otherwise), where $\eta$ ranges over all quadratic characters. Such a Mellin transform corresponds to a measure on the line which in a neighborhood of $\xi=0$ is of the form 
$$C_0(\xi) + |\xi|^{\frac{d}{2}-1} \sum_{\eta\in \widehat{F^\times/(F^\times)^2}} C_\eta(\xi) \cdot \eta(\xi),$$ 
unless $d=2$, in which case the pole at $\chi=|\bullet|^{-1}$ is double, and the corresponding singular term is of the form $C_1(\xi) \cdot \log|\xi|$.

A similar argument works in the Archimedean case, where double poles appear whenever the product $|\xi|^{\frac{d}{2}-1} \eta(\xi)$ is a smooth function of $\xi$. Here, the above Tate integral maps continuously into the appropriate ``Paley--Wiener space'' in the language of \cite[Remark 2.1.6]{SaTransfer1}, with the location and multiplicity of poles determined by the characters $|\xi|^{\frac{d}{2}-1} \eta(\xi)$, which corresponds to the Fr\'echet space of measures as in the statement of the proposition.
\end{proof}

Our final task will be to determine the image $\mathcal S(V/H)$ of the subspace $\mathcal S(V)\hookrightarrow \mathcal S(\tilde V, (d-1)[E])$. This will be completed in the next section. We start with the following observation:

\begin{theorem}\label{pushfthm}
 The space $\mathcal S(V/H)$ contains the space $\mathcal S(\c_X^*)$ of Schwartz measures on $\c_X^*$.
 
 Moreover, in the expression \eqref{pushforwards-all} for the pushforward of a measure $f\in \mathcal S(\tilde V, (d-1)[E])$, the coefficients $\frac{C_\eta}{d\xi}(0)$ depend only on the image of $f$ in the fiber 
$$\overline{\mathcal S(\tilde V, (d-1)[E])}=\Meas^\infty(V\smallsetminus\{0\})^{\Gm, |\bullet|^d}$$
(see \eqref{stalks}). In particular, by Lemma \ref{lemmastalks}, for all $f\in \mathcal S(V)$ these coefficients will lie in a one-dimensional subspace of $\CC^{\widehat{F^\times/(F^\times)^2}}$.
 
 Let $(a_\eta)_\eta$ be a vector spanning this one-dimensional subspace. Then $\mathcal S(V/H)$ is the space of those measures of the form 
 \begin{equation}\label{pushforwards} C_0(\xi) + |\xi|^{\frac{d}{2}-1} C_\sing(\xi) \sum_{\eta\in \widehat{F^\times/(F^\times)^2}} a_\eta \cdot \eta(\xi),
 \end{equation}
 where $C_0, C_\sing$ are Schwartz measures, and the same modification as in Proposition \ref{propbigimage} applies to the case where $|\xi|^{\frac{d}{2}-1}\eta(\xi)$ is smooth.
\end{theorem}

\begin{proof}
 As we have seen, the complement of the origin is smooth and surjective over $\c_X^*$, hence the image of $\mathcal S(V\smallsetminus\{0\})$ is equal to $\mathcal S(\c_X^*)$.
 
 Hence, the germs of the measures $C_\eta$ at $0$ depend only on the image of $f$ in the stalk $\mathcal S_E(\tilde V, (d-1)[E])$. 
 
 In the non-Archimedean case, ``germ of $C_\eta$'' means simply the value of $\frac{C_\eta}{d\xi}$ at $0$, and, as we saw in \eqref{stalks}, the stalk $\mathcal S_E(\tilde V, (d-1)[E])$  is identified with $\Meas^\infty(V\smallsetminus\{0\})^{\Gm, |\bullet|^d}$, so the claim follows.
 
 In the Archimedean case, we will show that the $H$-coinvariants of the stalk $\mathcal S_0(V)$ are generated over the stalk $C_0^\infty(\c_X^*)$ by any measure which is nonvanishing at the origin.\footnote{We define coinvariants of Fr\'echet spaces by dividing by the \emph{closure} of the space generated by elements of the form $f-h\cdot f$.}
 
 For this, consider the descending filtration of the stalk $\mathcal S_0(V)$ which defines its topology, i.e., $F^n\mathcal S_0(V)=$ the germs of smooth measures $f=\Phi dv$ (where $dv$ is a Haar measure) such that all partial derivatives of $\Phi$ of order $<n$ vanish at the origin. By the $\Gm$-action on $V$, this filtration corresponds to a grading on the dense subspace of $\Gm$-finite germs.  Notice that $X\times X$ admits a $G$-invariant measure, and therefore the Haar measure $dv$ is $H$-invariant; thus, we can choose such a measure to identify the $H$-modules of functions and measures. The graded piece $F^n\mathcal S_0(V)/F^{n+1}\mathcal S_0(V)$ is then identified with $\Sym_\RR^n(V^*)\otimes_\RR \CC$, and therefore the $H$-coinvariants of the stalk are
  \begin{eqnarray}\label{sympowers}\nonumber \mathcal S_0(V)_H = \lim_{\underset{n}\leftarrow} \left(\Sym_\RR^n(V^*)\otimes_\RR \CC\right)^{H} &=& \lim_{\underset{n}\leftarrow} \left(\Sym_\CC^n(V^*\otimes_\RR \CC)\right)^{H_\CC} = \\ &=&\begin{cases} \CC[[\xi]], & \mbox{ if } F=\RR,\\ \CC[[\xi, \bar\xi]],& \mbox{ if } F=\CC,\end{cases}\end{eqnarray}
 where we have treated $H$ as a real group, so that $H_\CC$ denotes its complexification.

 The space $\CC[[\xi]]$, if $F=\RR$, and $\CC[[\xi, \bar\xi]]$, if $F=\CC$, is naturally identified with the stalk $C^\infty_0(\c_X^*)$ at zero of the ring of smooth functions. This stalk is a ring that acts on $\mathcal S_0(V)_H$, and the isomorphism \eqref{sympowers} is equivariant with respect to the action of this stalk. Thus, the above calculation shows that the $H$-coinvariants of $\mathcal S_0(V)$ are freely generated over $C^\infty_0(\c_X^*)$ by the germ of any element $\Phi dv$ with $\Phi(0)\ne 0$. 
 
 Thus, the germs of pushforwards will also be generated, over $C^\infty_0(\c_X^*)$ and up to smooth measures, by the germ of the pushforward of any such measure $\Phi dv$. Consider such a measure with $\Phi(v)$ constant (and $\ne 0$) close to the origin. The pushforward map is $\Gm$-equivariant (with respect to the quadratic action on $\c_X^*$), hence, in terms of the expression \eqref{pushforwards-all}, the germ of the pushforward of such a measure $\Phi dv$ will be of the form  
 $$ C_0(\xi) + |\xi|^{\frac{d}{2}-1} \sum_{\eta\in \widehat{F^\times/(F^\times)^2}} C_\eta(\xi) \cdot \eta(\xi),$$
 where the measures $C_\eta$ are \emph{constant} around $\xi=0$. Thus, the image of the Haar \emph{stalk} $\mathcal S_E(\tilde V, (d-1)[E])_\Haar$ (in the notation of Lemma \ref{lemmastalks}) in the singular quotient of the stalk of $\mathcal S(V/H)$ at zero (i.e., ignoring the term $C_0$) is of the form 
 \begin{equation}\label{Haarstalk} |\xi|^{\frac{d}{2}-1} C_\sing(\xi) \sum_{\eta\in \widehat{F^\times/(F^\times)^2}} a_\eta \cdot \eta(\xi),\end{equation}
 and in particular is completely determined by the coefficients $\alpha_\eta$, which depend only on the image of an element in the \emph{fiber} $\overline{\mathcal S_E(\tilde V, (d-1)[E])}_\Haar$.
\end{proof}

Therefore, we are left with computing the ratio between the coefficients $a_\eta$, which correspond to the singular part of the pushforward of a measure on $V$ which restricts to a Haar measure in a neighborhood of the origin.

\section{Determination of the germs}\label{sec:germs}

\subsection{Reduction to the basic cases}

We will actually not compute the ratio of the coefficients $a_\eta$ explicitly in all cases, but rather prove, by reducing to an $\SL_2$- or $\PGL_2$-example, that they match the contributions of the Kloosterman germs under the transfer operator from the Kuznetsov formula. 
 The cases $d=$ even and $d=$ odd will be quite different, as we will see. We fix throughout the isomorphism $\xi:\c_X^*\xrightarrow\sim \mathbbm A^1$ with $0\in \c_X^*$ mapping to $0$, thus viewing the map $V\to V\sslash H=\mathfrak c_X^*$ as a quadratic form.

The main result of this subsection is Proposition \ref{reductiontobasic}, which says that pushforwards of Schwartz measures for a $d$-dimensional split quadratic space (under the quadratic map) are equal to twisted pushforwards for a two- or three-dimensional quadratic space (of the same parity as $d$); this will complete the proof of Theorem \ref{Xtheorem}. 

The two- or three-dimensional quadratic space $V_2$ is obtained from $V$ by choosing a maximal isotropic subspace $M\subset V$ and a hyperplane $M'\subset M$; then $V_2 = M'^\perp/M'$. Let us go through the argument carefully:

Fix such a maximal isotropic subspace $M$. Since $V$ is split, the orthogonal complement $M^\perp$ is either equal to $M$ (when $d$ is even), or contains $M$ as a hyperplane (when $d$ is odd). The quotient $V/M^\perp$ is isomorphic to the linear dual $M^*$ through the quadratic form, and the parabolic $P\subset \SO(V)$ stabilizing $M$ surjects to $\GL(M^*)$. 

The integration (pushforward) map $\mathcal S(V)\to \CC$ factors through surjective pushforward maps:
\begin{equation}\label{sequence} \mathcal S(V)\to \mathcal S(V/M^\perp) \to \mathcal S(\P M^*) \to \CC.\end{equation}

Let $M'\in \P M^*$, identified (and denoted by the same letter) with a hyperplane in $M$. Let $P'\subset P$ be the stabilizer of the flag $M'\subset M\subset V$. 
The space $\mathcal S(\P M^*)$, considered as a representation of $P$, can be identified, up to a scalar which we fix, with the (unnormalized) induced representation
$\Ind_{P'}^P (\delta_{P'/U_P})$, where $\delta_{P'/U_P}$ denotes the modular character of the image of $P'$ in the Levi quotient of $P$. By Frobenius reciprocity, the $P$-equivariant map $\mathcal S(V)\to \mathcal S(\P M^*)$ is given by a $(P', \delta_{P'/U_P})$-equivariant functional. The Lemma that follows determines this functional:

\begin{lemma}
 Let $\Phi$ be a Schwartz function on $V$, and $dv$ a Haar measure. Then, for suitable Haar measures,
 \begin{equation}\label{PV} \int_V \Phi(v) dv = \int_{\P M^*} \left(\int_{\Gm} \int_{M^\perp} \Phi(av + v_1) dv_1 |a|^{\dim M} d^\times a \right) dv.
 \end{equation}
\end{lemma}

Notice that the expression in brackets, viewed as a function of $v\in M^*\smallsetminus\{0\}$, is $(\Gm, |\bullet|^{-d})$-equivariant, hence $dv$ denotes an invariant measure on $\P M^*$, valued in the dual of the line bundle of $(\Gm, |\bullet|^{-d})$-equivariant functions on $M^*\smallsetminus\{0\}$. More precisely, under the action of $P'$, the expression in brackets is $\delta_{P'/U_P}$-equivariant, and $dv$ 
is an invariant measure on $\P M^*$, valued in the line bundle dual to the one induced from this character of $P'$. 

\begin{proof} This lemma is just a reformulation of the sequence \eqref{sequence}.\end{proof}

Let us reformulate the inner integral of \eqref{PV}: Fix $M'\in \P M^*$,  understood again as a hyperplane in $M$. Its preimage in $V$ under the rational map $V\to \P M^*$ is equal to $M'^\perp \smallsetminus M^\perp$. Fix a nonzero vector $v\in M^*$ in the line corresponding to $M'$; then the functional $L: av\mapsto a$ is a linear functional on the one-dimensional space of multiples of $v$ in $M^*$ or, equivalently, a functional 
$$L: M'^\perp \to \Ga.$$ 

The quotient $V_2:=M'^\perp/M'$ is a nondegenerate quadratic space of dimension $2$ or $3$ (same parity as $V$). Fix a Haar measure $dv'$ on $M'$, and let $\Phi\mapsto \Phi_2$ be the corresponding pushforward map (integration over cosets of $M'$ against $dv'$)
$$ \mathcal F (M'^\perp) \twoheadrightarrow \mathcal F(V_2),$$
where $\mathcal F$ denotes the spaces of Schwartz functions. Then the inner integral of \eqref{PV} can be written as
\begin{equation} \int_{V_2} \Phi_2(v_2) L(v_2)^{\dim M-1} dv_2,\end{equation}
for a suitable Haar measure $dv_2$.

Let us explicate this integral:

\begin{itemize}
 \item If $V_2$ is a two-dimensional (split) quadratic space, then $\dim M=\frac{d}{2}$, and there are coordinates $(x,y)$ such that the quadratic form is $\xi=xy$ and the functional $L$ is $L = x$, so the integral reads:
 \begin{equation}\label{twodim} \int_{V_2} \Phi_2(x,y) |x|^\frac{d-2}{2} dx dy.\end{equation}

 \item If $V_2$ is a three-dimensional (split) quadratic space, then $\dim M=\frac{d-1}{2}$, and there is an isomorphism $V_2\simeq\mathfrak{sl}_2$ with quadratic form $\xi = -\det$ and  $L\begin{pmatrix} A & B \\ C & -A \end{pmatrix} = C$, so the integral reads:
 \begin{equation}\label{threedim} \int_{V_2} \Phi_2\begin{pmatrix} A & B \\ C & -A \end{pmatrix} |C|^{\frac{d-3}{2}} dA dB dC.\end{equation}
\end{itemize}

In either case, these integrals can be disintegrated against the quadratic form, i.e., written as iterated integrals $\int_{\c_X^*} \int$ with the interior integral taken over the fibers of the map $V_2\to \c_X^*$, but we need to choose a section $\sigma: \c_X^* \to V_2$, since the integrand is not invariant over the fibers. Choose this section $\sigma$ so that its image is contained in an affine line of the hyperplane $L=1$; then it is necessarily contained in the affine line $\sigma(0)+\bar M$ (where $\bar M=M/M'$, the image of $M$ in $V_2$); explicitly:
\begin{itemize}
 \item $\sigma(\xi) = (1,\xi)$ in the coordinates above when $V_2$ is two-dimensional;
 \item $\sigma(\xi) = \begin{pmatrix}  & \xi \\ 1 &  \end{pmatrix}$ when $V_2$ is three-dimensional.
\end{itemize}

Let $B_2\subset \SO(V_2)$ be the stabilizer of the isotropic line $\bar M=M/M'$, and let $\delta_2^{-1}$ be the absolute value of the character by which it acts on $\bar M$. Then:

\begin{lemma}\label{rewrite}
 The expressions \eqref{twodim} and \eqref{threedim} can be written:
 \begin{equation}\label{onV2} \int_{\c_X^*} \int_{B_2} \Phi_2(\sigma(\xi) b) \delta_2(b)^{\dim M-1} db d\xi,
 \end{equation}
 for a suitable right Haar measure on $B_2$.
\end{lemma}

\begin{proof}
 In the coordinates above, in the two-dimensional case, $\bar M$ is the line $(0,*)$ and $B_2 = \SO(V_2) = \{(a,a^{-1})\}$, acting by the character $a^{-1}$ on $\bar M$. The integral \eqref{twodim} can be written 
 $$ \int \Phi_2(a, a^{-1} \xi ) |a|^{\dim M-1} d^\times a d\xi.$$
 
 In the three-dimensional case, $\bar M$ corresponds to the subspace $A=C=0$ in \eqref{threedim}, and $B_2$ is the upper-triangular Borel subgroup of $\SO(V_2) = \PGL_2$, acting on this line via the inverse of its modular character (remember that it acts on the right), and we can write \eqref{threedim} as 
 $$ \int_{V_2} \Phi_2 \left(\begin{pmatrix} a^{-1} & -a^{-1} b \\  &  1\end{pmatrix} \begin{pmatrix}  & \xi \\ 1 &  \end{pmatrix} \begin{pmatrix} a & b \\  & 1 \end{pmatrix} \right) |a|^{\dim M -1} d^\times a db d\xi.$$
\end{proof}

By this lemma, we have a new integration formula for $V$ in terms of the quadratic form, that includes $B_2$-orbital integrals on $V_2$, twisted by the character $\delta_2(b)^{\dim M-1}$. Thus, the integration formula for a $d$-dimensional quadratic space involves a twisted integration formula for a $2$- or $3$-dimensional quadratic space:

\begin{corollary}
  Let $\Phi$ be a Schwartz function on $V$, and $dv$ a Haar measure. Let $K\subset P$ be any compact subgroup such that $K\to P'\backslash P=\P M^*$ is surjective. Then, for suitable (right) Haar measures,
 \begin{equation}\label{PV2} \int_V \Phi(v) dv = \int_{\c_X^*} \left(\int_{K} \int_{B_2}  \int_{M'} \Phi((\sigma(\xi) b + v')k )  \delta_2(b)^{\dim M-1} dv' db dk\right) d\xi.
 \end{equation}
\end{corollary}

Here, we have lifted the section $\sigma$ to $M'^\perp$ and the group $B_2$ to $P'$, by choosing a section $V_2\to M'^\perp$.
 
\begin{proof}
 The integral over $K$ replaces the integral over $\P M^*$ in \eqref{PV}; since it is a compact integral of a smooth function, it can be moved to the interior, and the result follows by applying \eqref{onV2}.
\end{proof}
 
\begin{remark}
 Representation-theoretically, the two inner integrals of \eqref{PV2} represent a $(P', \delta_{P'/U_P})$-equivariant functional, hence a morphism 
 $$ \mathcal S(V)\to \Ind_{P'}^{\SO(V)}(\delta_{P'/U_P}) = \Ind_P^{\SO(V)} \Ind_{P'}^P(\delta_{P'/U_P}).$$
 The integral over $K$ corresponds to the quotient $\Ind_{P'}^P(\delta_{P'/U_P})\to \CC$ (the trivial representation), so the expression in brackets can be seen as a morphism 
 $$ \mathcal S(V)\to \Ind_P^{\SO(V)} (\CC).$$
 By the invariance of the left hand side of \eqref{PV2}, this morphism is $\SO(V)$-invariant, hence has image in the trivial subrepresentation of $\Ind_P^{\SO(V)}(\CC)$.
 
 Let $\ell_\xi(\Phi)$ be the $\SO(V)$-invariant functional represented by the expression in brackets of \eqref{PV2}.
Comparing \eqref{PV2} with the integration formula of Theorem \ref{thmintegrationformula-linear} (for the special case $H=\SO(V)$), we get:
 \begin{equation}\ell_\xi(\Phi) = |\xi|^{\frac{d}{2}-1} O_\xi(\Phi),
 \end{equation}
 where the $O_\xi$'s are orbital integrals on the fibers over $\xi\ne 0$, against invariant measures $d\dot{g}_\xi$ obtained by identifying all nondegenerate $\SO(V)$-orbits over the algebraic closure, and choosing volume forms as in Theorem \ref{thmintegrationformula-linear}.
\end{remark}
 
Thus, we arrive at the following result about the coefficients $a_\eta$ of the expression \eqref{pushforwards}:

\begin{proposition}\label{reductiontobasic}
 If $d$ is even, we have $a_\eta=0$ except for $\eta=1$, and there is an equality between the space of pushforward measures for $V\xrightarrow{\xi} \mathbbm A^1$ and the measures on $\mathbbm A^1$ of the form 
 \begin{equation}\label{twistedpfT} \xi\mapsto  \left(\int_{\Gm} \Phi_2(a, a^{-1}\xi)|a|^\frac{d-2}{2} d^\times a\right)   d\xi,
 \end{equation}
 where $\Phi_2$ varies among Schwartz functions on $\mathbbm A^2$.
  
 If $d$ is odd, there is an equality between the space of pushforward measures for $V\xrightarrow{\xi} \mathbbm A^1$ and the measures on $\mathbbm A^1$ of the form
 \begin{equation}\label{twistedpfG} \xi\mapsto \left(\int_{B_\ad} \Phi_2\left(\Ad(b^{-1}) \begin{pmatrix}
                                              & \xi \\ 1 
                                             \end{pmatrix}
 \right)\delta_2(b)^\frac{d-3}{2} db\right) d\xi,
 \end{equation}
 where $\Phi_2$ varies among Schwartz functions on $\mathfrak{sl}_2$, $B_\ad$ denotes the upper triangular Borel subgroup of $\PGL_2$, $\delta_2$ is its modular character, and $db$ is a right Haar measure.
\end{proposition}

Here, sticking with standard notation, we have denoted by $\Ad$ the \emph{left} adjoint representation of $\PGL_2$ on $\mathfrak{sl}_2$; but recall that our convention is that $G$ acts on the right on $X\times X$, hence $H$ acts on the right on $V$, and this convention is extended to the group $\SO(V)$.

\begin{remarks}
\begin{enumerate}
 \item  In other words, the germs are reduced to \emph{twisted versions} of the infinitesimal versions of the \emph{basic cases} $X=\Gm\backslash \PGL_2$ and $X=\SL_2=\SO_3\backslash \SO_4$. Indeed, the linearizations of those two are, respectively, $\mathbbm A_2/\Gm$ and $\mathfrak{sl}_2/\PGL_2$, and the latter can also be replaced by $\mathfrak{sl}_2/B_\ad$, because the affine quotients $\mathfrak{sl}_2\sslash \PGL_2$ and $\mathfrak{sl}_2\sslash B_\ad$ are the same. Putting an appropriate character on $\Gm$ or $B_\ad$, we obtain the germs for the general case. This fact will be used to relate those germs to the Kloosterman germs of the Kuznetsov formula, under the transfer operator.
 \item As we saw in \eqref{twodim}, \eqref{threedim}, the measures \eqref{twistedpfT}, \eqref{twistedpfG} can be considered as twisted pushforwards of the Haar measures $\Phi_2 dv$ (where $dv$ is a Haar measure on $\mathbbm A^2$, resp.\ $\mathfrak{sl}_2$), dual to the twisted pullback maps:
 \begin{equation}\label{twistedpbT}
  \Psi \mapsto \tilde\Psi(x, y) = \Psi(xy) |x|^\frac{d-2}{2},
 \end{equation}
 resp.
 \begin{equation}\label{twistedpbG}
  \Psi \mapsto \tilde\Psi\begin{pmatrix} A & B \\ C & -A \end{pmatrix} = \Psi(A^2+BC) |C|^{\frac{d-3}{2}}.
 \end{equation}
\end{enumerate}
\end{remarks}

\begin{proof}
 Let $\ell_\xi(\Phi)$ be the $\SO(V)$-invariant functional represented by the expression in brackets of \eqref{PV2}. By that formula, the pushforward $f$ of $\Phi dv$ can be written 
 $$ f(\xi) = \ell_\xi(\Phi) d\xi = \int_{K} \int_{B_2}  \Phi_2^K(\sigma(\xi) b)  \delta_2(b)^{\dim M-1} db dk,$$
 where $\Phi_2^K \in  \mathcal F(V_2)$ is the Schwartz function 
 \begin{equation}\label{PhiK} v_2\to \int_{K} \int_{M'} \Phi((v_2+v')k)  dv' dk.
 \end{equation}

 As we have seen in Lemma \ref{rewrite}, this is equal to the expressions \eqref{twistedpfT}, \eqref{twistedpfG} in the two cases, applied to the function $\Phi_2^K$. We just need to argue that these spaces of pushforward measures obtained from a function of the form $\Phi_2^K$ is the same as the space obtained from an arbitrary Schwartz function $\Phi_2\in \mathcal F(V_2)$.  The pushforward map $\mathcal S(M'^\perp)\to \mathcal S(V_2)$ is surjective, and 
 starting from an arbitrary Schwartz measure $\Phi_2 dv_2\in \mathcal S(V_2)$ we can choose a preimage $\Phi_1 dv_1\in \mathcal S(M'^\perp)$. 
Without loss of generality (in terms of the output of \eqref{twistedpfT}, \eqref{twistedpfG}), we will assume that $\Phi_1$ is $K\cap P'$-invariant.

 Notice that we have freedom in choosing $K$, as long as the map $K\to P'\backslash P=\P M^*$ is surjective. Identify $M^*$ as a subspace of $V$ through an isotropic splitting of the quotient $V\to M^*$, so that we have a direct sum decomposition $V=M^\perp\oplus M^*$, and choose $K$ in the Levi subgroup $\GL(M^*)\subset P$. Any element of $K$ fixing the line in $M^*$ corresponding to $M'$ has to belong to $K\cap P'$. Thus, for any two $v_1, v_2\in M'^\perp$ with nonzero image in $M^\perp$, the relation $v_1\cdot k=v_2$ for some $k\in K$ implies that $v_2\in v_1\cdot (K\cap P')$. Hence, the map of topological quotients
 $$ M'^\perp/K\cap P' \to V/K,$$
 surjective by our assumption on $K$, is also injective. Thus, $\Phi_1$, a $K\cap P'$-invariant function on $M'^\perp$, is the restriction of a unique $K$-invariant function $\Phi$ on $V$; in particular, the average $\Phi_2^K$ (defined in terms of $\Phi$, and using probability measure on $K$) is equal to the function $\Phi_2$ that we started from. 
 
 In the non-Archimedean case, it is immediate to see that if $\Phi_1$ is smooth, so is $\Phi$, because smoothness means that they are locally constant on the topological quotient $V/K$. In the Archimedean case, taking $K=\SO_n(\RR)$ when $F=\RR$ and $K=\operatorname{U}_n(\RR)$ when $F=\CC$ (where $n=\dim M$), the quotient $M^*/K$ can be identified with $\RR_{\ge 0}$ through the distance function from the origin. For any one-dimensional subspace $Fv$ of $M^*$, any $K\cap \GL_1(F)$-invariant smooth function on $Fv$ is the restriction of a smooth radial function on $M^*$. 
 
 Finally, in the case $d=$ even, to show that the coefficients $a_\eta$ with $\eta\ne 0$ vanish, we can use the same argument as in the proof of Proposition \ref{propbigimage} in order to analyze the twisted orbital integrals \eqref{twistedpfT}. In this case, without any need to pass to a blowup, 
 the Mellin transform in the variable $\xi$ can be written as a product of two Tate integrals:
$$\int \chi^{-1}(\xi)
 \left(\int \Phi_2(a, a^{-1}\xi)|a|^\frac{d-2}{2} d^\times a\right)   d\xi = \iint \Phi_2(x, y) \chi^{-1}(xy) |x|^\frac{d}{2} |y| d^\times x d^\times y,$$
 which is a holomorphic multiple of the product of local Dirichlet $L$-functions
 $$ L(\chi^{-1}, \frac{d}{2}) L(\chi^{-1},1).$$
 Hence as in the proof of  Proposition \ref{propbigimage}, the measure \eqref{twistedpfT} is a linear combination of the form 
 $$C_0(\xi) + |\xi|^{\frac{d}{2}-1} C_1(\xi),$$ 
 with $C_0$ and $C_1$ Schwartz measures, with the usual logarithmic modification when $|\xi|^{\frac{d}{2}-1}$ is smooth at zero.
\end{proof}

By Proposition \ref{linearization},  this completes the proof of Theorem \ref{Xtheorem}, which we state here more precisely. Notice that the precise local ($F$-analytic) isomorphism between a neighborhood of $0\in \c_X^*$ and a neighborhood of $[\pm 1]\in \C_X$, mentioned in Proposition \ref{linearization}, is not important, since the germs are invariant under any $F$-analytic automorphism. In particular, fixing isomorphisms $\c_X^*\simeq\mathbbm A^1_\xi$ and $\C_X=\mathbbm A^1_c$ with $c_{\pm 1}$ corresponding to $[\pm 1]$, we can take $\xi = c-c_{\pm 1}$. 

\begin{theorem}\label{Xtheorem2}
 There is a canonical isomorphism $\C_X:= X\times X\sslash G\simeq A_X\sslash W_X$, and the map $X\times X\to \C_X$ is smooth away from the preimages of $[\pm 1]$, where $[\pm 1]$ denote the images of $\pm 1\in A_X$ in $A_X\sslash W_X$. 
 
 In particular, there are two distinguished closed $G$-orbits $X_1 = X^\diag$ and $X_{-1}$ (over $[\pm 1]$, respectively); if $d_{\pm 1}$ denote their codimensions, then $d_1=\dim X$ and 
 $$d_{-1} = \epsilon \left<2\rho_{P(X)},\check\gamma\right>-d_1+2,$$ 
 where $\check\gamma$ is the spherical coroot, $2\rho_{P(X)}$ is the sum of roots in the unipotent radical of $P(X)$, and
$$ \epsilon = \begin{cases} 1, \mbox{ when the spherical root is of type $T$ (dual group $\SL_2$)};\\
2, \mbox{ when the spherical root is of type $G$ (dual group $\PGL_2$).}
\end{cases}$$
In the case of root of type $G$, $d_1=d_{-1}$.

 The space $\mathcal S(X\times X/G)$ consists of measures on $\C_X(\simeq \mathbbm A^1_c)$ which are smooth and of rapid decay, together with their polynomial derivatives (compactly supported in the non-Archimedean case) away from neighborhoods of $c_{\pm 1}$ (the coordinates of the points $[\pm 1]$), while in neighborhoods of $c_{\pm 1}$ they are of the form \eqref{twistedpfT} --- when the spherical root is of type $T$ --- or \eqref{twistedpfG} --- when the spherical root is of type $G$---, with $\xi = c-c_{\pm 1}$ and $d=d_{\pm 1}$.
\end{theorem}

\begin{remark}
 Recall that, up to a linear combination of coefficients, the singularities of the measures of the form \eqref{twistedpfT}, \eqref{twistedpfG} have been described explicitly in Theorem \ref{pushfthm}. In the Archimedean case, there is a natural Fr\'echet topology on the space of these measures, and by Proposition \ref{propbigimage} the map from $\mathcal S(X\times X)$ is continuous; hence, the quotient topology on $\mathcal S(X\times X/G)$ coincides with the natural Fr\'echet topology on the space of such measures.
\end{remark}

\subsection{Completion of the proof of the main theorem} \label{sspfmainthm}

We are now ready to prove Theorem \ref{maintheorem}. The remaining step of the proof relies on a result from \cite{SaTransfer2}, which I state here:

\begin{proposition}\label{fromtransfer}
 Let $G^*=\SL_2$. The operator
 $$f \mapsto  |\zeta|^{-\frac{1}{2}-s} \cdot \left( |\bullet|^{\frac{3}{2}+s} \psi(\bullet) d^\times\bullet\right) \star  f(\zeta)$$ 
 defines an isomorphism 
 \begin{equation}\label{SL2isom} 
 \mathcal S^-_{L(\Ad,\frac{1}{2}-s)} (N_\psi\backslash G^*/N_\psi) \xrightarrow\sim \mathcal S\left(\frac{\SL_2}{B_\ad,\delta^{\frac{1}{2}+s}}\right)
 \end{equation}
 away from the poles of the local $L$-functions $L(\eta, \frac{1}{2} -s )$, where $\eta$ ranges over all quadratic characters; here, the space on the right is the space of measures on the affine line $\mathbbm A^1$ (with coordinate $\zeta$), which are smooth away from $\zeta=\pm 2$, of rapid decay (together with their derivatives) at infinity, and in a neighborhood of $\zeta=\pm 2$, setting $\xi=\zeta\mp 2$, are of the form \eqref{twistedpfG}, with $\frac{d-3}{2}=-\frac{1}{2}-s$.
\end{proposition}

\begin{proof}
 This is \cite[Proposition 8.3.3]{SaTransfer2}, in the special case $\chi = \delta^s$.
\end{proof}

We will need an analogous result for the transform of the Kloosterman germ when $G^*=\PGL_2$. A special case of this result was proven in \cite{SaBE1}; therefore, I will confine myself to sketching the proof of the general case:

\begin{proposition}\label{frombeyond}
 Let $G^*=\PGL_2$. The operator
 $$f \mapsto  \left( |\bullet|^{\frac{1}{2}-s_1} \psi(\bullet) d\bullet\right) \star  \left( |\bullet|^{\frac{1}{2}-s_2} \psi(\bullet) d\bullet\right) \star f$$ 
 takes elements of $\mathcal S(N_\psi\backslash G^*/N_\psi)$ into a space $\mathcal M$ of measures on the affine line which in the neighborhood of $\xi = 1$ are of the form 
 \begin{equation} \label{at1}
  (\Phi_1(\xi) +\Phi_2(\xi) |\xi-1|^{s_1+s_2-1}) d\xi,
 \end{equation}
 where $\Phi_1, \Phi_2$ are smooth functions, with the second summand replaced, when $|\xi-1|^{s_1+s_2-1}$ is smooth, by $\Phi_2(\xi)|\xi-1|^{s_1+s_2-1} \log|\xi-1|$.

 Moreover, the subspace $\mathcal S(F^\times)\subset \mathcal S(N_\psi\backslash G^*/N_\psi)$ gets mapped to smooth measures in a neighborhood of $\xi=1$, and the transform descends to an isomorphism
 $$ \mathcal S(N_\psi\backslash G^*/N_\psi)/\mathcal S(F^\times) \xrightarrow\sim \mathcal M_1/\mathcal S(F)_1,$$
 between the ``Kloosterman stalk'' and the ``singular stalk of $\mathcal M$ at $\xi=1$'', that is, the quotient of $\mathcal M_1 =$ the stalk of $\mathcal M$ at $\xi=1$ by $\mathcal S(F)_1=$ the stalk of smooth measures.
\end{proposition}

More precise information will be obtained about the image of the space $\mathcal S(N_\psi\backslash G^*/N_\psi)$ and its extension $\mathcal S^-_{L_X}(N_\psi\backslash G^*/N_\psi)$ under the transfer operator in Theorem \ref{maintheorem2}.

\begin{proof}
 When $s_1=s_2=\frac{1}{2}$, this is a special case of \cite[Theorem 5.1]{SaBE1}, more precisely, the matching of the short exact sequences (5.3) and (5.4). 
 The same arguments work in the general case; they rely on the fact that the Kloosterman germ can be explicitly written as a Fourier transform in the case of $G^*=\PGL_2$, see \cite[Proposition 4.8]{SaBE1}. Thus, I will leave the verification to the reader, mentioning only that, after application of the first convolution, by $\left( |\bullet|^{\frac{1}{2}-s_2} \psi(\bullet) d\bullet\right)$, the Kloosterman stalk contributes a singularity of the form 
 $$\Phi(\xi) \psi^{-1}(\xi^{-1}) |\xi|^{s_2-\frac{1}{2}} d\xi$$ 
 at zero, which after the second convolution, by $\left( |\bullet|^{\frac{1}{2}-s_1} \psi(\bullet) d\bullet\right)$, gives rise to a singularity 
 $$\Phi_2(\xi) |\xi-1|^{s_1+s_2-1} d\xi$$ 
 in a neighborhood of $\xi=1$.
\end{proof}

I repeat the statement of Theorem \ref{maintheorem}, for the convenience of the reader.
Recall that $G^*$ is such that its dual group is $\check G_X$, that is: $G^*=\PGL_2$ when the spherical root of $X$ is of type $T$, and $G^*=\SL_2$ when the spherical root is of type $G$, and that we have defined in \S \ref{ssnotation} an enlarged space 
$\mathcal S^-_{L_X} (N_\psi\backslash G^*/N_\psi)$
of test measures for the Kuznetsov formula of $G^*$, determined by the $L$-value associated to $X$.

\begin{theorem}\label{maintheorem2}
Let $\C_X=(X\times X)\sslash G$. There is an isomorphism $\C_X \simeq \mathbbm A^1$, and the map $X\times X\to \mathbbm A^1$ is smooth away from the preimage of two points of $\mathbbm A^1$, that we will call singular. We fix the isomorphisms as follows:
\begin{itemize}
 \item When $\check G_X=\SL_2$, we take the set of singular points to be $\{0, 1\}$, with $X^\diag\subset X\times X$ mapping to $1 \in \C_X\simeq \mathbbm A^1$.
 \item When $\check G_X = \PGL_2$, we take the set of singular points to be $\{-2,2\}$, with $X^\diag\subset X\times X$ mapping to $2 \in \C_X\simeq \mathbbm A^1$.
\end{itemize}

Then, there is a continuous linear isomorphism: 
 \begin{equation}
  \mathcal T: \mathcal S^-_{L_X}(N_\psi\backslash G^*/N_\psi) \xrightarrow\sim \mathcal S(X\times X/G),
 \end{equation}
 given by the following formula:
 
\begin{itemize}
 \item When $\check G_X=\SL_2$ with $L_X=L(\Std, s_1) L(\Std,s_2)$ with $s_1\ge s_2$, 
 \begin{equation}
    \mathcal Tf(\xi) =  |\xi|^{s_1-\frac{1}{2}}  \left( |\bullet|^{\frac{1}{2}-s_1} \psi(\bullet) d\bullet\right) \star  \left( |\bullet|^{\frac{1}{2}-s_2} \psi(\bullet) d\bullet\right) \star f(\xi).
 \end{equation}
 \item When $\check G_X=\PGL_2$ with $L_X=L(\Ad, s_0)$,
 \begin{equation}\label{transferG}
    \mathcal Tf(\zeta) =  |\zeta|^{s_0-1}  \left( |\bullet|^{1-s_0} \psi(\bullet) d\bullet\right) \star  f(\zeta).
 \end{equation}
\end{itemize}

\end{theorem}

\begin{remark}\label{remarkconstants}
 The points $s_1, s_2, s_0$ are determined by the geometry, according to the following formulas:   
 \begin{equation}
  s_1+s_2=\frac{\dim X}{2};
 \end{equation}
 \begin{equation}
  s_1 = \frac{\left<\check\gamma, \rho_{P(X)}\right>}{2},
 \end{equation}
 and therefore, by \eqref{codimensionseq},
 \begin{equation}
  s_1-s_2= \frac{\dim X_{-1}}{2} -1;
 \end{equation}
 \begin{equation}\label{s0}
  s_0 = \left<\check\gamma, \rho_{P(X)}\right> = \frac{\dim X-1}{2} = \frac{\dim X_{-1}-1}{2},
 \end{equation}
 the last one by \eqref{equalcodim}.
 
 \end{remark}

\begin{proof}
 Let us start with the case $G^*=\SL_2$. Proposition \ref{fromtransfer} settles this case, but let me first observe, since similar arguments will be needed in the case of $G^*=\PGL_2$, that apart from the Kloosterman germ, the effect of Fourier convolutions can be calculated explicitly. Indeed, the space $\mathcal S^-_{L_X}(N_\psi\backslash G^*/N_\psi)$ can be thought of as the space of sections of a cosheaf over $\P^1(F)$; in a neighborhood of infinity, its elements have the form
 $$ f(\zeta) = |\zeta|^{1-s_0} \Phi(\zeta^{-1}) d^\times \zeta,$$
 where $\Phi$ is a smooth function. 
 If we consider the subspace $\mathcal S^-_{L_X}(N_\psi\backslash G^*/N_\psi)^0$ of Schwartz sections in the complement of $0\in \P^1$, it is immediate to see that the transfer operator of \eqref{transferG} defines a continuous isomorphism between this space and the space $\mathcal S(\mathbbm A^1)$ of usual Schwartz measures on the $F$-points of the affine line. Indeed, 
 \begin{eqnarray} \label{Fouriercalculation}  \mathcal Tf(\zeta) =  |\zeta|^{s_0-1}  \left( |\bullet|^{1-s_0} \psi(\bullet) d\bullet\right) \star  \left(|\zeta|^{1-s_0} \Phi(\zeta^{-1}) d^\times \zeta\right) = \nonumber \\
  = |\zeta|^{s_0-1} d^\times \zeta \cdot \int |z|^{1-s_0} \psi(z)  \left | \frac{\zeta}{z}\right|^{1-s_0} \Phi(\frac{z}{\zeta}) dz = d\zeta \cdot \int \psi(\zeta u) \Phi(u) du,
\end{eqnarray}
 which is the usual Fourier transform of a Schwartz function.  
 
 There remains to determine the behavior of the stalk at zero under this transform, i.e., the behavior of the ``Kloosterman germs''. 
 Applying Proposition \ref{fromtransfer} with $s_0=\frac{1}{2}-s$, we get that the operator $\mathcal T$ maps the space $\mathcal S^-_{L_X}(N_\psi\backslash G^*/N_\psi)$ isomorphically to $\mathcal S(\frac{\SL_2}{B_\ad, \delta_2^{1-s_0}})$, in the notation of that proposition. This is precisely the space of measures on the affine line $\mathbbm A^1$ (with coordinate $c$), which are smooth away from $c=\pm 2$, of rapid decay (together with their derivatives) at infinity, and in a neighborhood of $c=\pm 2$, setting $\xi=c\mp 2$, are of the form \eqref{twistedpfG}, with $\frac{d-3}{2}=s_0-1$. By Theorem \ref{Xtheorem2}, this is precisely the space $\mathcal S(X\times X/G)$.

 In the case $G^*=\PGL_2$, we start again with the subspace $\mathcal S^-_{L_X}(N_\psi\backslash G^*/N_\psi)^0$ of Schwartz sections in the complement of zero; those, now, are of the form 
 \begin{equation}\label{recall} 
 f(\xi)=  \left(\Phi_1(\xi^{-1}) |\xi|^{\frac{1}{2}-s_1} + \Phi_2(\xi^{-1}) |\xi|^{\frac{1}{2}-s_2}\right) d^\times \xi
 \end{equation}
 in a neighborhood of infinity,
 with the suitable logarithmic modification when $|\xi|^{s_1-s_2}$ is smooth. (All $\Phi_i$'s, here and below, denote smooth functions.) 
 
 The transfer operator $\mathcal T$ in this case is the composition of two multiplicative convolutions, followed by multiplication by the factor $|\xi|^{s_1-\frac{1}{2}}$. I claim that the first convolution, by the measure 
 $$\left( |\bullet|^{\frac{1}{2}-s_2} \psi(\bullet) d\bullet\right) = \left( |\bullet|^{\frac{3}{2}-s_2} \psi(\bullet) d^\times \bullet\right),$$ takes the space $\mathcal S^-_{L_X}(N_\psi\backslash G^*/N_\psi)^0$ to the space of measures which at infinity are of the form $\Phi_3(\xi^{-1}) |\xi|^{\frac{1}{2}-s_1} d^\times \xi$, while at zero are of the form $\Phi_4(\xi) |\xi|^{\frac{1}{2}-s_2} d\xi$, and otherwise smooth. Indeed, for the summand 
 $$ \Phi_2(\xi^{-1}) |\xi|^{\frac{1}{2}-s_2} d^\times \xi$$ 
 of \eqref{recall}, computing as in \eqref{Fouriercalculation} we see that it is mapped to a Schwartz measure times $|\xi|^{\frac{1}{2}-s_2}$. To determine the image of the germ represented by 
 $$ \Phi_1(\xi^{-1}) |\xi|^{\frac{1}{2}-s_1} d^\times \xi,$$
 we first notice that, by elementary properties of Fourier transform and a calculation like \eqref{Fouriercalculation}, its image under this convolution will be a smooth, but not necessarily compactly supported measure, times the factor $|\xi|^{\frac{1}{2}-s_2}$. Thus, we need to determine its behavior at infinity.
 
 In the non-Archimedean case, the stated behavior follows immediately from the functional equation of Tate integrals, which can be written as 
 $$ \left( |\bullet|^s \psi(\bullet) d^\times \bullet\right) \star \chi =  \gamma(\chi, 1-s,\psi) \cdot \chi$$
 (see \cite[2.12]{SaTransfer1}); this formula contains the usual Fourier transform of the expression $|\bullet|^{s-1} \chi^{-1}(\bullet)$, viewed as a generalized function by meromorphic continuation of the Tate zeta integrals, and makes sense away from the poles and zeroes of the gamma factor. Thus, if $f$ is a smooth measure, supported away from $0$, which at infinity is a multiple of $|\xi|^{\frac{1}{2}-s_1} d^\times \xi$, it will remain a multiple of this character after convolution, away from the poles or zeroes of the gamma factor $\gamma(s_2-s_1,\psi)$. In the non-Archimedean case, this gamma factor is zero when $s_1=s_2$, and infinite when $s_2 = s_1+1$, but here we have assumed that $s_1\ge s_2$. When $s_1=s_2$, instead of the asymptotics \eqref{Fouriercalculation}, we have a logarithmic term. To deal with that, and with the Archimedean case, one needs to consider the Mellin transforms of the measure $f$ before and after convolution, which live in a certain ``Paley--Wiener'' space of meromorphic functions, with poles determined by the asymptotics. I point the reader to \cite[Theorem 1.6]{Igusa} and \cite[\S 2.1]{SaTransfer1} for the relevant arguments.

 Finally, again by a calculation as in \eqref{Fouriercalculation}, the second convolution, by $\left( |\bullet|^{\frac{1}{2}-s_1} \psi(\bullet) d\bullet\right)$, takes this space to the space of measures which are of rapid decay (compactly supported, in the non-Archi\-medean case), of the form $\left(\Phi_5(\xi) |\xi|^{\frac{1}{2}-s_1} + \Phi_6(\xi) |\xi|^{\frac{1}{2}-s_2}\right) d\xi$ in a neighborhood of $\xi=0$, and otherwise smooth, and multiplication by the factor $|\xi|^{s_1-\frac{1}{2}}$ (by our convention that $s_1=\max(s_1,s_2)$) turns the germs at zero to 
 $$\left(\Phi_5(\xi)  + \Phi_6(\xi) |\xi|^{s_1-s_2}\right) d\xi =\left( \Phi_5(\xi)  + \Phi_6(\xi) |\xi|^{\frac{d_{-1}}{2}-1}\right) d\xi.$$

 There remains to examine the effect of the transfer operator to the ``Kloosterman germs'', i.e., to the stalk of $\mathcal S^-_{L_X}(N_\psi\backslash G^*/N_\psi)$ at zero. By Proposition \ref{frombeyond}, this stalk contributes an extra summand of 
 $$\Phi_7(\xi) |\xi-1|^{s_1+s_2-1} d\xi = \Phi_7(\xi) |\xi-1|^{\frac{d_1}{2}-1} d\xi$$
 in a neighborhood of $\xi=1$. Notice that multiplication by the factor $|\xi|^{s_1-\frac{1}{2}}$, which is smooth at $\xi=1$, does not alter this singularity. By Theorem \ref{Xtheorem2} and Proposition \ref{reductiontobasic}, this is the singularity of the elements of the space $\mathcal S(X\times X/G)$ at $\xi=1$. (Here, we apply the statement of Proposition \ref{reductiontobasic} about the coefficients $a_\eta$, instead of the more implicit expression \eqref{twistedpfT}.)

\end{proof}

\appendix

\section{Index of notation} \label{index}

I list some notation used through various sections. Notation used only locally is not included here.

\begin{description}
\item[$0\in\c_X^*$] Denotes the image of $0\in\a_X^*$. When an isomorphism of these spaces with $\mathbbm A^1$ is chosen, this point should map to $0$.
\item[$\lbrack\pm 1\rbrack$] The images of the points $\pm 1\in A_X$ in $A_X\sslash W_X$.
\item[$A$] The universal Cartan of $G$, i.e., the quotient of any Borel subgroup $B$ by its unipotent radical $N$.
\item[$A_X$] The universal Cartan of $X$, p.\ \pageref{refAX}.
\item[$\mathring\a_X^*$] The $W_X$-stable subset of $\a_X^*$ where $W_X$ acts freely. Equal to $\a_X^*\smallsetminus\{0\}$ in rank one.
\item[$B$] Stands for a Borel subgroup of $G$; typically, the choice of Borel subgroup does not matter for the statements, so we do not fix one.
\item[$\c_X^*$] The quotient $\a_X^*\sslash W_X$.
\item[$\mathring\c_X^*$] The quotient $\mathring\a_X^*/W_X$.
\item[$\C_X$] The quotient $X\times X\sslash G$ (diagonal action of $G$). After Proposition \ref{GITquotients}, this is also identified with the quotient $A_X\sslash W_X$ (in rank one).
\item[$\mathring\C_X$] The complement of $[\pm 1]\in A_X\sslash W_X = \C_X$.
\item[$dx, d^\times x$] The additive Haar measure on $F$ which is self-dual with respect to a fixed additive character $\psi$, and the multiplicative Haar measure on $F^\times$ given by $d^\times x = \frac{dx}{|x|}$. 
\item[$\delta$] Used for modular characters of various groups, defined as the quotient of the right by the left Haar measure.
\item[$\mathcal F(X)$] The space of Schwartz functions on the $F$-points of a smooth variety $X$, p.\ \pageref{refS}.
\item[$\hat\g^*$, $\tilde\g^*$] The polarization and Springer--Grothendieck resolution of $\g^*$, p.\ \pageref{refhatg}.
\item[$\g_X^*$] A normal cover of the image of the moment map, p.\ \pageref{refgX}.
\item[$\check G_X$] The dual group of the spherical variety $X$.
\item[$G^*$] The split group ($\PGL_2$ or $\SL_2$) whose dual group is equal to $\check G_X$, p.\ \pageref{refGstar}.
\item[$J$] Knop's abelian group scheme over $\c_X^*$, p.\ \pageref{refJ}.
\item[$J^0$] The open subsceme of identity components in the fibers of $J$, p.\ \pageref{refJ0}.
\item[$J_X$] A certain scheme that is birational to $J\bullet T^*X$, introduced on p.\ \pageref{refJX}.
\item[$\hat\kappa_X$, $\tilde\kappa_X$] Knop's sections, p.\ \pageref{refKsec}.
\item[$L_X$] The $L$-value attached to the spherical variety $X$, p.\ \pageref{refLX}.
\item[$L(X)$] The Levi quotient of the parabolic $P(X)$, p.\ \pageref{refPX}; also, a Levi subgroup, when the choice of Levi does not matter.
\item[$L_1$] The kernel of the map $L(X)\to A_X$, p.\ \pageref{refL1}.
\item[$\Lambda_X$] The character group of $A_X$, $\Lambda_X=\Hom(A_X,\Gm)$; analogous notation used for character groups of Borel orbits, p.\ \pageref{refLambdaY}. 
\item[$M_1\bullet M_2$] For schemes over $\c_X^*$, denotes $M_1\times_{\c_X^*} M_2$.
\item[$M^\rs$] For a scheme over $\a_X^*$ or $\c_X^*$, the preimage of $\mathring\a_X^*$, resp.\ \ $\mathring\c_X^*$, p.\ \pageref{refrs}.
\item[$\mu$] The moment map $T^*X\to \g^*$.
\item[$\mu_\inv$] The invariant moment map $T^*X\to \c_X^*$, p.\ \pageref{refmuinv}.
\item[$N$] The unipotent radical of a Borel subgroup $B$ of $G$; or, the upper triangular unipotent subgroup of $G^*=\PGL_2, \SL_2$, identified with the additive group $\Ga$.
\item[$\mathcal N(H)$] The normalizer of a subgroup $H$.
\item[$N_Y X, N_Y^* X$] The normal, conormal bundle of $Y$ in $X$.
\item[$N_Y, N^*_Y$] Specifically for the ambient variety $X\times X$, this is shorthand for $N_Y(X\times X)$, $N_Y^*(X\times X)$.
\item[$P_S$] For a set $S$ of simple roots, and having fixed a Borel subgroup $B$, the parabolic $P_S\supset B$ generated by the root spaces of the roots $-\alpha$, $\alpha\in S$.
\item[$P(X)$] The class of parabolics stabilizing the open Borel orbit, p.\ \pageref{refPX}.
\item[$\mathcal R(G)$] The radical of a group $G$.
\item[$\mathcal S(X)$] The space of Schwartz measures on the $F$-points of a smooth variety $X$, p.\ \pageref{refS}.
\item[$\mathcal S(X/G)$] The pushforward of $\mathcal S(X)$ under $X\to X\sslash G$, p.\ \pageref{refSmodG}.
\item[$\mathcal S(N_\psi\backslash G^*/N_\psi)$] The twisted pushforward for the Kuznetsov quotient, p.\ \pageref{refSKuz}.
\item[$\mathcal S^-_{L_X} (N_\psi\backslash G^*/N_\psi)$] The enlarged space of Kuznetsov test measures determined by the $L$-value $L_X$, p.\ \pageref{refenlarged}.
\item[$\mathcal S_Y(\bullet), \overline{\mathcal S_Y(\bullet)}$] The stalks and fibers of certain cosheaves of Schwartz measures over a closed subspace $Y$, p.\ \pageref{refstalkfiber}.
\item[$T_Y X, T^*_Y X$] The restriction of the tangent, cotangent bundle of $X$ to the subspace $Y$.
\item[$\widehat{T^*X}$, $\widetilde{T^*X}$] The polarization and Springer--Grothendieck lift of $T^*X$, p.\ \pageref{refhatTX}.
\item[$\widehat{T^*X}^\bullet$, $\widetilde{T^*X}^\bullet$] The distinguished irreducible components of $\widehat{T^*X}$, $\widetilde{T^*X}$, p.\ \pageref{refTXbullet}.
\item[$U_P$] The unipotent radical of a parabolic $P$.
\item[$W_X$] The Weyl group of the spherical variety $X$, p.\ \pageref{refAX}.
\item[$\mathring X$] The open Borel orbit on $X$ (for some choice of Borel subgroup).
\item[$X_\emptyset, X_\emptyset^\bullet$] The affine boundary degeneration of $X$, and its open $G$-orbit, p.\ \pageref{refboundary}.
\item[$X\sslash G$] The invariant-theoretic quotient $\Spec F[X]^G$.
\item[$\psi$] A fixed, non-trivial additive character $F\to \CC^\times$; also identified with a character of the upper triangular unipotent subgroup $N\subset G^*=\PGL_2$ or $\SL_2$.
\item[$\mathcal Z(G)$] The center of a group $G$.
\end{description}

\bibliographystyle{alphaurl}
\bibliography{biblio}

\end{document}